\documentclass[a4paper]{article}
\pdfoutput=0 

\usepackage{lipsum}
\usepackage{amsfonts}
\usepackage{graphicx}
\usepackage{epstopdf}
\usepackage{algorithmic}
\usepackage{xcolor}
\usepackage{amsmath, amssymb, theorem, latexsym}
\usepackage{rotating}
\usepackage{multirow}
\usepackage{url}
\usepackage{enumerate}

\ifpdf
  \DeclareGraphicsExtensions{.eps,.pdf,.png,.jpg}
\else
  \DeclareGraphicsExtensions{.eps}
\fi

\newcommand{\TheTitle}{On spectral properties of the Bloch-Torrey operator in two dimensions}

\newcommand{\TheAuthors}{D. S. Grebenkov and B. Helffer}

\title{{\TheTitle}} 

\author{Denis S. Grebenkov\thanks{Laboratoire de Physique de la Mati\`ere Condens\'ee, 
 CNRS--Ecole Polytechnique, University Paris-Saclay, 91128 Palaiseau, France 
 (denis.grebenkov@polytechnique.edu).}
 \and
 Bernard Helffer\thanks{Laboratoire de Math\'ematiques Jean Leray, Universit\'e de Nantes
 2 rue de la Houssini\`ere 44322 Nantes, France,
 and
 Laboratoire de Math\'ematiques, Universit\'e Paris-Sud, CNRS, Univ. Paris Saclay, France
 (bernard.helffer@univ-nantes.fr).}
}

\usepackage{amsopn}

\ifpdf
\hypersetup{
  pdftitle={\TheTitle},
  pdfauthor={\TheAuthors}
}
\fi

\def\B{{\mathcal B}}
\def\R{{\mathbb R}}

\def\Kappa{{\mathcal K}}

\def\curv{{\mathfrak c}}

\def\C{{\mathbb C}}

\def\t{{\rho}}
\def\A{{\mathcal A}}
\def\L{{\mathcal L}}
\def\Ai{\textrm{Ai}}
\def\Re{\textrm{Re}}
\def\Im{\textrm{Im}}
\newtheorem{theorem}{Theorem}[section]
\newtheorem{lemma}[theorem]{Lemma}

\newtheorem{proposition}[theorem]{Proposition}

\newtheorem{remark}[theorem]{Remark}

\newtheorem{example}[theorem]{Example}

\newtheorem{assumption}[theorem] {Assumption}
\newcommand{\pa}{\partial}
\newcommand{\sign}{{\rm sign\,}}

\newcommand*{\qed}{\hfill\ensuremath{\square}}%

\makeatletter
\@addtoreset{equation}{section}
\renewcommand{\theequation}{\thesection.\arabic{equation}}
\makeatother

\begin{document}

\maketitle

\begin{abstract}
We investigate a two-dimensional Schr\"odinger operator, $-h^2 \Delta
+iV(x)$, with a purely complex potential $iV(x)$.  A rigorous
definition of this non-selfadjoint operator is provided for bounded
and unbounded domains with common boundary conditions (Dirichlet,
Neumann, Robin and transmission).  We propose a general perturbative
approach to construct its quasimodes in the semi-classical limit.  An
alternative WKB construction is also discussed.  These approaches are
local and thus valid for both bounded and unbounded domains, allowing
one to compute the approximate eigenvalues to any order in the small
$h$ limit.  The general results are further illustrated on the
particular case of the Bloch-Torrey operator, $-h^2\Delta + ix_1$, for
which a four-term asymptotics is explicitly computed.  Its high
accuracy is confirmed by a numerical computation of the eigenvalues
and eigenfunctions of this operator for a disk and circular annuli.
The localization of eigenfunctions near the specific boundary points
is revealed.  Some applications in the field of diffusion nuclear
magnetic resonance are discussed.
\end{abstract}

{\bf Keywords}: 
Transmission boundary condition, spectral theory, Bloch-Torrey equation, semi-classical analysis, WKB

{\bf AMS}:
 35P10,    
 47A10,    
 47A75     

\section{Introduction}

In a previous paper \cite{GHH}, we have analyzed in collaboration with
R. Henry one-dimensional models associated with the complex Airy
operator $-\frac{d^2}{dx^2} + i gx$ on the line, with $g\in \R\,$.  We
revisited the Dirichlet and Neumann realization of this operator in
$\mathbb R^+$ and the main novelty was to consider a transmission
problem at $0\,$.  In higher dimensions, an extension of the complex
Airy operator is the differential operator that we call the
Bloch-Torrey operator or simply the BT-operator
\begin{equation*} 
- D\Delta + i g x_1\,,
\end{equation*}
where $\Delta = \partial^2/\partial x_1^2 + \ldots +
\partial^2/\partial x_n^2$ is the Laplace operator in $\R^n$, 
and $D$ and $g$ are real parameters.  More generally, we will study
the spectral properties of some realizations of the differential
Schr\"odinger operator
\begin{equation}
\A_h^\# = -h^2 \Delta + i\,V(x)\, ,
\end{equation}
in an open set $\Omega$, where $h$ is a real parameter and $V(x)$ a
real-valued potential with controlled behavior at $\infty$, and the
superscript $\#$ distinguishes Dirichlet (D), Neumann (N), Robin (R),
or transmission (T) conditions.  More precisely we discuss

\begin{enumerate}
\item 
the case of a bounded open set $\Omega$ with Dirichlet, Neumann or
Robin boundary condition;
\item 
the case of a complement $\Omega:=\complement \overline{\Omega}_-$ of
a bounded set $ \Omega_-$ with Dirichlet, Neumann or Robin boundary
condition;
\item
the case of two components $ \Omega_-\cup\Omega_+$, with $\Omega_-
\subset \overline{\Omega}_- \subset \Omega$ and $\Omega_+ = \Omega
\backslash \overline{\Omega}_-$, with $\Omega$ bounded and
transmission conditions at the interface between $\Omega_-$ and
$\Omega_+$;
\item
the case of two components $ \Omega_-\cup\, \complement\,
\overline{\Omega}_-\,$, with $ \Omega_-$ bounded and transmission
conditions at the boundary;
\item 
the case of two unbounded components $\Omega_-$ and $\Omega_+$
separated by a hypersurface with transmission conditions.
\end{enumerate}
In all cases, we assume that the boundary is $C^\infty$ to avoid
technical difficulties related to irregular boundaries (see
\cite{Grisvard}).  
Roughly speaking (see the next section for a precise definition), the
state $u$ (in the first two items) or the pair $(u_-,u_+)$ in the last
items should satisfy some boundary or transmission condition at the
interface.  In this paper, we consider the following situations:
\begin{itemize}
\item 
the Dirichlet condition: $u_{|\partial \Omega} =0\,$;
\item 
the Neumann condition: $\partial_\nu u_{|\partial \Omega} =0\,$,
where $\partial _\nu = \nu\cdot \nabla$, with $\nu$ being the outwards
pointing normal;
\item 
the Robin condition: $h^2 \partial_\nu u_{|\partial \Omega} = -
\Kappa u_{|\partial \Omega}\,$, where $\Kappa \geq 0$ denotes the
Robin parameter;
\item 
the transmission condition: 
$$h^2 \partial_\nu u_{+\, |\partial \Omega_-} = h^2 \partial_\nu u_{-\,|\partial \Omega_-} 
= \Kappa (u_{+\,|\partial \Omega_- }-u_{-\, |\partial \Omega_- })\,, $$
where $\Kappa \geq 0$ denotes the transmission parameter, and the
normal $\nu$ is directed outwards $\Omega_-$.
\end{itemize}
From now on $\Omega^\#$ denotes $\Omega$ if $\# \in \{D,N,R\}$ and
$\Omega_-$ if $\#=T\,$. $L^{2}_{\#}$ will denote $L^2(\Omega)$ if $\#
\in \{D,N,R\}$ and $L^2(\Omega_-) \times L^2(\Omega_+)$ if $\#=T\,$.\\

In \cite{GHH}, we have analyzed in detail various realizations of the
complex Airy (or Bloch-Torrey) operator $A_0^\#:= -
\frac{d^2}{d\tau^2} + i \tau$ in the four cases corresponding to
Dirichlet, Neumann, and Robin on the half-line $\mathbb R^+$ or for
the transmission problem on the whole line $\mathbb R$ (in what
follows, $\mathbb R^\#$ will denote $\mathbb R^+$ if $\#\in \{D,N,R\}$
and $\mathbb R$ if $\#=T$).  The boundary conditions read
respectively:
\begin{itemize}
\item $u(0)=0\,$;
\item $u'(0)=0\,$;
\item $u'(0) = \kappa\, u (0)\,$;
\item $u'_-(0) =u'_+(0) = \kappa \, (u_+(0) - u_-(0))$
\end{itemize}
(with $\kappa \geq 0$ in the last items).  For all these cases, we
have proven the existence of a discrete spectrum and the completeness
of the corresponding generalized eigenfunctions.  Moreover, there is
no Jordan block (for the fourth case, this statement was proven only
for $\kappa$ small enough). \\

In this article, we start the analysis of the spectral properties of
the BT operator in dimensions 2 or higher that are relevant for
applications in superconductivity theory
\cite{Alm,aletal10,aletal11,aletal12}, in fluid dynamics \cite{Mar},
in control theory \cite{BHHR}, and in diffusion magnetic resonance
imaging \cite{Grebenkov07,Gr1} (and references therein).  We will
mainly focus on
\begin{itemize}
\item definition of the operator, 
\item construction of approximate eigenvalues in some asymptotic regimes,
\item  localization of quasimode states near certain boundary points,
\item numerical simulations.
\end{itemize}
In particular, we will discuss the semi-classical asymptotics $h\to
0\,$, the large domain limit, the asymptotics when $g\rightarrow 0$ or
$+\infty\,$, the asymptotics when the transmission or Robin parameter
tends to $0\,$.  Some other important questions remain unsolved like
the existence of eigenvalues close to the approximate eigenvalues (a
problem which is only solved in particular situations).  We hope to
contribute to this point in the future.\\

When $g = 0$, the BT-operator is reduced to the Laplace operator for
which the answers are well known.  In particular, the spectrum is
discrete in the case of bounded domains and equals $[0,+\infty)$ when
one or both components are unbounded.  In the case $g \neq 0\,$, we
show that if there is at least one boundary point at which the normal
vector to the boundary is parallel to the coordinate $x_1$, then there
exist approximate eigenvalues of the BT-operator {\it suggesting} the
existence of eigenvalues while the associated eigenfunctions are
localized near this point.  This localization property has been
already discussed in physics literature for bounded domains
\cite{deSwiet94}, for which the existence of eigenvalues is trivial.
Since our asymptotic constructions are local and thus hold for
unbounded domains, the localization behavior can be conjectured for
exterior problems involving the BT-operator.

Some of these questions have been already analyzed by Y. Almog (see
\cite{Alm} and references therein for earlier contributions), R. Henry
\cite{Hen1, Hen2} and Almog-Henry \cite{AlHe} but they were mainly
devoted to the case of a Dirichlet realization in bounded domains in
$\mathbb R^2$ or particular unbounded domains like $\mathbb R^2$ and
$\mathbb R^2_+$, these two last cases playing an important role in the
local analysis of the global problem.\\ 
Different realizations of the operator $\mathcal A_h$ in $\Omega$ are
denoted by $\A_h^D\,$, $\A_h^N\,,$ $\A_h^R\,$ and $\A_h^T$.  These
realizations will be properly defined in Section
\ref{sec:definition} under the condition that, when $\Omega$ is
unbounded, there exists $C > 0$ such that
\begin{equation}\label{condsup}
|\nabla V (x) | \leq C \sqrt{1 + V(x)^2}\,.
\end{equation}
Our main construction is local and summarized in the following
\begin{theorem}\label{th:app}
Let $\Omega\subset \mathbb R^2$ as above,
$V\in\mathcal{C}^\infty(\overline{\Omega} ; \mathbb{R})$ and $x^0\in
\partial \Omega^\#$ such that%
\footnote{
As noticed in \cite{AlHe}, a point satisfying the second condition in
\eqref{defPaPerp} always exists when $\partial \Omega^\#$ is bounded.}
\begin{equation}\label{defPaPerp}
\nabla V (x^0) \neq 0\,,\quad  \nabla V(x^0) \wedge \nu(x^0) = 0 \,,
\end{equation}
where $\nu (x^0)$ denotes the outward normal on $\partial\Omega$ at
$x^0$\,.\\
Let us also assume that, in the local curvilinear coordinates, the
second derivative of the restriction of $V$ to the boundary at $x^0$
(denoted as $2\, v_{20}$) satisfies
\begin{equation*}
v_{20} \neq 0\,.
\end{equation*}
For the Robin and transmission cases, we also assume that for some
$\kappa >0$
\begin{equation}\label{scaKappa}
\Kappa =  h^{\frac 43} \kappa\,.
\end{equation}
If $\mu_0^\#$ is a simple eigenvalue of the realization ``$\#$'' of
the complex Airy operator $-\frac{d^2}{dx^2} + ix$ in $L^{2}_{\#}\,$,
and $\mu_2$ is an eigenvalue of the Davies operator
\break $-\frac{d^2}{dy^2} + i y^2$ on $L^2(\mathbb R)$, then there exists
an approximate pair $(\lambda_h^\#,u_h^\#)$ with $u_h^\#$ in the
domain of $\mathcal A_h^\#$, such that
\begin{equation}
\label{eq:lambdah_aux} 
\lambda_h ^\#= i \, V(x^0)  + h^{\frac 23} \sum_{j\in \mathbb N} \lambda_{2j}^\#\,  h^{\frac j 3}  +  \mathcal O (h^\infty)\,,
\end{equation}
\begin{equation} 
(\A_h^\# -\lambda_h^\#)\, u_h^\# = \mathcal O(h^\infty) \, \mbox{ in } L^{2}_{\#}(\Omega)\,,\quad ||u_h^\#||_{L^2} \sim  1\,,
\end{equation}
where
\begin{equation}
\lambda_0^\# = \mu_0^\# \,|\, v_{01}|^\frac 23  \exp \left(i\frac{\pi}{3} \sign \, v_{01}\right)  \,,\quad 
\lambda_2 =\mu_2  |v_{20}|^\frac 12 \exp \left(i \frac{\pi}{4} \sign v_{20} \right)\,,
\end{equation}
with $v_{01}:= \nu \cdot \nabla V (x^0)\,$.
\end{theorem}
In addition, we will compute $\lambda_4^\#$ explicitly (see the
Appendix) in the four types of boundary conditions and also describe
an alternative WKB construction to have a better understanding of the
structure of the presumably corresponding eigenfunctions.  We will
also discuss a physically interesting case when $\kappa $ in
\eqref{scaKappa} depends on $h$ and tends to $0\,$.

The proof of this theorem provides a general scheme for quasimode
construction in an arbitrary planar domain with smooth boundary
$\partial\Omega$.  In particular, this construction allows us to
retrieve and further generalize the asymptotic expansion of
eigenvalues obtained by de Swiet and Sen for the Bloch-Torrey operator
in the case of a disk \cite{deSwiet94}.  The generalization is
applicable for any smooth boundary, with Neumann, Dirichlet, Robin, or
transmission boundary condition.  Moreover, since the analysis is
local, the construction is applicable to both bounded and unbounded
components.

The paper is organized as follows.  In Sec. \ref{sec:definition}, we
provide rigorous definitions and basic properties of the BT-operator
in bounded and unbounded domains, with Dirichlet, Neumann, Robin, and
Transmission conditions.  Sec. \ref{sec:former} recalls former
semi-classical results for a general operator $-h^2 \Delta + iV(x)$.
In Sec.~\ref{sec:quasimodesP}, we provide preliminaries for
semi-classical quasimode constructions in the two-dimensional case.
The construction scheme is detailed in Sec.~\ref{sec:quasimodes}.  In
particular, the four-terms asymptotics of the approximate eigenvalues
is obtained and we prove the main theorem.  In
Sec.~\ref{sec:kappa_scaling} we consider other scaling regimes for the
Robin or transmission parameter.  In Sec.~\ref{sec.WKB} we propose an
alternative construction for the first approximate eigenvalue using
WKB quasi-mode states.  In Sec.~\ref{sec:examples}, we illustrate
general results for simple domains such as disk and annulus.
Sec.~\ref{sec:numerics} describes numerical results in order to check
the accuracy of the derived four-terms asymptotics of eigenvalues of
the BT-operator in simple domains such as a disk, an annulus, and the
union of disk and annulus with transmission boundary condition.  We
also illustrate the localization of eigenfunctions near circular
boundaries of these domains.  Since a direct numerical computation for
unbounded domains (e.g., an exterior of the disk) is not possible, we
approach this problem by considering an annulus with a fixed inner
circle and a moving away outer circle.  We check that the localization
of some eigenfunctions near the inner circle makes them independent of
the outer circle.  We therefore conjecture that the BT-operator has
some discrete spectrum for the exterior of the disk.  More generally,
this property is conjectured to hold for any domain in $\R^n$ (bounded
or not) with smooth boundary which has points whose normal is parallel
to the gradient direction.  Finally, we briefly discuss in
Sec. \ref{sec:NMR} the importance of the obtained results in the field
of diffusion nuclear magnetic resonance (see \cite{Grebenkov17} for
further details).

{ \bf Acknowledgments.\\} We thank Raphael Henry who collaborated with
us in \cite{GHH} and in the preliminary discussions for the present
paper.  The second author would also like to thank Yaniv Almog and
Didier Robert for useful discussions.

\section{Definition of the various realizations of the Bloch-Torrey operator}
\label{sec:definition}

\subsection{The case of a bounded open set $\Omega$}

This is the simplest case.  For the analysis of the Dirichlet (resp.
Neumann) realization $\A_h^D$ (resp. $\A_h^N$) of the BT-operator, the
term $V(x)$ is simply a bounded non self-adjoint perturbation of the
Dirichlet (resp. Neumann) Laplacian. \\ We have for three boundary
conditions:
\begin{itemize}
\item 
For the Neumann case, the form domain $\mathcal V$ is $H^1(\Omega)$ and
(if $\Omega$ is regular) the domain of the operator is $\{ u\in
H^2(\Omega)\,,\, \partial_\nu u_{/\partial \Omega} =0\}$.  The
quadratic form reads
\begin{equation}
\label{eq:formDN}
\mathcal V\ni u \mapsto  q_V(u) := h^2\,  || \nabla u||^2_{ \Omega} 
  + i\, \int_{ \Omega} V(x)\,  | u(x) |^2 \, dx \,.
\end{equation} 

\item
For the Dirichlet case, the form domain is $H_0^1(\Omega)$ and (if
$\Omega$ is regular) the domain of the operator is
$H^2(\Omega)\mathcal \cap H_0^1(\Omega)$.  The quadratic form is given
by (\ref{eq:formDN}).

\item
For the Robin case (which is a generalization of the Neumann case),
the form domain is $H^1(\Omega)$ and (if $\Omega$ is regular) the
domain of the operator $\A_h^R$ is $\{ u\in H^2(\Omega)\,,\, - h^2
\partial_\nu u_{/\partial \Omega} = \Kappa u_{/\partial \Omega} \}\,$,
where $\Kappa$ denotes the Robin coefficient, and $\nu$ is pointing
outwards.  The quadratic form reads
\begin{equation}
\label{eq:formR}
u \mapsto  q_V(u) := h^2 \,|| \nabla u||^2_{ \Omega} 
  + i \int_{ \Omega} V(x)  | u(x) |^2 \, dx + \Kappa \int_{\partial \Omega} |u|^2 ds \,.
\end{equation}
The Neumann case is retrieved for $\Kappa=0\,$.
\end{itemize}

\subsection{The case of a bounded set in $\mathbb R^n$ and its complementary set with transmission condition at the boundary}

We consider $ \Omega_-\cup\,\complement\, \overline{\Omega}_-\,$, with
$ \Omega_-$ bounded in $\mathbb R^n$ and $\partial \Omega_-$
connected.  In this case the definition of the operator is similar to
what was done for the one-dimensional case in \cite{GHH}.  However, we
start with a simpler case when $ \Omega_- \subset \overline{ \Omega}_-
\subset \Omega$ with $\Omega$ bounded and $ \Omega_+ = \Omega
\setminus \overline{\Omega}_-$ (with Neumann boundary condition
imposed on the exterior boundary $\partial\Omega$).  After that, we
explain how to treat the unbounded case with $\Omega = \R^n$ and
$\Omega_+ = \complement\, \overline{ \Omega}_-\,$.  Note that the case
of a complement $\Omega:=\complement \overline{\Omega}_-$ of a bounded
set $ \Omega_-$ with Dirichlet, Neumann or Robin boundary condition
can be treated along the same lines, the transmission problem being
the most complicated case.

\subsubsection{Transmission property in the bounded case} 
\label{sec:bounded_trans}

To treat the difficulties one by one, we start with the situation when
$ \Omega_- \subset \overline{\Omega}_- \subset \Omega\,$, $ \Omega_+ :
= \Omega \setminus \overline{\Omega}_-\,$, and $\Omega$ bounded and
connected (e.g., a disk inside a larger disk).\\ We first introduce
the variational problem, with the Hilbert space
\begin{equation*}
\mathcal H = L^2( \Omega_-)\times L^2( \Omega_+)
\end{equation*}
and the form domain
\begin{equation*}
\mathcal V := H^1( \Omega_-) \times H^1( \Omega_+)\,.
\end{equation*}
The quadratic form reads on $\mathcal V$
\begin{equation}
\begin{split}
u=(u_-,u_+) \mapsto  q_V(u) &:= h^2 || \nabla u_-||^2_{ \Omega_-} + h^2 || \nabla u_+||^2_{ \Omega_+} 
 + \,\Kappa \, || (u_--u_+)||^2_{L^2(\partial  \Omega_-)} \\
& \quad   + i \int_{ \Omega_-} V(x)  | u_-(x) |^2 \, dx +  i \int_{ \Omega_+} V(x)  | u_+(x)|^2\, dx\, \, , \\
\end{split}
\end{equation}
where $\Kappa$ is a positive parameter of the transmission problem,
and $h>0$ is a semi-classical parameter whose role will be explained
later and which can be thought of as equal to one in this section.
The dependence of $\Kappa$ on $h>0$ will be discussed later.  We
denote by $\mathfrak a_V$ the associated sesquilinear form:
\begin{equation*}
\mathfrak a_V(u,u)=q_V (u)\,.
\end{equation*}
The potential $V(x)$ is assumed to be real (and we are particularly
interested in the example $V(x) = g x_1$).  In this case, one gets
continuity and coercivity of the associated sesquilinear form on
$\mathcal V$ (after a shift of the quadratic form by adding a
constant).  This is true for any $\Kappa$ without assumption on its
sign.  The trace of $u_-$ and $u_+$ on $\partial \Omega_-$ is indeed
well defined for $(u_-,u_+)\in \mathcal V$. \\

Applying Lax-Milgram's theorem to the shifted form, we first get that
the solution of the variational problem associated with $\mathfrak
a_V$, $(u_-,u_+)$, should satisfy $\Delta u_-\in L^2(
\Omega_-)$ and $\Delta u_+ \in L^2(\Omega_+)$.  Together with
$(u_-,u_+)\in \mathcal V$ this permits to define the Neumann condition
(via the Green formula) for both $u_-$ and $u_+$ in $H^{-\frac 12}
(\partial \Omega_-)$, and in addition for $u_+$ in $H^{-\frac 12} (\pa
\Omega)$.  Indeed, to define $\partial_\nu u_-$ as a linear form on
$H^\frac 12 (\partial \Omega_-)$, we use that for any $v\in
H^1(\Omega_-)$,
\begin{equation}\label{Greenfor}
- \int_{ \Omega_-}\Delta u_- \, v \, dx = \int_{ \Omega_-} \nabla u_-\cdot \nabla v \, dx + \int_{\partial  \Omega_-}  \partial_\nu u_- \, v  \, d\sigma\,,
\end{equation}
and the existence of a continuous right inverse for the trace from
$H^\frac 12 (\partial \Omega_-)$ into $H^1( \Omega_-)\,$.  Here the
normal $\nu$ is oriented outwards $\Omega_-$ and when $u_-$ is more
regular ($u_-\in H^2(\Omega_-)$),  we have $\partial_\nu u_- =
\nu\cdot \nabla u_-$.  In a second step we get the Neumann condition
for $u_+$ on $\partial \Omega$,
\begin{equation}\label{condNeu1}
\pa_\nu u_+ =0\mbox{ on } \pa \Omega\,,
\end{equation}
and the transmission condition on $\partial  \Omega_-$
\begin{equation} \label{Condter}
\begin{array}{l l}
\begin{split}
\pa_\nu u_- & = \pa_\nu u_+  \\
 h^2 \pa_\nu u_- & =  \Kappa \, \bigl(u_+ - u_- \bigr)\,  \\
\end{split}  &
\mbox{ on } \partial  \Omega_-\,,  \\
\end{array}
\end{equation}
which is satisfied in $H^{-\frac 12}(\pa \Omega_-)$.  We keep here the
previous convention about the outwards direction of $\nu$ on $\partial
\Omega_-$.
\\
Finally, we observe that the first traces of $u_-$ and $u_+$ on $\pa
 \Omega_-$ belong to $H^\frac 12 (\pa  \Omega_-)$.  Hence by
\eqref{Condter}, the second traces of $u_-$ and $u_+$ are
in $H^\frac 12(\pa \Omega_-)$.  But now the regularity of the Neumann
problem in $ \Omega_-$ and $ \Omega_+$ implies that
\begin{equation*}
(u_-,u_+) \in H^2( \Omega_-) \times H^2( \Omega_+)\,.
\end{equation*}
Here we have assumed that all the boundaries are regular. 

\begin{remark}
One can actually consider a more general problem in which the two
diffusion coefficients $D_-$ and $D_+$ in $\Omega_-$ and $\Omega_+$
are different.  The transmission condition reads
\begin{equation*}
D_{-} \partial_\nu u_{-} = D_{+} \partial_\nu u_{+} = \Kappa  (u_{+} - u_{-})\quad  \mbox{ on } \partial \Omega_-\,.
\end{equation*}
If we take $D_{-} = D_{+} = D = h^2$, we recover the preceding case.
In the limit $D_{+} \to \infty$, we can consider the particular case
where $u_{+}$ is identically $0$ and we recover the Robin condition on
the boundary $\partial \Omega_-$ of the domain $\Omega_-$.
\end{remark}

\subsubsection{The unbounded case with bounded transmission boundary}
\label{sec:Def_unbouded}

In the case $ \Omega_+ = \complement \overline{\Omega}_-$ (i.e.,
$\Omega =\mathbb R^n$), we have to treat the transmission problem
through $\partial \Omega_-$ with the operator $-h^2 \Delta + i V(x)$
on $L^2( \Omega_-)\times L^2(\Omega_+)$.  Nothing changes at the level
of the transmission property because $\partial \Omega_-$ is bounded.
However, the variational space has to be changed in order to get the
continuity of the sesquilinear form.  Here we have to account for the
unboundedness of $V$ in $ \Omega_+$.  For this purpose, we introduce
\begin{equation}
\mathcal V :=\{ (u_-,u_+)\in \mathcal H \,,\, |V|^\frac 12 u_+ \in L^2( \Omega_+)\}\,.
\end{equation}
If $V$ has constant sign outside a compact, there is no problem to get
the coercivity by looking separately at $\Re\, \mathfrak a_V(u,u)$ and $\Im\,
\mathfrak a_V (u,u)$.  When $V$ does not have this property (as it is in the
case $V(x)=x_1$), one cannot apply Lax-Milgram's theorem in its
standard form.  We will instead use the generalized Lax-Milgram
Theorem as presented in \cite{AH} (see also \cite{GHH}).

\begin{theorem} \label{LMn}
Let $\mathcal V$ denote a Hilbert space and let $\mathfrak a$ be a continuous
sesquilinear form on $\mathcal V\times \mathcal V$.  If $\mathfrak a$ satisfies,
for some $\Phi_1, \Phi_2 \in \mathcal L (\mathcal V)\,,$ and some $\alpha >0$, 
\begin{equation}
\label{lm5n}
|\mathfrak a(u,u)| + |\mathfrak a(u, \Phi_1(u))| \geq \alpha\, \|u\|_\mathcal V^2\,,\quad \forall u\in \mathcal V\,,
\end{equation}
\begin{equation}
\label{lm5ne}
|\mathfrak a(u,u)| + |\mathfrak a( \Phi_2(u),u)| \geq \alpha\, \|u\|_\mathcal V^2\,,\quad \forall u\in \mathcal V\,,
\end{equation}
then $ A\in \mathcal L(\mathcal V)$ defined by
\begin{equation}\label{lm2}
\mathfrak a (u,v) = \langle Au\,,\,v\rangle_\mathcal V\,,\, \forall u \in \mathcal V\,,\quad \forall v \in \mathcal V\,,
\end{equation}
is a continuous isomorphism from $\mathcal V$ onto $\mathcal V$.
\end{theorem}
We now consider two Hilbert spaces $\mathcal V$ and $\mathcal H$ such
that $\mathcal V\subset \mathcal H$ (with continuous injection and
dense image).  Let $\A$ be defined by
\begin{equation}\label{lm3}
D(\A) =\{u\in \mathcal V\;|\; v\mapsto \mathfrak a (u,v) \mbox{ is continuous on } \mathcal V \mbox{ in the
  norm of } \mathcal H\}
\end{equation}
and
\begin{equation}
\label{lm4}
\mathfrak a (u,v) = \langle \A u\,,\,v \rangle_\mathcal H \quad \forall u \in D(\A) \mbox{ and } \forall v\in \mathcal V\,.
\end{equation}
Then we have 
\begin{theorem}
\label{LaxMilgramv2}
Let $\mathfrak a$ be a continuous sesquilinear form satisfying
\eqref{lm5n} and \eqref{lm5ne}.  Assume further that $\Phi_1$ and
$\Phi_2$ extend into continuous linear maps in $\mathcal L (\mathcal
H)\,$.  Let $\A$ be defined by \eqref{lm3}-\eqref{lm4}.  Then
\begin{enumerate}
\item $\A$ is bijective from $D(\A)$ onto $\mathcal H\,$.
\item  $D(\A)$ is dense in both $\mathcal V$ and $\mathcal H$\,.
\item $\A$ is closed.
\end{enumerate}
\end{theorem}

\begin{example}
We can use on $\mathcal V\,$ the multipliers
\begin{equation*}
\Phi_1(u_-,u_+) = \Phi_2(u_-,u_+) = \left(u_-\,,\, \frac{V}{\sqrt{1+ V^2}}~ u_+\right).
\end{equation*}
We first observe that, for some $C>0$,
\begin{equation*}
\Re \, \mathfrak a_V (u,u) \geq \frac 1 C  \left(||\nabla u_-||^2 + || \nabla u_+||^2\right) - C \left(  || u_-||^2 + || u_+||^2\right)\,.
\end{equation*}
To obtain the generalized coercivity (after a shift of the quadratic
form), we now look at $\Im \, \mathfrak a_V (u, \Phi_1(u))$ and get,
for some $\hat C>0$,
\begin{equation*}
\Im \, \mathfrak a_V (u, \Phi_1(u)) \geq  \int_{ \Omega_+} |V(x)| |u_+|^2\, dx  - \hat C\, ( ||u||^2 + || \nabla u||^2)\,.
\end{equation*}
Note that this works (see \cite{AH}) for general potentials $V(x)$
satisfying \eqref{condsup}.
\end{example}

Note also that the domain of the operator $\A^T$ associated with the
sesquilinear form is described as follows
\begin{equation}
\begin{array}{ll}
D& := \{ u \in \mathcal V \,,\, (- h^2\Delta + i V) u_- \in L^2( \Omega_-)\,,\, (- h^2 \Delta + i V) u_+ \in L^2( \Omega_+)\\
& \quad \hspace{ 3cm}
 \mbox{and transmission condition on } \partial  \Omega_- \}\,.
\end{array}
\end{equation}
It is clear that this implies $u_- \in H^2( \Omega_-)$.  The question
of showing that $u_+ \in H^2( \Omega_+)$ is {\it a priori} unclear.
By using the local regularity, we can show that for any $\chi$ in
$C^\infty_0 (\overline{ \Omega_+})$,
\begin{equation*}
(-h^2 \Delta + i V) (\chi u) \in L^2(\mathbb R^n)\,,
\end{equation*}
and consequently $\chi u \in H^2(\mathbb R^n)$. \\
In order to show that $u_+ \in H^2(\Omega_+)$, one needs to introduce
other techniques and additional assumptions.  For example, using the
pseudodifferential calculus, it is possible to prove (see \cite{Rob}),
that $u_+\in H^2 (\Omega_+)$ and $V u_+\in L^2(\Omega_+)$ under the
stronger condition that for any $\alpha \in \mathbb N^n$, there exists
$C_\alpha$ such that
\begin{equation}\label{condsup2}
|D_x^\alpha V(x) | \leq C_\alpha  \sqrt{1 + V(x)^2}\,,\quad \forall x \in \mathbb R^n\,.
\end{equation}

\subsubsection{The case of two unbounded components in $\mathbb R^2$ separated by a curve}

The case of two half-spaces is of course the simplest because we can
come back to the one-dimensional problem using the partial Fourier
transform.  The analysis of the resolvent should however be detailed
(see Henry \cite{Hen1} who treats the model of the half-space for the
BT operator with Neumann or Dirichlet conditions).  In fact, we
consider the quadratic form
\begin{equation*}
\begin{split}
q(u) & = h^2 \int_{x_1 < 0} |\nabla u_-(x)|^2\, dx + i \int_{x_1 < 0} \ell(x) |u_-(x) |^2 \, dx \\ 
 & +  h^2 \, \int_{x_1 > 0} |\nabla u_+(x)|^2\, dx + i\, \int_{x_1 > 0} \ell(x) |u_+(x) |^2 \, dx \\
 & + \Kappa\, \int | u_-(0,x_2) - u_+(0,x_2)|^2 \, dx_2\,, \\
\end{split}
\end{equation*}
where $x\mapsto \ell (x)$ is a nonzero linear form on $\mathbb R^2$:
\begin{equation*}
\ell (x) = \alpha x_1 + \beta x_2\,.
\end{equation*}
Here, we can also apply the general Lax-Milgram theorem in order to
define a closed operator associated to this quadratic form.  The
extension to a more general curve should be possible under the
condition that the curve admits two asymptotes at infinity.\\

In this section, we have described how to associate to a given
sesquilinear form $\mathfrak a$ defined on a form domain $\mathcal V$
an unbounded closed operator $\mathcal A$ in some Hilbert space
$\mathcal H$.  We will add the superscript $\#$ with $ \#\in \{ D, N,
R, T\}$ in order to treat simultaneously the different cases.  The
space $\mathcal H^\#$ will be $L^2(\Omega)$ when $\# \in \{D, N,R\}$
and will be $L^2(\Omega_-)\times L^2(\Omega_+)$ in the case with
transmission $\#=T$. $\mathcal V^\#$ will be respectively
$H_0^1(\Omega)$, $H^1(\Omega)$, $H^1(\Omega)\,$, and
$H^1(\Omega_-)\times H^1(\Omega_+)$.  The corresponding operators are
denoted $\mathcal A^\#_h$ with $ \#\in \{D,N,R,T\}$.

\section{Former results}
\label{sec:former}

\subsection{Spectral results for bounded domains}

For bounded domains, there are standard theorems, coming back to Agmon
\cite{Agm}, permitting to prove the non-emptiness of the spectrum and
moreover the completeness of the ``generalized'' eigenfunctions%
\footnote{
By this we mean elements in the kernel of $(\mathcal
A_h^\#-\lambda)^k$ for some $k\geq 1\,$.}.
In the case $V(x) = gx_1$ (here we can think of $g\in \mathbb C)$, the
limit $g\longrightarrow0$ can be treated by regular perturbation
theory.  In particular, Kato's theory \cite{Kato} can be applied, the
spectrum being close (modulo $\mathcal O (g)$) to the real axis.  It
is interesting to determine the variation of the lowest real part of
an eigenvalue.

For the Dirichlet problem, the Feynman-Hellmann formula gives the
coefficient in front of $g$ as $i \int_\Omega x_1|u_0(x)|^2\, dx$,
where $u_0$ is the first $L^2(\Omega)$-normalized eigenfunction of the
Dirichlet Laplacian.  In fact, using the standard Kato's procedure we
can look for an approximate eigenpair $(\lambda,u)$ in the form:
\begin{equation}\label{aa1}
u = u_0+ i \,g \, u_1 + g^2 \, u_2+ \dots 
\end{equation}
and 
\begin{equation}\label{aa2}
\lambda =\lambda_0 + i \,g\, \lambda_1 + g^2\, \lambda_2 + \dots
\end{equation}
Developing in powers of $g$, we get for the coefficient in front of
$g$:
\begin{equation}\label{aa3}
(-\Delta  -\lambda_0) u_{1}  =  - x_1 u_0 + \lambda_{1}  u_0\,,
\end{equation}
and $\lambda_{1}$ is chosen in order to solve \eqref{aa3}
\begin{equation}\label{aa4}
\lambda_{1}= \int_\Omega x_1 |u_0(x)|^2\, dx\,.
\end{equation}
We then take 
\begin{equation} \label{aa5}
u_{1} = -  (-\Delta -\lambda_0)^{(-1,reg)}\, \left( (x_1-\lambda_{1}) u_0\right)\,,
\end{equation}
where $(-\Delta -\lambda_0)^{(-1,reg)}$ is the regularized resolvent,
defined on the vector space generated by $u_0$ as
\begin{equation*}
(-\Delta -\lambda_0)^{(-1,reg)} u_0=0\,,
\end{equation*}
and as the resolvent on the orthogonal space to $u_0\,$. \\
To look at the coefficient in front of $g^2$, we write
\begin{equation}\label{aa6}
(-\Delta  -\lambda_0) u_{2}  = ( x_1-\lambda_{1} )\, u_{1} + \lambda_{2}  u_0\,,
\end{equation}
and get 
\begin{equation*}
\lambda_{2} =  - \int_{\Omega} (x_1 -\lambda_{1})  u_{1} (x) u_0(x) \, dx\,,
\end{equation*}
from which
\begin{equation*}
\lambda_{2} =  \bigl\langle (-\Delta -\lambda_0)^{(-1,reg)} ( (x_1-\lambda_0)u_0) \,|\,  ( (x_1-\lambda_0)u_0)\bigr \rangle_{L^2(\Omega)} >0 \,.
\end{equation*}
The effect of the perturbation is thus to shift the real part of the
``first'' eigenvalue on the right.\\
 
The limit $g\rightarrow +\infty$ for a fixed domain, or the limit of
increasing domains (i.e. the domain obtained by dilation by a factor
$R\to +\infty$) for a fixed $g$ can be reduced by rescaling to a
semi-classical limit $h\to 0$ of the operator $\A_h$ with a fixed
potential $V(x)$.
In this way, the BT-operator appears as a particular case (with $V(x)
= x_1$) of a more general problem.  We can mention (and will discuss)
several recent papers, mainly devoted to the Dirichlet case,
including: Almog \cite{Alm}, Henry \cite{Hen1} (Chapter 4),
Beauchard-Helffer-Henry-Robbiano \cite{BHHR} (analysis of the 1D
problem), Henry \cite{Hen2}, Almog-Henry \cite{AlHe} and in the
physics literature \cite{deSwiet94,Grebenkov07} (and references
therein).

\subsection{Spectral results for unbounded domains}

In the case of unbounded domains with bounded transmission boundary as
defined in Sec. \ref{sec:Def_unbouded}, there is no compact resolvent.
We note indeed that the pairs $(u_-,u_+)$ with $u_-=0$ and $u_+\in
C_0^\infty( \Omega_+)$ belong to the domain of the operator.  It is
easy to construct a sequence of $L^2$ normalized $u_+^{(k)}$ in
$C_0^\infty( \Omega_+)$ which is bounded in $H^2( \Omega_+)$, with
support in $(-R,+R) \times \mathbb R^{n-1}$, and weakly convergent to
$0$ in $L^2( \Omega_+)$.  This implies that the resolvent cannot be
compact.

The noncompactness of the resolvent does not exclude the existence of
eigenvalues.  Actually, when $\Kappa=0\,$, the spectral problem is
decoupled into two independent problems: the Neumann problem in $
\Omega_-$ which gives eigenvalues (the potential $ix_1$ in $ \Omega_-$
is just a bounded perturbation, as discussed in
Sec. \ref{sec:bounded_trans}) and the Neumann problem for the exterior
problem in $ \Omega_+$ with $-\Delta + i g x_1$ for which the question
of existence of eigenvalues is more subtle if we think of the model of
the half-space analyzed in Almog \cite{Alm} or \cite{Hen1}.  We will
see that in the semi-classical limit (or equivalently $g\rightarrow
+\infty$) the points of $\pa \Omega_-$ at which the normal vector to
$\pa \Omega_-$ is parallel to $(1,0,\dots,0)$, play a particular role.

\subsection{Semi-classical results}

In order to treat simultaneously various problems we introduce
$\Omega^\#$ with $\# \in \{D,N,R,T\}$ and $\Omega^D=\Omega\,$,
$\Omega^N=\Omega\,$,  $\Omega^R=\Omega\,$ and $\Omega^T = \Omega_-\,$.\\
R. Henry \cite{Hen2} (see also \cite{AlHe}) looked at the Dirichlet
realization of the differential operator
\begin{equation}\label{defAh}
\mathcal A_h^D := - h^2 \Delta+ i \, V(x)\,,
\end{equation}
in a fixed bounded domain $\Omega$, where $V$ is a real potential and
$h$ a semi-classical parameter that goes to $0$.\\
Setting $V(x) =x_1$, one gets a problem considered by de Swiet and Sen
\cite{deSwiet94} in the simple case of a disk but these authors
mentioned a possible extension of their computations to more general
cases.

For a bounded regular open set, R. Henry in \cite{Hen2} (completed by
Almog-Henry \cite{AlHe}, see below) proved the following
\begin{theorem}\label{thmNSAschro1}
Let $V\in\mathcal{C}^\infty(\overline{\Omega} ; \mathbb{R})$ be such that,
for every $x\in\overline{\Omega}\,$,
\begin{equation} \label{notdeg}
\nabla V(x)\neq0\,.
\end{equation}
Then, we have
\begin{equation}\label{limSpect1}
 \varliminf\limits_{h\to0}\frac{1}{h^{2/3}}\inf \bigl\{\Re\, \sigma(\mathcal A _h^D) \bigr\} \geq  \frac{|a_1|}{2}J_m^{2/3}\,,
\end{equation}
where $\mathcal A _h^D$ is the operator defined by (\ref{defAh}) with
the Dirichlet condition, $a_1<0$ is the rightmost zero of the Airy
function $Ai\,$, and 
\begin{equation}\label{defJm}
  J_m = \min_{x\in\partial\Omega_\perp}|\nabla V(x)|\,,
\end{equation}
where $$\partial \Omega_\perp= \{ x\in \partial \Omega\,,\, \nabla V (x) \wedge \nu (x) =0\}\,.$$
\end{theorem}
This result is essentially a reformulation of the result stated by
Y. Almog in \cite{Alm}.

\begin{remark}
The theorem holds in particular when $V(x)=x_1$ in the case of the
disk (two points) and in the case of an annulus (four points).  Note
that in this application $J_m =1$.
\end{remark}
A similar result can be proved for the Neumann case.
\begin{remark}
To our knowledge, the equivalent theorems in the Robin case and the
transmission case are open.  We hope to come back to this point in a
future work.
\end{remark}
A more detailed information is available in dimension $1$ (see
\cite{BHHR}) and in higher dimension \cite{AlHe} under some additional
assumption on $\partial \Omega_\perp$.  The authors in \cite{AlHe}
prove the existence of an approximate eigenvalue.  Our main goal is to
propose a more general construction which will work in particular for
the case with transmission condition.

\begin{remark}[Computation of the Hessian]  
For a planar domain, let us denote by $(x_1(s),x_2(s))$ the
parameterization of the boundary by the arc length $s$ starting from
some point, $t(s) = (x_1'(s), x'_2(s))$ is the normalized oriented
tangent, and $\nu(s)$ is the outwards normal to the boundary at $s$.
Now we compute at $s=0$ (corresponding to a point $x^0=x(0)\in
\partial \Omega^\#_\perp$, where $\nabla V\cdot t(0) =0\,$),
\begin{equation*}
\begin{split}
\left(\frac{d^2}{ds^2} V(x_1(s), x_2(s))\right)_{s=0}& = \langle t(0) | {\rm Hess V} (x_1(0),x_2(0))\,|  t(0) \rangle \\
 &   - \curv(0) \left( \nabla V (x_1(0),x_2(0)) \cdot \nu(0)\right) \,,  \\
\end{split}
\end{equation*}
where we used $ t'(s) = - \curv(s) \nu(s)\,$, $\curv(s)$ representing
the curvature of the boundary at the point $x(s)$.
\end{remark}

\begin{example}
When $V(x_1,x_2) = x_1$, we get
\begin{equation*}
\left(\frac{d^2}{ds^2} V(x_1(s), x_2(s))\right)_{s=0}  =  - \curv(0)  (e_1 \cdot \nu (0))  \,,
\end{equation*}
with $e_1= (1,0)\,$.\\
In the case of the disk of radius $1$, we get 
\begin{equation}\label{ex1}
\left( \frac{d^2}{ds^2} V(x_1(s), x_2(s)) \right)_{s=0} \bigl( e_1\cdot  \nu(0) \bigr) = - 1\,,
\end{equation}
for $(x_1,x_2) = (\pm 1,0)$.
\end{example}

Let us now introduce a stronger assumption for $\# \in \{N,D\}$.\\
\begin{assumption}\label{nondeg3}
At each point $x$ of $\partial \Omega^\#_\perp$, the Hessian of
$V_{/\partial \Omega}$ is
\begin{itemize}
\item positive definite if $\partial_\nu  V < 0\,$,
\item negative definite if $\partial_\nu  V > 0\,$,
\end{itemize}
with $\nu$ being the outwards normal and $\partial_\nu V:=\nu\cdot
\nabla V$.
\end{assumption}

Under this additional assumption%
\footnote{
We actually need this assumption only for the points $x$ of $\partial
\Omega_\perp$ such that $|\nabla V(x)|= J_m\,$.},
the authors in \cite{AlHe} (Theorem 1.1) prove the equality in
\eqref{limSpect1} by proving the existence of an eigenvalue
near each previously constructed approximate eigenvalue, and get a
three-terms asymptotics.

\begin{remark}
Note that this additional assumption is verified for all points of
$\partial \Omega_\perp$ when $V(x)=x_1$ and $\Omega$ is the disk.  In
fact, for this model, there are two points $(-1,0)$ and $(1,0)$, and
formula \eqref{ex1} gives the solution.
\end{remark}

Y. Almog and R. Henry considered in \cite{Alm,Hen2,AlHe} the Dirichlet
case but, as noted by these authors in \cite{AlHe}, one can similarly
consider the Neumann case.

Without Assumption \ref{nondeg3}, there is indeed a difficulty for
proving the existence of an eigenvalue close to the approximate
eigenvalue.  This is for example the case for the model operator
\begin{equation*}
-h^2 \frac{d^2}{dx^2} - h^2 \frac{d^2}{dy^2} + i (y-x^2)\,,
\end{equation*}
on the half space.  The operator is indeed not sectorial, and Lemma
4.2 in \cite{AlHe} is not proved in this case.  The definition of the
closed operator is questionable.  One cannot use the technique given
in a previous section because the condition \eqref{condsup} is not
satisfied.  The argument used by R. Henry in \cite{Hen1} for the
analysis of the Dirichlet BT-operator in a half space $\mathbb R^2_+$
(based on \cite{ReSi} (Theorem X.49) and \cite{Ich}) can be extended
to this case. \\
This problem occurs for the transmission problem in which the model
could be related to
\begin{equation*}
-h^2 \frac{d^2}{dx^2} - h^2 \frac{d^2}{dy^2} + i  (y+x^2)\,,
\end{equation*}
on the whole space $\mathbb R^2$ with transmission on $y=0\,$.  This
case will not be treated in this paper.

\subsubsection{On the growth of semi-groups}

In the case of Dirichlet and Neumann realizations, one can study the
decay of the semi-group $\exp(-t\A_h^\#)$ relying on the previous
results and additional controls of the resolvent (see \cite{Hen1},
\cite{AlHe}).  When the domain is bounded, the potential is a bounded
perturbation of self-adjoint operators.  In this case, the control of
the resolvent when $\Im \lambda$ tends to $\pm \infty$ is
straightforward, with the decay as $\mathcal O (1/|\Im \lambda|)$.
Applying the Gearhardt-Pr\"uss theorem (see for example in \cite{H}),
the decay is
\begin{equation*}
\mathcal O_\epsilon \biggl(\exp \bigl( - t (1-\epsilon)  \inf_{\lambda \in \sigma(\A_h^\#)} \{\Re \lambda \} \bigr) \biggr)\quad \forall \epsilon >0\,,
\end{equation*}
where $\sigma(\A)$ denotes  the spectrum of $\A$.  In this case, $\sigma
(\A_h^\#)$ is not empty and the set of generalized eigenfunctions is
complete (see \cite{Agm}).\\
In the unbounded case, the situation is much more delicate.  The
spectrum $\sigma (\A_h^\#)$ can be empty and one has to control the
resolvent as $|\Im \lambda| \rightarrow +\infty\,$.  The behavior of
the associate semi-group can be super-exponential when $\sigma
(\A_h^\#)$ is empty. Moreover, it is not granted that $\inf_{\lambda
\in \sigma(\A_h^\#)} \{\Re \lambda \}$ gives the decay rate of the
semi-group. \\

\section{Quasimode constructions -- Preliminaries}
\label{sec:quasimodesP}

Let us present in more detail the situation considered in Theorem
\ref{th:app}.

\subsection{Local coordinates}

Choosing the origin at a point $x^0$ at which $\nabla V (x^0)
\wedge \nu (x^0)=0$, we replace the Cartesian coordinates $(x_1,x_2)$
by the standard local variables $(s,\t)$, where $\t$ is the signed
distance to the boundary, and $s$ is the arc length starting from
$x^0$.  Hence
\begin{itemize}
\item 
In the case of one component, $\t=0$ defines the boundary $\partial
\Omega$ and $\Omega$ is locally defined by $\t>0\,$.

\item 
In the case of two components, $\t=0$ defines $\partial \Omega_-$,
while $\t < 0$ and $\t > 0$ correspond, in the neighborhood of
$\partial \Omega_-$, respectively to $ \Omega_-$ and $ \Omega_+\,$.
\end{itemize} 
In the $(s,\t)$ coordinates, the operator reads
\begin{equation}\label{defLhV} 
\mathcal \A_h =- h^2 a^{-1}\partial_s(a^{-1}\partial_s)- h^2 a^{-1}\partial_\t(a\partial_\t) + i \, \widetilde V (s,\t)  \,,
\end{equation}
with
\begin{equation*}
\widetilde V (s,\t):= V (x_1(s,\t),x_2(s,\t))\, ,
\end{equation*}
where
\begin{equation}
a(s,\t) = 1- \curv(s) \, \t \,,
\end{equation}
$\curv(s)$ representing the curvature of the boundary at $x(s,0)$.
Once the formal quasimodes are constructed in local coordinates, one
can return to the initial coordinates by using a standard Borel
procedure with cut-off functions, see Remark \ref{rem:cutoff}.\\
For future computation, we also rewrite (\ref{defLhV}) as
\begin{equation}\label{a1b}
 \A_h=- h^2 a^{-2}\partial_s^2  +  h^2 a^{-3}(\partial_s a)\, \partial_s  - h^2 \partial_\t^2 
- h^2a^{-1}(\partial_\t a) \,\partial_\t  + i \,  \widetilde  V(s,\t) \,.
\end{equation}
The boundary conditions read
\begin{itemize}
\item Dirichlet condition
\begin{equation}
 u (s,0) =0\,,
\end{equation}
\item Neumann condition 
\begin{equation}
\partial_\t u (s,0) =0\,,
\end{equation}
\item   Robin condition  with parameter $\Kappa$
\begin{equation}
h^2 \partial_\t u (s,0) = \Kappa u (s,0)  \,,
\end{equation}
\item Transmission condition with parameter  $\Kappa$
\begin{equation} 
\left\{
\begin{array}{l}
\partial_\t u_+ (s,0) = \partial_\t u_- (s,0)\,,\\
h^2 \partial_\t u_+ (s,0) = \Kappa \, \bigl(u_+(s,0) - u_-(s,0)\bigr)\,.
\end{array}
\right.
\end{equation}
\end{itemize}
In the last two cases, the link between $\Kappa$ and $h$ will be given
later in \eqref{assT}.

We omit the tilde of $\tilde{V}$ in what follows.\\

We recall that the origin of the coordinates is at a point $x^0$ such that
\begin{equation*}
\nabla V (x^0) \neq 0 \quad \mbox{ and } \quad \nabla V (x^0) \wedge \vec \nu (x^0) =0\,.
\end{equation*}
Hence we have
\begin{equation}\label{a2}
\frac{\partial V}{\partial s} (0,0) =0\,,
\end{equation}
and 
\begin{equation}\label{a3}
\frac{\partial V}{\partial \t} (0,0) \neq 0\,.
\end{equation}
We also assume in our theorem that 
\begin{equation}\label{a4}
\frac{\partial^2 V}{\partial s^2}(0,0) \neq 0\,.
\end{equation}
Hence we have the following Taylor expansion
\begin{equation}
V(s,\t) \sim \sum_{j,k} v_{jk} s^j \t^k\,,
\end{equation}
where
\begin{equation}
\label{eq:Vjk}
v_{jk} = \frac{1}{j! ~k!} ~ \left(\frac{\partial^{j+k}}{\partial s^j \partial \t^k} V(s,\t)\right)_{s=\rho=0}\, ,
\end{equation}
with
\begin{equation}\label{condV}
v_{00} = V(0,0)\,,\quad   v_{10} =0\,,\quad   \, v_{01} \neq 0\,,\quad  v_{20} \neq 0\,,
\end{equation}
corresponding to the assumptions of Theorem \ref{th:app}.

\subsection{The blowing up argument}

Approximating the potential $V$ near $x^0$ by the first terms of its
Taylor expansion $v_{00} + v_{01} \t + v_{20} s^2\,$, a basic model
reads
\begin{equation*}
-h^2 \frac{d^2}{ds^2} - h^2 \frac{d^2}{d\t^2} + i \, (v_{01} \t+ v_{20} s^2)  \quad \textrm{on the half space}~\{ \t >0\} ,
\end{equation*}
in the case when $\#\in \{D,N,R\}$, and on $\mathbb R^2$ when $\#=T$,
which is reduced by a natural scaling
\begin{equation}\label{dil1}
(s,\t) = (h^{\frac 12} \sigma , h^{\frac 23} \tau)
\end{equation}
to
\begin{equation*}
h \biggl(- \frac{d^2}{d\sigma^2} + i v_{20} \sigma^2\biggr) + h^\frac 23  \biggl(- \frac{d^2}{d\tau^2} + i v_{01} \tau \biggr)\,,
\end{equation*}
whose definition and spectrum can be obtained by separation of
variables in the four cases.\\

\subsubsection{Expansions}

In the new variables $(\sigma,\tau)$ introduced in \eqref{dil1}, the
expansion is
\begin{equation}
\widehat V_h (\sigma,\tau) :=V(h^\frac 12 \sigma,h^\frac 23 \tau) \sim \sum_{m\geq 0} h^\frac{m}{6} \left (\sum_{ 3k +4p =m} v_{kp} \sigma^k\tau^p\right) \,.
\end{equation}
In particular, the first terms are
\begin{equation}
\widehat V_h (\sigma,\tau)
  \, = v_{00} + h^\frac 23 \, v_{01} \tau + h v_{20} \sigma^2
+ h^\frac 7 6 v_{11} \sigma \tau+ h^{\frac 43}  v_{02} \tau^2 + h^{\frac 32}v_{30} \sigma^3+ \mathcal O (h^\frac 53)\,.
\end{equation}
\\
Similarly, we consider the dilation of $a (s,\rho)$
\begin{equation}
\widehat a_h (\sigma,\tau):=a(h^\frac 12 \sigma,h^\frac 23 \tau)  = 1- h^\frac 23 \tau ~  \curv(h^\frac 12 \sigma) ,
\end{equation}
which can be expanded in the form
\begin{equation}
\widehat a_h (\sigma,\tau) \sim  1 - h^\frac 23 \tau \left(\sum_\ell\frac{1}{\ell!}\curv^{(\ell)}(0) \sigma^\ell  h^{\frac \ell 2}\right) \,.
\end{equation}
In the $(\sigma,\tau)$ coordinates, we get
\begin{equation}\label{a6}
 \widehat{\A}_h=- h \widehat a_h^{-2}\partial_\sigma^2  +  h^\frac 32 \widehat a_h^{-3}\widehat{(\partial_s a)}_h \, \partial_\sigma  
-h^\frac 23 \partial_\tau^2 - h^\frac 43  \widehat a_h^{-1} \widehat{(\partial_\t a)}_h \,\partial_\tau  + i\,  \widehat V_h(\sigma, \tau)\,.
\end{equation}
We note that
\begin{equation*}  
\widehat{(\partial_s a)}_h (\sigma, \tau) = - h^\frac 23 \tau \curv'( h^{\frac 12}\sigma)
\quad \mbox{ and } \quad 
\widehat{(\partial_\t a)}_h (\sigma,\tau)= - \curv  (h^\frac 12 \sigma)\,.
\end{equation*}
We rewrite $\widehat{\A}_h$ by expanding in powers of $h^\frac 16$:
\begin{equation}
\label{eq:Ah_series}
\widehat{\A}_h \sim  i \, v_{00} + h^\frac 23\,  \sum_{j\geq 0}  h^\frac j6 \L_{j} (\sigma, \tau, \partial_\sigma, \partial_\tau)
\end{equation}
or, equivalently, as
\begin{equation}
\label{eq:Ah_seriesbis}
\begin{split}
h^{-\frac 23} (\widehat{\A}_h - i v_{00})& = - \partial_\tau^2 - h^\frac 13  \widehat a_h^{-2}\partial_\sigma^2  
- h^\frac 23  \widehat a_h^{-1} \widehat{(\partial_\t a)}_h \,\partial_\tau  \\ 
& \quad +  h^\frac 56  \widehat a_h^{-3}\widehat{(\partial_s a)}_h \, \partial_\sigma 
+ i\,  h^{-\frac 23} (\widehat V_h(\sigma, \tau) - v_{00})\\
 & \sim \,  \sum_{j\geq 0}  h^\frac j6 \L_{j} (\sigma, \tau, \partial_\sigma, \partial_\tau)\,, \\
\end{split}
\end{equation}
where the first terms are given by
\begin{equation}
\begin{array}{ll}
\L_0 & = -\partial_\tau^2 + i \, \, v_{01}\, \tau\,,\\
\L_1 & = 0\,,\\
\L_2 &= -\partial_\sigma^2 + i\, v_{20} \, \sigma^2\,,\\
\L_3 &= i \,v_{11}\, \sigma \tau\,,\\
\L_4 &= \curv(0) \, \partial_\tau +  i\,v_{02}\, \tau^2\,,\\
\L_5 & = i\,v_{03}\,  \sigma^3\,.
\end{array}
\end{equation}
For any $j\geq 0$, each $\mathcal L_j$ is a differential operator of
order $\leq 2$ with polynomial coefficients of degree which can be
controlled as a function of $j\,$.  In particular these operators
preserve the vector space $\mathcal S (\mathbb R_\sigma) \otimes
\mathcal S^\#$.  The Fr\'echet space $\mathcal S^\#$ denotes $\mathcal
S (\overline {\mathbb R_+})$ in the case when $\#\in \{D,N,R\}$ and
$\mathcal S (\overline{\mathbb R_-})\times
\mathcal S (\overline{\mathbb R_+})$ when $\#=T$.

\subsubsection{Parity}\label{sssparity}

Note also that we have 
\begin{lemma}
\begin{equation}\label{sym}
\check {(\mathcal L_j  f)} = (-1)^{j}\mathcal L_j \check f \, ,
\end{equation}
where $\check f (\tau,\sigma) = f (\tau,-\sigma)\,.$
\end{lemma}
{\bf Proof}\\
This is a consequence of
\begin{equation} \label{jsym}
\L_{j} (\sigma, \tau, \partial_\sigma, \partial_\tau)=(-1)^j \,\L_{j} (-\sigma, \tau, -\partial_\sigma, \partial_\tau)
\end{equation}
that can be seen by observing \eqref{eq:Ah_seriesbis}.  We will see
that each term in the right hand side of
\eqref{eq:Ah_seriesbis} satisfies \eqref{sym}.\\
First, denoting $\hat h = h^\frac 16$, we can rewrite
\begin{equation}
\widehat a_h (\sigma,\tau) \sim  1 - \hat h^ 4  \tau \left(\sum_{\ell \geq 0} \frac{1}{\ell!}\curv^{(\ell)}(0) \sigma^\ell  \hat h^{3\ell}\right) \,,
\end{equation}
and expanding in powers of $\hat h$, we see that the coefficient in
front of $\hat h^\ell$ has the parity of $\ell$ in $\sigma$.  The same
is true for $\widehat a_h (\sigma,\tau)^{-2}$.  Hence the coefficient
in front of $h^\frac j6$ in $h^\frac 13 \widehat a_h
(\sigma,\tau)^{-2} \partial_\sigma^2$ satisfies \eqref{jsym}.  \\
We now look at $ h^\frac 56 \widehat a_h^{-3}\widehat{(\partial_s
a)}_h$ and write
\begin{equation*}
h^\frac 56 \widehat{(\partial_s a)}_h (\sigma, \tau) \partial_ \sigma = - \hat h ^9  \curv'( \hat h ^3 \sigma) \partial _\sigma\, .
\end{equation*}
It is clear from this formula that the second term in the right hand
side of \eqref{eq:Ah_seriesbis} satisfies \eqref{jsym}.\\
The third term $-\partial_\tau^2$ clearly satisfies \eqref{jsym}.  For
the forth term $- h^\frac 23 \widehat a_h^{-1} \widehat{(\partial_\t
a)}_h \,\partial_\tau$, it is enough to use the previous expansions
and to observe that
\begin{equation*}
\widehat{(\partial_\t a)}_h (\sigma,\tau)= - \curv  (\hat h^3 \sigma)\,.
\end{equation*}
Finally, we consider 
\begin{equation*}
 i\,  h^{-\frac 23} (\widehat V_h(\sigma, \tau) - v_{00}) \sim  i \, \sum_{m\geq 4} \hat h^{m -4} 
\left (\sum_{ 3k +4p =m} v_{kp}\, \sigma^k \, \tau^p\right)\,,
\end{equation*}
and we observe that $k$ and $m$ should have the same parity. \qed  \\
This lemma will be useful for explaining cancellations in the
expansion of the quasimode.

\subsubsection{Boundary or transmission conditions}
In these local coordinates, the boundary conditions read
\begin{itemize}
\item the Dirichlet condition
\begin{equation}
u (\sigma,0) =0\,,
\end{equation}
\item the Neumann condition 
\begin{equation}
\partial_\tau u (\sigma,0) =0\,,
\end{equation}
\item the Robin condition 
\begin{equation}
\partial_\tau u (\sigma,0) = \Kappa h^{-\frac 43}  u (\sigma,0)  \,,
\end{equation}
\item the transmission condition 
\begin{equation}\label{assT0}
\partial_\tau u_- (\sigma,0)= \partial_\tau u_+ (\sigma,0)\,,\quad  \partial_\tau u_+ (\sigma,0) = 
 \Kappa \, h^{-\frac 43} \bigl( u_+(\sigma,0) - u_-(\sigma,0)\bigr) \,.
\end{equation}
\end{itemize}
Depending on the physical problem, the Robin or Transmission parameter
$\Kappa$ can exhibit different scaling with $h$.  Here we assume the
scaling
\begin{equation} \label{assT}
\Kappa = \kappa \, h^\frac 43\, ,
\end{equation}
so that the Robin or transmission conditions in the variables
$(\sigma,\tau)$ are independent of $h$ and read
\begin{equation}\label{assR0a}
  \partial_\tau u (\sigma,0) = \kappa \, u(\sigma,0)  \,,
\end{equation}
and
\begin{equation}\label{assT0a}
\partial_\tau u_- (\sigma,0)= \partial_\tau u_+ (\sigma,0)\,,\quad  \partial_\tau u_+ (\sigma,0) = 
\kappa \bigl( u_+(\sigma,0) - u_-(\sigma,0)\bigr) \,.
\end{equation}
In Sec. \ref{sec:large_domain}, we justify this scaling by considering
the transmission problem in dilated domains, while other scalings are
discussed in Sec. \ref{sec:kappa_scaling}.  We denote by $\L_0^\#$ the
realization of $\L_0$ with $\#= D,N, R, T$ for Dirichlet, Neumann,
Robin, or Transmission condition.  We recall that the Hilbert space
$L^2_\#$ denotes $L^2(\mathbb R_+)$ in the case when $\#\in
\{D,N,R\}$, and $L^2(\mathbb R_-)\times L^2(\mathbb R_+)$ when $\#=T$.
For the complex harmonic oscillator $\mathcal L_2$ we consider (with
the same notation) the self-adjoint realization on $L^2(\mathbb
R_\sigma)$.

\subsection{Comparison with the large domain limit}
\label{sec:large_domain}

We assume that $0\in \Omega_-$ and we dilate $\Omega_-$ and $\Omega$
by the map $(x_1,x_2)\mapsto (S x_1, S x_2)$ ($S>0$ supposed to be
large) and get $\Omega_-^S$ and $\Omega^S$.\\
It remains to check how the transmission problem for $\Omega^S$ with
$V(x)=x_1$ is modified by dilation.  If we start from the form
\begin{equation*}
u\mapsto  \int_{\Omega^S} |\nabla u|^2 dx + i\int_{\Omega^S}  x_1 |u(x)|^2\,dx  + \kappa_S \int_{\partial \Omega_-^S} |u_+- u_-| ^2 \,ds_S \,,
\end{equation*}
with a transmission coefficient $\kappa_S\,$, we get by the change of
coordinates $x=S y$, for $v(y)= u(S y)$\,,
\begin{equation*}
  \int_{\Omega} |\nabla_y v|^2 dy + i \, S^3 \int_{\Omega}  y_1 |v(y)|^2\,dy  + \kappa_S\, S \,  \int_{\partial  \Omega_-}  |v_+- v_-| ^2 \,ds \,.
\end{equation*}
Dividing by $S^3$, we get
\begin{equation*}
 \frac{1}{S^3} \int_{\Omega} |\nabla_y v|^2 dy + i  \int_{\Omega}  y_1 |v(y)|^2\,dy  + \kappa_S\, S^{-2} \,  \int_{\partial  \Omega_-}  |v_+- v_-| ^2 \,ds \,.
\end{equation*}
In order to treat this problem as semi-classical, we set
\begin{equation*}
h^2 = \frac{1}{S^3}\,,\qquad \Kappa =  \kappa_S\,  S^{-2}\,, 
\end{equation*}
Hence we get
\begin{equation*}
\kappa = \kappa_S\,,
\end{equation*}
and our assumption \eqref{assT} on $\Kappa$ corresponds to what we get
by rescaling from the problem in $\Omega_S\,$ with $\kappa_S$
independent of $S\,$.\\

For this application, Theorem \ref{th:app} gives the following
\begin{theorem}
For $S >0$, let $V_S (x) = S \, V (S^{-1}x)$, with the potential $V$
defined on $\Omega$ satisfying the conditions of Theorem \ref{th:app},
and $\kappa_S$ is independent of $S$.  Then, with the notation of
Theorem \ref{th:app}, one can construct a quasimode $\lambda_S^\#$ of
the $\#$ realization of the operator $-\Delta + i V_S$ in
$\Omega^\#_S$ such that
\begin{equation}
\label{eq:lambdah_auxR} 
\lambda_S ^\#= i \, S V(x^0)  +  \sum_{j\in \mathbb N} \lambda_{2j}^\#\,  S^{-\frac j 2}  +  \mathcal O (S^{-\infty})\,,
\end{equation}
as $S \rightarrow +\infty$.
\end{theorem}

This theorem can also be applied to $V(x)=x_1$, in which case $V_S$ is
independent of $S$.
\begin{remark}
More generally, one can consider
\begin{equation*}
V_S (x) = S^m \, V (S^{-1}x)\,,
\end{equation*} 
with $m >-2$.  In this case, we get $\kappa = \kappa_S\, S^{1-m}$.  If
$\kappa$ is independent of $S$ or tends to $0$ as $S\rightarrow
+\infty$, one can apply the semi-classical analysis of the previous
sections.
\end{remark}

\section{The quasimode construction. Proof of the main theorem}
\label{sec:quasimodes}

\subsection{The form of the quasimode}

In what follows, we assume in the Robin or transmission cases that
$\kappa$ is independent of $h$ (see \eqref{assT}).  We now look for a
quasimode $u_h^{app,\#}$ that we write in the $(\sigma,\tau)$
variables in the form:
\begin{equation}\label{expa1}
u_h^{app,\#}  \sim  d (h) \left( \sum_{j\geq 0} h^\frac j6 u_j^\#(\sigma, \tau)\right) \,,
\end{equation}
associated with an approximate eigenvalue
\begin{equation}\label{expa2} 
\lambda_h^{app,\#}\sim i \, v_{00} + h^\frac 23  \sum_{j\geq 0} h^\frac j6 \lambda^\#_j \,.
\end{equation}
Here $d(h) \sim d_0 h^{-\frac{7}{12}}$ with $d_0 \neq 0$ chosen such
that, coming back to the initial coordinates, the $L^2$-norm of the
trial state equals $1$. \\

Note that the $u_j^\#$ are in the domain of $\mathcal L_j^\#$
if we take the condition $\#$ (with $\#\in \{ N, D, R, T\}$).\\
Note also that we do not assume {\it a priori} that the $\lambda_j^\#$
for $j$ odd are $0$ as claimed in our theorem.  \\
As will be seen in the proof, we can choose
\begin{equation}\label{qm3}
u_j^\# (\sigma, \tau) = \phi_j^\# (\sigma) \psi_0^\#(\tau) \,,\quad  j=0, 1, 2\,,
\end{equation}
and
\begin{equation}\label{qm4}
u_j^\# (\sigma, \tau) = \phi_j^\# (\sigma) \psi_0^\#(\tau) +  \sum_{\ell=1} ^{N_j} \phi_{j,\ell} ^\#(\sigma) \psi_{j,\ell}^\# (\tau)\,,\quad j\geq  3\,,
\end{equation}
with $\phi_{j,\ell}^\# \in \mathcal S (\mathbb R) $ and $\psi_{j,\ell}
\in \mathcal S^\#$ to be specified below.\\
Moreover, we have
\begin{equation}\label{qm5}
\mathcal L_0^\# \psi_0^\# = \lambda_0^\# \psi_0^\#\,,
\end{equation}
\begin{equation}\label{qm6}
 \mathcal L_2  \phi_0^\# (\sigma) = \lambda_2^\# \phi_0^\# (\sigma)\,,
\end{equation}
with
\begin{equation}\label{qm8}
 \langle \psi_{j,\ell}^\# \,, \, \bar  \psi_0^\# \, \rangle_{L^2_\#}=0\,,
\end{equation}
and 
\begin{equation}\label{qm9}
  \langle \psi_{0}^\# \, , \, \bar  \psi_0^\# \rangle_{L^2_\#} \neq 0\,.
\end{equation}
The construction will consist in expanding $(\widehat {\mathcal A_h} -
\lambda_h^{app,\#}) u_{h}^{app,\#}$ in powers of $h^\frac 16$ and
finding the conditions of cancellation for each coefficient of this
expansion.\\
 
\begin{remark} \label{rem:cutoff}
If we succeed in this construction and come back to the initial
coordinates, using a Borel procedure to sum the formal expansions and
multiplying by a cutoff function in the neighborhood of a point $x^0$
of $\partial \Omega^\#$, we obtain an approximate spectral pair
localized near $x^0$ (i.e. $\mathcal O (h^\infty)$ outside any
neighborhood of $x^0$).  The Borel procedure consists in choosing a
cutoff function $\theta$ (with $\theta =1$ in a small neighborhood of
$0$ and a sequence $H_n$ such that $\beta \mapsto \sum_j \beta^j
\lambda_j \theta ( \beta /H_j)$ converges in $C^\infty ([0, \beta_0])$
for some $\beta_0>0$.  We then define
\begin{equation*}
  \lambda_h^\# = i\, v_{00} + h^\frac 23  \sum_{j\geq 0} \beta^j  \lambda_j \theta ( \beta /H_j)\,,
\end{equation*}
with $\beta = h^\frac 16$.\\
This $\lambda_h^{\#}$ is not unique but the difference between two
different choices is $\mathcal O (h^\infty)$.  A similar procedure can
be used to define a quasimode state $u_h^{\#}$ strongly localized near
$x^0$ (see \cite{AlHe,HK,HKR} for more details).
\end{remark}

\begin{remark}
We emphasize that the above construction is not sufficient (the
problem being non self-adjoint) for proving the existence of an
eigenvalue with this expansion.  The construction is true for any
regular domain (exterior or interior) under the conditions
\eqref{a2}-\eqref{a4}.  When $V(x)=x_1$, we recover in this way the
condition that the curvature does not vanish at $x^0$.  We recall that
this construction can be done near each point where $\nabla V (x^0)
\wedge \nu (x^0) =0$.  The candidates for the spectrum are determined
by ordering different quasimodes and comparing their real parts.  We
guess that the true eigenfunctions will have the same localization
properties as the constructed quasimode states.
\end{remark}

\subsection{Term $j=0\,$}

Identifying the powers in front of $h^\frac 16$, after division by
$d(h)$,  one gets the first equation corresponding to
$j=0$\,. \\ We consider four boundary conditions.\\

\subsubsection*{Neumann and Dirichlet cases}
For the Neumann boundary condition, one has
\begin{equation}
\L_0^N\,  u_0^N = \lambda_0^N u_0^N \,,\quad  \partial_\tau  u_0^N(\sigma,0)=0\,,
\end{equation}
and we look for a solution in the form
\begin{equation}
u_0^N(\sigma,\tau) = \phi_0^N(\sigma) \psi_0 ^N(\tau)\,.
\end{equation}
At this step, we only look for a pair $(\lambda_0^N,\psi_0^N)$ with
$\psi_0^N$ non identically $0$ such that
\begin{equation}
 (-\partial_\tau^2 + i \, v_{01} \tau)\,  \psi_0^N (\tau) = \lambda_0 ^N \psi_0^N (\tau) \mbox{ in } \mathbb R^+\,,\quad  (\psi_0^N)'(0)=0\,.
\end{equation}
We recall from \eqref{condV} that $\, v_{01} \neq 0$ so we have the
standard spectral problem for the complex Airy operator in the half
line with Neumann condition at $0\,$.  The spectral theory of this
operator is recalled in \cite{GHH}.  The spectrum consists of an
infinite sequence of eigenvalues $( \lambda^{N,(n)})_{n\geq 1}$
(ordered by increasing real part) that can be expressed through the
zeros $a'_n$ ($n\geq 1$) of the derivative of the Airy function
$\Ai'(z)$:
\begin{equation}
\label{eq:lambda0N}
 \lambda^{N,(n)} = - a'_n \, |\, v_{01}|^\frac 23 \,  \exp \left(\frac{i\pi}{3}  \,{\rm sign}\, \, v_{01}\right) \,.
\end{equation}
Different choices of $n$ will determine the asymptotic expansion of
different approximate eigenvalues of the original problem.  If we are
interested in controlling the decay of the associated semi-group, we
choose $\lambda_0^N= \lambda^{N,(1)}$ which corresponds to the
eigenvalue with the smallest real part.

One can similarly treat the Dirichlet problem (like in \cite{AlHe}).
In this case, one has
\begin{equation}
\L_0^D \, u_0 ^D= \lambda_0^D u_0^D \mbox{ in } \mathbb R_+ \,,\quad   u_0^D(\sigma,0)=0\,,
\end{equation}
and we look for a solution in the form
\begin{equation}
u_0^D(\sigma,\tau) = \phi_0^D(\sigma) \psi_0 ^D(\tau)\,,
\end{equation}
where $\psi_0^D(\tau)$ satisfies
\begin{equation}
\L_0^D \, \psi_0^D = \lambda_0^D \psi_0^D \mbox{ in } \mathbb R_+ \,,\quad   \psi_0^D(0)=0\,.
\end{equation}
The spectral theory of this operator is also recalled in \cite{GHH}.
The spectrum consists of an infinite sequence of eigenvalues $(
\lambda^{D,(n)})_{n\geq 1}$ (ordered by increasing real part) that can
be expressed through the zeros $a_n$ ($n\geq 1$) of the Airy function
$\Ai(z)$:
\begin{equation}
\label{eq:lambda0D}
 \lambda^{D,(n)} = - a_n \, |\, v_{01}|^\frac 23 \,  \exp \left(\frac{i\pi}{3}  \,{\rm sign}\, \, v_{01}\right) \,.
\end{equation}

One can show (see \cite{Hen1} for a proof by analytic dilation) that
\begin{equation}\label{normcond}
 \int_0^{+\infty} \psi_0^N (\tau)^2\, d\tau \neq 0\, \quad \mbox{ and } \quad   \int_0^{+\infty} \psi_0^D (\tau)^2\, d\tau \neq 0\,.
\end{equation}
This is also a consequence of the completeness of the eigenfunctions
of the complex Airy operator in the half-line with Neumann or
Dirichlet boundary condition.  This property is true for any
eigenvalue $\lambda_0^\#$ of $\mathcal L_0^\#$. \\

For $n\geq 1$, the eigenfunctions $\psi_0^N=\psi^{N,(n)}$ and
$\psi_0^D=\psi^{D,(n)}$ are specifically translated and complex
dilated Airy functions:
\begin{eqnarray}
\label{eq:psi0N}
\psi^{N,(n)}(\tau) &=& c^N_n \, \Ai\left(a_n' + \tau \, |\, v_{01}|^{\frac 13} \exp\left(\frac{i\pi}{6}\, \sign \, v_{01}\right) \right) 
\quad \mbox{ for }  \tau \geq 0\,,\;\\
\label{eq:psi0D}
\psi^{D,(n)}(\tau) &=& c^D_n \,  \Ai\left(a_n + \tau \, |\, v_{01}|^{\frac 13} \exp\left(\frac{i\pi}{6}\, \sign \, v_{01}\right) \right)  
\quad \mbox{ for }  \tau \geq 0\,,\;
\end{eqnarray}
where the normalization constants $c^N_n$ and $c^D_n$ can be fixed by
choosing the following normalization that we keep throughout the
paper:
\begin{equation*}
 \int\limits_0^\infty \psi_0^\# (\tau)^2 d\tau =1\,.
\end{equation*}
These coefficients are computed explicitly in Appendix
\ref{sec:lambda4}  (see \eqref{eq:CN}, \eqref{eq:CD}) .

\subsubsection*{Robin case} 
For the Robin boundary condition, one has
\begin{equation}
\L_0^R\,  u_0^R = \lambda_0^R u_0^R \,,\quad  \partial_\tau u_0^R(\sigma,0)= \kappa\,  u_0^R(\sigma,0)\,,
\end{equation}
and we look for a solution in the form
\begin{equation}
u_0^R(\sigma,\tau) = \phi_0^R(\sigma) \psi_0^R(\tau)\,,
\end{equation}
where the function $\psi_0^R(\tau)$ satisfies
\begin{equation}
(-\partial_\tau^2 + iv_{01} \tau ) \psi_0^R(\tau) = \lambda_0^R \psi_0^R(\tau) \, \mbox{ in } \, \R_+, 
\quad (\psi_0^R)'(0) = \kappa \, \psi_0^R(0)\,.
\end{equation}
This one-dimensional problem was studied in \cite{GHH}.  In
particular, the spectrum consists of an infinite sequence of
eigenvalues $(\lambda^{R,(n)})_{n\geq 1}$ (ordered by increasing real
part) that can be expressed as
\begin{equation}
\label{eq:lambda0R}
\lambda^{R,(n)} (\kappa) = - a_n^R(\kappa) \, |\, v_{01}|^{\frac23} \, \exp\biggl(\frac{\pi i}{3}\, \sign v_{01}\biggr) \,,
\end{equation}
where $a_n^R(\kappa)$ is a solution of the equation
\begin{equation}\label{eigenvalueR}
\exp\biggl(\frac{\pi i}{6}\, \sign v_{01}\biggr) \, \Ai'\bigl(a_n^R(\kappa)\bigr) -  \frac{\kappa}{|v_{01}|^{\frac13}} \, \Ai\bigl( a_n^R(\kappa) \bigr) = 0\,,
\end{equation}
and $\kappa \geq 0$ denotes the Robin parameter.%
\footnote{
In \cite{GHH}, we discussed the complex Airy operator with $v_{01} =
-1$, see Eq. (3.25).}\\
Except for the case of small $\kappa$, in which the eigenvalues are
close to the eigenvalues of the Neumann problem, it does
not seem easy to localize all the solutions of \eqref{eigenvalueR} in
general.  Note that from \eqref{eigenvalueR}, we deduce that
\begin{equation}
(\lambda^{R,(n)})'(0) =  - (a_n^R)'(0) \, |\, v_{01}|^{\frac23} \, \exp\biggl(\frac{\pi i}{3}\, \sign v_{01}\biggr)\,,
\end{equation}
where
\begin{equation}
\label{eq:anRprime}
 (a_n^R)'(0) =  \frac{1}{ a_n^R(0) |v_{01}|^{\frac13}} \, \exp\biggl(- \frac{\pi i}{6}\, \sign v_{01}\biggr) \ne 0\,.
\end{equation}
Nevertheless it is proven in \cite{GHH} that the zeros of the function
in (\ref{eigenvalueR}) are simple and that there is no Jordan block.
So as can be deduced from the next lemma, any eigenfunction satisfies
$\int \psi_0^R(\tau)^2\, d\tau \neq 0\,$.  We consequently fix the
normalization of $\psi_0^R$ by imposing
\begin{equation}
\int\limits_0^\infty \psi_0^R(\tau)^2\, d\tau =1\,.
\end{equation}
For $n\geq 1$, the associated eigenfunction $\psi_0^R=\psi^{R,(n)}$ reads
\begin{equation}
\label{eq:psi0R}
\psi^{R,(n)}(\tau) = c^R_n \, \Ai \biggl(a_n^R(\kappa) + \tau |v_{01}|^{\frac13} 
\exp\biggl(\frac{\pi i}{6}\, \sign v_{01}\biggr) \biggr)  \quad (\tau \geq 0),
\end{equation}
where $c^R_n$ is the normalization constant given by (\ref{eq:CR}).

\subsubsection{Transmission case}  In the transmission case, one gets, with
$\psi_0^{T} = (\psi_0^-,\psi_0^+)$,
\begin{equation}
\begin{array}{l}
\L_0 \psi_0^- =\lambda_0^{T} \psi_0^-  \mbox{ in } \mathbb R_-\,, \quad \L_0 \psi_0^+ =\lambda_0^{T} \psi_0^+ \mbox{ in } \mathbb R_+\,,\quad \\
\partial_\tau \psi_0^- (\sigma,0)= \partial_\tau \psi_0^+ (\sigma,0)\,,\quad \partial_\tau \psi_0^+ (\sigma,0) 
= \kappa \,  \bigl( \psi_0^+(\sigma,0) - \psi_0^-(\sigma,0)\bigr) \,.
\end{array}
\end{equation}
The existence of $\lambda_0^{T}$ has been proved in \cite{GHH}.  In
addition, the eigenvalue (of the smallest real part) is simple (no
Jordan block) for $\kappa \geq 0 $ small.  We can use the explicit
computations in \cite{GHH} or the following abstract lemma by
Aslayan-Davies for a closed operator $A$ \cite{AD}:
\begin{lemma}\label{abstlemma}
If $f$ and $f^*$ are the normalized eigenvectors of $A$ and $A^*$
associated with the eigenvalues $\lambda$ and $\bar \lambda$
respectively, and if the the spectral projector $P$
has rank $1$, then $\langle f\,,\, f^*\rangle \neq 0$ and
\begin{equation*}
||P|| = \frac{1}{|\langle f\,,\, f^*\rangle|}\,. 
\end{equation*}
\end{lemma}
The proof that $P$ has rank $1$ for the case $V(x) = x_1$ is given in
\cite{GHH} but only for $\kappa \geq 0\,$.  In general, we make
the assumption
\begin{assumption}\label{AassT}
$\lambda_0^{T}(\kappa)$ is simple (no Jordan block)\,.
\end{assumption} 
Under this assumption, we have 
\begin{equation}
\int_{-\infty}^\infty \psi_0^{T} (\tau)^2 \, d\tau :=   \int_{-\infty}^0 \psi_{0}^-(\tau)^2\,d\tau +  \int^{+\infty}_0 \psi_{0}^+(\tau)^2\,d\tau  \neq 0\,.
\end{equation}

The explicit form of the eigenfunctions $\psi^{T, (n)}$ ($n\geq 1$)
can be obtained from the analysis provided in \cite{Gr1,GHH}:
\begin{equation}
\label{eq:psi0T}
\begin{split}
\psi^{+, (n)}(\tau) & = - c^T_n\,  \bar \delta \,\Ai'\bigl(a_n^-(\kappa)\bigr) \, \Ai\left(a_n^+(\kappa) + \tau \, |\, v_{01}|^{\frac 13} \delta\right) , \\
\psi^{-, (n)}(\tau) & =  c^T_n \, \delta \, \Ai'\bigl(a_n^+(\kappa)\bigr)\, \Ai\left(a_1^-(\kappa) - \tau \, |\, v_{01}|^{\frac 13} \bar \delta \right) , \\
\end{split}
\end{equation}
where $c^T_n$ is a normalization constant (to satisfy
\eqref{normpsi0}), $\delta = \exp\left(\frac{i\pi}{6}\, \sign
\, v_{01}\right)$, and
\begin{equation}
\label{eq:ankappa}  
a_n^\pm(\kappa) = \hat \lambda^{T, (n)} \bigl(\kappa/ |\, v_{01}|^{\frac13} \bigr) \, \exp\left(\pm \frac{2\pi i}{3}~ \sign \, v_{01}\right)  ,
\end{equation}
where the $\hat \lambda^{T,(n)}(\check \kappa)$ are the eigenvalues of
the complex Airy operator $- \frac{d^2}{dx^2} + ix$ on the line with
transmission condition at $0$, with coefficient $$\check \kappa =
\kappa/|v_{01}|^\frac 13\,.$$  They are defined implicitly as
complex-valued solutions (enumerated by the index $n = 1,2,\ldots$) of
the equation \cite{Gr1,GHH}
\begin{equation}
\label{eq:an_equation}
2\pi\, \Ai'\bigl(e^{2\pi i/3} \hat \lambda^{T,(n)}(\check \kappa)\bigr)\, \Ai'\bigl(e^{-2\pi i/3} \hat \lambda^{T,(n)}(\check \kappa) \bigr) = - \check \kappa  \,.
\end{equation}
The eigenvalues $\hat \lambda^{T,(n)}(\check \kappa)$ are ordered
according to their increasing real parts: 
\begin{equation*}
\Re\{\hat \lambda^{T,(1)}(\check \kappa)\} \leq \Re\{\hat \lambda^{T,(2)}(\check \kappa)\} \leq \ldots 
\end{equation*} \\
Note that $\psi^{-, (n)}(0_-) \neq \psi^{+,(n)} (0_+) $.  The
associated eigenvalue is
\begin{equation}
\label{eq:lambda0T}
\lambda^{T,(n)}(\kappa)  = \hat \lambda^{T,(n)} \bigl(\kappa/|v_0|^\frac 13\bigr) \, |\, v_{01}|^{\frac23} \,.
\end{equation}
In what follows, $(\lambda_0^T(\kappa), \psi_0^T)$ denotes an
eigenpair $(\lambda^{T,(n)}(\kappa), \psi^{T,(n)})$ corresponding to a
particular choice of $n\geq 1$. \\

{\bf Summary at this stage.}
For $\# \in \{D,N,R,T\}$, we have constructed $u_0^\#$ in the form
\eqref{qm3}.  At this step $\phi_0^\# (\sigma)$ remains ``free''
except that it should not be identically $0\,$.  We have chosen
$\lambda_0^\#$ as an eigenvalue of $\mathcal L_0^\#$ (assuming that it
is simple, with no Jordan block) and $\psi_0^\#$ is the associated
eigenfunction of $\mathcal L_0^\#$, which belongs to $\mathcal S^\#\,$
and permits, according to Lemma~\ref{abstlemma}, to have the
normalization
\begin{equation}\label{normpsi0}
\int_{\mathbb R^\#}  \psi_0^\#(\tau) ^2 d\tau =1\, .
\end{equation}

From now on, we do not mention (except for explicit computations) the
reference to Dirichlet, Neumann, Robin or Transmission condition when
the construction is independent of the considered case.

\subsection{Term $j=1$}

The second equation (corresponding to $j=1$) reads
\begin{equation}\label{eq:j=1a}
 (\L_0^\# -\lambda_0)\, u_1^\#= \lambda_1\, u_0^\#\,.
\end{equation}
We omit sometimes the superscript $\#$ for simplicity.\\
The guess is that $\lambda_1=0\,$.  To see if it is a
necessary condition, one can take the scalar product (in the $\tau$
variable) with $\bar \psi_0$ (to be understood as the element in ${\rm
Ker} (\L_0^* -\bar \lambda_0)$).  We take the convention that the
scalar product is antilinear in the second argument.  This leads to
\begin{equation*}
\left(\int \psi_0^2(\tau) d\tau\right) \,  \lambda_1 \,  \phi_0(\sigma)=0\,,
\end{equation*}
the integral being on $\mathbb R^+$ for Dirichlet, Neumann or Robin,
and on $\mathbb R$ in the transmission case.  From
Eq. (\ref{normpsi0}), we get then
\begin{equation*}
 \lambda_1 \phi_0(\sigma)=0\,,
\end{equation*}
and by the previous condition on $\phi_0(\sigma)$
\begin{equation}
\label{eq:lambda1}
 \lambda_1=0\,.
\end{equation}
Hence, coming back to \eqref{eq:j=1a}, we choose
\begin{equation}
 u_1^\#(\sigma,\tau) = \phi_1^\#(\sigma) \psi_0^\# (\tau)\,,
\end{equation}
where $\phi_1^\#$ remains free at this step.\\

\subsection{Term $j=2\,$}

The third condition (corresponding to $j=2$) reads
\begin{equation}\label{a10}
 (\L_0^\# -\lambda_0) \, u_2 + \L_2 \, u_0 = \lambda_2 \, u_0\,.
\end{equation}
To find a necessary condition, we take the scalar product (in the
$\tau$ variable) with $\bar \psi_0$.  In this way we get (having in
mind \eqref{normpsi0})
\begin{equation*}
  \langle \L_2\, u_0\,,\, \bar \psi_0 \rangle = \lambda_2\,   \phi_0(\sigma)\,.
\end{equation*}
Computing the left hand side, we get
\begin{equation*}
  (-\partial_\sigma^2 +i\, v_{20}\, \sigma^2 ) \, \phi_0(\sigma) = \lambda_2\,  \phi_0(\sigma)\,.
\end{equation*}
From Assumption \eqref{condV}, we know that $v_{20} \neq 0\,$.  Hence
we are dealing with an effective complex harmonic oscillator whose
spectral analysis has been done in detail (see Davies \cite{Dav} or
the book by Helffer \cite{H}).  The eigenvalues can be explicitly
computed (by analytic dilation) and there is no Jordan block.
Moreover the system of corresponding eigenfunctions is complete.  This
implies that $(\lambda_2,\phi_0)$ should be a spectral pair for $
(-\partial_\sigma^2 + i \, v_{20} \, \sigma^2 )\,$.\\
The eigenpairs of the quantum harmonic oscillator are well known:
\begin{equation}
\label{eq:phi0_k}
\lambda_2^{(k)} = \gamma (2k-1),  \quad 
\phi_0^{(k)}(\sigma) = \frac{\gamma^{\frac 14} \, e^{-\gamma \sigma^2/2}\,  H_{k-1}(\gamma^{\frac 12} \sigma)}{\pi^{\frac 14} \sqrt{ 2^{k-1}\, (k-1)! }}\,  
  \quad (k = 1,2,\ldots),
\end{equation}
where $\gamma^s = |v_{20}|^{\frac s 2} \exp\bigl(\frac{\pi s i}{4} \,
\sign v_{20}\bigr)$ (for $s=\frac 14, \frac 12, 1$), $H_k(z)$ are
Hermite polynomials, and the prefactor ensures that
\begin{equation*}
\int\limits_{-\infty}^\infty \phi_0^{(k)}(\sigma)^2\,  d\sigma =1\,.
\end{equation*}
The eigenvalue with the smallest real part corresponds to $k = 1$ for
which
\begin{equation}
\label{eq:phi0} 
\phi_0^{(1)}(\sigma) = c_{\phi_0} \exp\left( - \lambda_2 \frac{\sigma^2}{2} \right)\,,
\end{equation}
while the corresponding eigenvalue is
\begin{equation}
\label{eq:lambda2}
\lambda_2^{(1)}= |v_{20}|^\frac 12 \exp \left(\frac{i\pi}{4}\, {\rm sign}~ v_{20}\right) \,,
\end{equation}
and $c_{\phi_0}$ ensures the normalization of $\phi_0^{(1)}(\sigma)$:
\begin{equation}
\label{eq:Cphi0} 
c_{\phi_0} =   |v_{20}|^\frac 18 \,\pi ^{-\frac 14} \exp \left(\frac{i\pi}{16}\, {\rm sign}~ v_{20}\right)\, .
\end{equation}
We do not need actually the specific expression of $\phi_0^\#=\phi_0$
and it is enough to know that $\phi_0^\# \in \mathcal S(\mathbb R)$.\\

Coming back to the solution of \eqref{a10}, which simply reads
\begin{equation}\label{a10a}
 (\L_0 -\lambda_0) \, u_2 = 0 \,,
\end{equation}
we consequently  look for $u_2^\# (\sigma,\tau)$ in the form
\begin{equation}\label{formphi2}
  u_2^\#(\sigma,\tau) = \phi_2^\#(\sigma) \, \psi_0 ^\# (\tau) \,,
\end{equation}
where $\phi_2^\#(\sigma)$ is free at this stage.\\

{\bf Summary at this stage.}
We note that the construction is conform with the general form
introduced in \eqref{qm3}.  At this stage, $(\lambda_0^\#,\psi_0^\#)$
is a spectral pair for $\mathcal L_0^\#\,$, $\lambda_1^\#=0\,$,
$u_1^\#(\sigma, \tau)=\phi_1^\#(\sigma) \psi_0^\#(\tau)$ (with
$\phi_1^\#$ free), $(\lambda_2^\#,\phi_0^\#)$ is a spectral pair for
$\mathcal L_2$ (actually independent of $\#$).  Note that $\phi_0^\#$
can be either odd, or even.

\subsection{ Term $j=3\,$}

The fourth equation corresponds to $j=3$ and reads
\begin{equation}\label{a12}
 (\L_0 -\lambda_0) u_3 + (\L_2-\lambda_2) u_1+ \L_3 u_0 = \lambda_3 u_0\,.
\end{equation}
Taking the scalar product (in $L^2\hat \otimes L^2_\#: = L^2(\mathbb
R_\sigma \times \mathbb R_\tau^+)$ for Dirichlet, Neumann and Robin,
and in $L^2\hat \otimes L^2_\# := L^2(\mathbb R_\sigma \times \mathbb
R_\tau^-) \times L^2(\mathbb R_\sigma \times \mathbb R_\tau^+) $ for
the transmission case) with $\bar u_0$ and having in mind our
normalizations of $\psi_0$ and $\phi_0$, we obtain
\begin{equation*}
  \langle \L_3 u_0\,,\, \bar u_0\ \rangle = \lambda_3  \,,
\end{equation*}
so $\lambda_3$ is determined by
\begin{equation}
\label{eq:lambda3}
  \lambda_3 = i  \, v_{11}  \left(\int \sigma \phi_0(\sigma)^2 d\sigma \right) \left(\int \tau \psi_0^\#(\tau)^2 d\tau \right) \,.
\end{equation}
Note that whatever the parity of $\phi_0$, $\phi_0^2$ is even, so
$\int \sigma \phi_0(\sigma)^2 d\sigma =0\,$.  Hence,
\begin{equation}
\lambda_3=0\,.
\end{equation}
We come back to \eqref{a12}, but now take the scalar product with
$\bar \psi_0$ in the $\tau$ variable.  So we get
\begin{equation*}
  \langle (\L_2-\lambda_2) u_1+ (\L_3 -\lambda_3)  u_0\,,\, \bar \psi_0 \rangle = 0\,.
\end{equation*}
Taking into account \eqref{condV} and the form of $u_0$ and $u_1$,
this reads
\begin{equation}\label{a18}
  (\L_2-\lambda_2)\, \phi_1 =  -  i \, v_{11} \, \sigma\,   \left(\int  \tau  \psi_0(\tau)^2 d\tau\right)\, \phi_0 \,.
\end{equation}
The right hand side is in the image of the realization of
$(\L_2-\lambda_2)\,$.  There is a unique $\phi_1$ solution of
\eqref{a18} satisfying
\begin{equation}
\int_\mathbb R \phi_1(\sigma) \phi_0(\sigma)\, d\sigma =0\,.
\end{equation}
 
\begin{remark}
Note that $\phi_0 \phi_1$ is odd.
\end{remark}

We can now solve \eqref{a12}.  We observe that
\begin{equation*}
(\L_2-\lambda_2) u_1+ ( \L_3-\lambda_3)  u_0 = ( (\L_2-\lambda_2) \phi_1) \psi_0 +  (\L_3-\lambda_3) u_0\, .
\end{equation*}
According to what we have done already, \eqref{a12}  has the form
\begin{equation*}
 \left( (\L_0-\lambda_0) u_3\right)(\sigma,\tau) =  g_3 (\tau)   f_3(\sigma)\,,
\end{equation*}
where 
\begin{equation*}
g_3(\tau) = (\tau  - c_3) \, \psi_0(\tau)   
\end{equation*}
is orthogonal to $\bar \psi_0$, i.e. 
\begin{equation*}
c_3= \int\, \tau \psi_0(\tau)^2\, d\tau\,,
\end{equation*}
and
\begin{equation*}
f_3(\sigma) = i \,v_{11\,} \sigma   \phi_0(\sigma)\,.
\end{equation*}
 
\begin{remark}
Note that $\phi_0 f_3$ is odd.
\end{remark}

We then write for $j=3$ the expression \eqref{qm4}, with $N_3=1$,
\begin{equation}
  u_3 (\sigma, \tau) = \phi_{3}(\sigma) \psi_0(\tau) + \phi_{3,1}(\sigma) \psi_{3,1} (\tau)\,,
\end{equation}
where $\psi_{3,1}$ is determined as the unique solution of the problem
\begin{equation}\label{5.52}
  (\L_0^\# -\lambda_0^\#)\, \psi_{3,1}= g_3 \, ,
\end{equation}
 which is orthogonal to $\bar \psi_0$, and
\begin{equation}
\phi_{3,1}(\sigma) = f_3(\sigma)\,.
\end{equation}
\begin{remark}
Note that $\phi_0 \phi_{3,1}$ is odd.
\end{remark}

{\bf Summary at this stage.}
We note that the construction is conform with the general form
introduced in \eqref{qm3}-\eqref{qm4}.  At this stage, $\phi_3^\#$ is
introduced, $\lambda_3^\# = 0$ and $\phi_1^\#$ are determined but
$\phi_2^\#$ and $\phi_3^\#$ remain free.  Note that $N_3=1$ in
\eqref{qm4}, $\phi_{3,1}^\#$ is determined in $\mathcal S (\mathbb R)$
and $\psi_{3,1}^\# $ is determined in $\mathcal S^\#$.

\subsection{Term $j=4\,$}

The fifth condition corresponds to $j=4\,$ and reads
\begin{equation}\label{a51}
 (\L_0-\lambda_0) u_4 + (\L_2-\lambda_2) u_2 + (\L_3 -\lambda_3)  u_1+ \L_4 u_0 = \lambda_4 u_0\,.
\end{equation}
We follow the same procedure as in the preceding step.  $\lambda_4$ is
determined by integrating \eqref{a51} after multiplication by $u_0\,$:
\begin{equation}
\label{eq:lambda4_def}
\begin{split}
\lambda_4 &= \langle ( \L_3 -\lambda_3) u_1+ \mathcal L_4 u_0\,,\,\bar u_0\rangle \\
& = i \, v_{11}\, \left( \int \sigma \phi_1(\sigma)\phi_0(\sigma) d \sigma \right)\, \left(\int \tau \psi_0(\tau)^2 d\tau \right)  \\  
&  + \curv  (0)  \int \psi_0'(\tau)\psi_0 (\tau) d \tau  + i\, v_{02}\, \int \tau^2 \psi_0(\tau)^2\, d\tau\,.  \\
\end{split}
\end{equation} 
$\phi_2$ is determined by integrating \eqref{a51} in the $\tau$
variable over $\mathbb R^\#$ after multiplication by $ \psi_0 \,$.  We
get
\begin{equation}\label{abc1}
(\mathcal L_2-\lambda_2) \phi_2 = - \langle (\L_3-\lambda_3)  u_1\,,\, \bar \psi_0\rangle_{L^2_\tau} - \langle \L_4 u_0\,,
\, \bar \psi_0\rangle_{L^2_\tau} + \lambda_4 \phi_0 := f_4 \,,
\end{equation}
where our choice of $\lambda_4$ implies the orthogonality of $f_4$ to
$\bar \phi_0$\, in $L^2_\#\,$.\\
There exists consequently a unique $\phi_2$ solution of \eqref{abc1}
that is orthogonal to $\bar \phi_0\,$.\\

We then proceed like in the fourth step, observing that $u_4$ should
satisfy, for some $N_4 \geq 1$,
\begin{equation}\label{a5a1}
 (\L_0-\lambda_0) u_4  = \sum_{\ell =1}^{N_4}  f_{4,\ell} (\sigma)\, g_{4,\ell}(\tau)\,,
\end{equation}
with $f_{4,\ell}$ in $\mathcal S(\mathbb R)$, $g_{4,\ell}$ in
$\mathcal S^\#$ and orthogonal to $\bar \psi_0$ in $L^2_\#$.  The
expression in the right hand side is deduced from our previous
computations of $u_0$, $u_2$ and $u_3$ and $\lambda_4$~. \\
We then look for a solution $u_4$ in the form
\begin{equation}
 u_4(\sigma,\tau) = \phi_4(\sigma)\, \psi_0(\tau) + \sum_{\ell=1}^{N_4} \phi_{4,\ell} (\sigma) \, \psi_{4,\ell}(\tau)\,,
\end{equation}
which is obtained by solving for each $\ell$
\begin{equation}\label{5.59}
(\L_0^\#- \lambda_0^\# ) \psi_{4,\ell}= g_{4,\ell}\,,\quad \int \psi_{4,\ell} (\tau)\, \psi_0(\tau)\,d \tau =0\,,
\end{equation}
with the suitable boundary (or transmission) condition at $0$ and
taking
\begin{equation*}
\phi_{4,\ell}   =  f_{4,\ell} \,.
\end{equation*}
\\
We now make explicit the computation of the right hand side in
\eqref{a5a1}.  Using our choice of $\lambda_4$ and $\phi_2$, we obtain
\begin{equation*}
\begin{split}
& -(\L_2-\lambda_2) u_2 -  (\L_3 -\lambda_3)  u_1-  \L_4 u_0 + \lambda_4 u_0\\
& = (-(\L_2-\lambda_2) \phi_2) \psi_0 - ( (\L_3 -\lambda_3)  \phi_1 \psi_0  - \L_4 \phi_0 \psi_0 + \lambda_4 \phi_0 \psi_0\\
& = f_{4,1} (\sigma)  (\tau -c_3)  \psi_0 (\tau)  +  f_{4,2} (\sigma) (\partial_ \tau \psi_0 - c_4\psi_0)  +  f_{4,3} (\sigma) (\tau^2 -c_5) \psi_0  \,, \\
\end{split}
\end{equation*}
with $c_4 =\int (\partial_ \tau \psi_0)(\tau) \psi_0 (\tau) d\tau$ and
$c_5 = \int \tau^2 \psi_0 (\tau) d\tau$\,.\\
Moreover the $f_{4,\ell}$ are even with respect with $\sigma$.\\
Hence we can take $N_4=3$ and 
\begin{equation}\label{a57}
\begin{split} 
g_{4,1} (\tau) & := (\tau -c_3)  \psi_0 (\tau) \,, \\ 
g_{4,2} (\tau) & := (\partial_ \tau \psi_0 - c_4 \psi_0)\,, \\
g_{4,3} (\tau) & := (\tau^2 -c_5)  \psi_0 (\tau)\,.  \\
\end{split}
\end{equation}
We do not provide explicit formula for the corresponding $\psi_{4,\ell}\,$. \\

{\bf Summary at this stage.}
At the end of this step we have determined the $\lambda_j^\#$ for
$j\leq 4\,$, the $\psi^\#_{j,\ell}$ and $\phi^\#_{4,\ell}$ for $3 \leq
j \leq 4$ and the $\phi^\#_j(\sigma)$ for $j\leq 2\,$.  Like in
\cite{HK}, this construction can be continued to any order.  This
achieves the proof of Theorem \ref{th:app}.

\subsection{Term $j=5\,$ and vanishing of the odd terms}

We first focus on the sixth step corresponding to the computation of
$\lambda_5$.  The sixth condition corresponds to $j=5\,$ and reads
\begin{equation}\label{a515}
 (\L_0-\lambda_0) u_5 + (\L_2-\lambda_2) u_3 + (\L_3 -\lambda_3)  u_2+ ( \L_4-\lambda_4)  u_1  + \L_5 u_0 = \lambda_5 u_0\,.
\end{equation}
$\lambda_5$ is determined by integrating \eqref{a515} after
multiplication by $\bar u_0\,$.  By our preceding constructions and
\eqref{sym}, we see that
\begin{equation*}
\sigma \mapsto  u_0(\sigma)\,  \bigl( (\L_2-\lambda_2) \,u_3 + (\L_3 -\lambda_3) \, u_2+ ( \L_4-\lambda_4) \, u_1  + \L_5 \, u_0\bigr) (\sigma)
\end{equation*}
is odd.  This immediately leads to $\lambda_5=0$\,.\\
 
With some extra work consisting in examining the symmetry properties
with respect to $\sigma$ and using Sec. \ref{sssparity}, we obtain
\begin{proposition}
In the formal expansion, $\lambda_j=0$ if $j$ is odd.
\end{proposition}

\subsection{Four-terms asymptotics}

Gathering \eqref{eq:Ah_series}, \eqref{eq:lambda0N},
\eqref{eq:lambda0D}, \eqref{eq:lambda0R}, \eqref{eq:phi0_k}, and
\eqref{eq:lambda4_def}, the four-terms asymptotics of approximate
eigenvalues reads for $n,k = 1,2,\ldots$
\begin{equation}  \label{eq:lambda_asympt}
\begin{split}
\lambda^{app,\#}_h := \lambda^{\#,(n,k)}_h 
& = i \, v_{00} + h^{\frac 23} |\, v_{01}|^{\frac 23} \mu_n^\# \exp\left(\frac{i\pi}{3} \sign \, v_{01}\right) \\
& + h (2k-1) |v_{20}|^{\frac 12} \exp\left(\frac{i\pi}{4} \sign
v_{20}\right) + h^{\frac 43} \lambda_4^{\#,(n)} + \mathcal O(h^{\frac 53}) \,, \\
\end{split}
\end{equation}
where $\mu_n^D = -a_n$, $\mu_n^N = -a'_n$, $\mu_n^R = -a_n^R(\kappa)$
(defined by (\ref{eigenvalueR})), and $\mu_n^T = -a_n^+(\kappa)$
(defined by (\ref{eq:ankappa})), while $\lambda_4^{\#,(n)}$ is
explicitly computed in Appendix \ref{sec:lambda4} (see
\eqref{eq:lambda4D}, \eqref{eq:lambda4N},
\eqref{eq:lambda4R}, and \eqref{eq:lambda4T} for Dirichlet, Neumann,
Robin, and Transmission cases), and the involved coefficients $v_{jk}$
of the potential $V(s,\t)$ are defined in (\ref{eq:Vjk}). \\

\begin{remark}
In the above construction, if we take $\phi_j^\# =0$ for $j\geq 3$, we
get an eigenpair $(\lambda^{app,\#}_h,u^{app,\#}_h)$ with
\begin{equation*}
u^{app,\#}_h =u_0^\#  + h^\frac 16 u_1^\# + h^\frac 13 u_2^\#\,,
\end{equation*}
such that
\begin{equation}\label{asymptaaa}
(\mathcal A_h^\# -\lambda^{app,\#}_h) u^{app,\#}_h = \mathcal O (h^{\frac 32})\,.
\end{equation}
To get in \eqref{asymptaaa} the remainder $\mathcal O (h^\frac 53)$,
one should continue the construction for two more steps.
\end{remark}

\begin{remark}
Note that the leading terms in the eigenvalue expansion do not contain
the curvature which appears only in $\lambda_4$ (see
Eq. (\ref{eq:lambda4N})) and is thus of order $h^{\frac 43}$.
\end{remark}

\section{Other scalings in the Robin or transmission problems}
\label{sec:kappa_scaling}

The scaling (\ref{assT}) of the transmission parameter $\Kappa$ with
$h$ was appropriate to keep the Robin or transmission condition for
the rescaled problem.  In biophysical applications, the transmission
condition reads
\begin{equation}
\label{eq:assT3}
D \partial_\nu u_+ =D \partial_\nu u_- = \Kappa \, \bigl(u_+ - u_- \bigr),
\end{equation}
where $D$ is the bulk diffusion coefficient, while the transmission
parameter $\Kappa$ represents the permeability of a membrane which is
set by the membrane properties and thus does not necessarily scale
with $h$.  Similarly, in the Robin boundary condition,
\begin{equation}
\label{eq:assT3a}
- D \partial_\nu u_- =\Kappa \, u_- ,
\end{equation}
which accounts for partial reflections on the boundary, $\Kappa$
represents partial reactivity or surface relaxivity which are set by
properties of the boundary. \\

We consider two practically relevant situations for the BT-operator \break  
$-D \Delta + i\,g\,x_1$:
\begin{itemize}
\item
When $D\to 0$ with fixed $g$, one can identify $h^2 = D$ and $V(x) =
g\,x_1$ so that the rescaled transmission condition in (\ref{assT0})
gives $\Kappa h^{-\frac43}$ which tends to $+\infty$ as $h\to 0$ if
$\Kappa$ is fixed.  In this limit, the transmission condition is
formally reduced to the continuity condition at the boundary:
$u_+(\sigma,0) = u_-(\sigma,0)$, together with the flux continuity in
the first relation of (\ref{assT0}).  In other words, the interface
between two subdomains is removed.  The construction of the previous
section seems difficult to control in this asymptotics and the
mathematical proof of the heuristics should follow other ways.

\item
When $g\to \pm \infty$ with fixed $D$, one can divide the BT-operator
and (\ref{eq:assT3}) by $g$ and then identify $h^2 = D/g$ and $V(x) =
x_1$.  In this situation, the rescaled transmission condition in
(\ref{assT0}) gives a parameter $\kappa= (\Kappa/D) h^{\frac23}$ which
tends to $0$ as $h\to 0$.  In this limit, the transmission condition
is reduced to two Neumann boundary conditions on both sides of the
interface: $\partial_\tau u_+(\sigma,0) = \partial_\tau u_-(\sigma,0)
= 0\,$.
\end{itemize}

We now discuss how the eigenvalue asymptotic expansion obtained for
rescaled $\Kappa$ can be modified for the second situation.  The
constructions of the previous section can be adapted and controlled
with respect to $\kappa$ for $\kappa $ small enough.  As observed
along the construction,  one can start with
(\ref{eq:lambda_asympt}) and then expand the factor $\mu^\#_n(\kappa)$
into Taylor series that results in the quasi-mode in the Robin or
Transmission case:
\begin{theorem}  \label{thm:scaling}
With the notation of Theorem \ref{th:app} except that in
\eqref{scaKappa} we assume
\begin{equation}  \label{eq:kappa_newscaling}
\kappa= \hat \kappa \, h^\frac 23\,,
\end{equation}
we have for $\# \in \{R,T\}$, $n,k = 1,2,\ldots$
\begin{equation}
\label{newscaling}
\begin{split}
\lambda^{\#,(n,k)}_h & = i \, v_{00} - h^{\frac 23} |\, v_{01}|^{\frac 23} a'_n  \exp\left(\frac{i\pi}{3} \sign \, v_{01}\right) \\
& + h\, (2k-1)\, |v_{20}|^{\frac 12} \exp\left(\frac{i\pi}{4} \sign v_{20}\right)   \\
& + h^{\frac 43} \left( \lambda_4^{N,(n)} - \hat \kappa\,  \frac{|v_{01}|^{\frac13}}{a'_n} 
\exp\left(\frac{i\pi}{6}\, \sign v_{01}\right) \right) + O(h^{\frac 53})\,,  \\
\end{split}
\end{equation}
where $\lambda_4^{N,(n)}$ is explicitly given in \eqref{eq:lambda4N},
and the involved coefficients $v_{jk}$ of the potential $V(s,\t)$ are
defined in (\ref{eq:Vjk}). 
\end{theorem}
Here, we have used that $\lambda_4^{\#,(n)}(\kappa)=\lambda_4^{N,(n)}$
for $\kappa=0$ (see Remark \ref{RemRkappa=0}).  The coefficient in
front of $\hat\kappa$ involves $(\mu_n^\#)'(0)$ that was computed
explicitly by differentiating the relation determining
$\mu_n^\#(\kappa)$ with respect to $\kappa$.  For the Robin case, we
used (\ref{eq:anRprime}) to get
\begin{equation}
\label{eq:muRprime}
(\mu_n^R)'(0) = - (a_n^R)'(0) = - \frac{1}{a'_n \, |v_{01}|^{\frac13}} \exp\left(-\frac{i\pi}{6} \sign \, v_{01}\right) ,
\end{equation}
with $a_n^R(0) = a'_n$. \\
Similarly, differentiating (\ref{eq:an_equation}) with respect to
$\kappa$ and using (\ref{eq:ankappa}), we got (see Appendix
\ref{sec:muT})
\begin{equation} 
\label{eq:muTprime}
(\mu_n^T)'(0) = - (a_n^+)'(0) = - \frac{1}{a'_n \, |v_{01}|^{\frac13}} \, \exp\left(-\frac{i\pi}{6}\, \sign v_{01}\right) ,
\end{equation}
with $a_n^+(0) = a'_n$.  The effect of Robin or transmission condition
appears only in the coefficient of $h^\frac 43$.

In order to control the construction with respect to $\kappa$, it is
enough to get an expression of the kernel of the regularized resolvent
for $z=\lambda^\#_0$.  Let us treat the Robin case and assume $v_{01}
=-1$.\\
As proven in \cite{GHH}, the kernel of the resolvent is given by
\begin{equation*}
\mathcal G^{-,R} (x,y\,;\lambda) =  \mathcal G _0 ^-(x,y\,;\lambda) + \mathcal G _1^{-,R}(x,y\,;\kappa, \lambda)\,  \quad \textrm{for}~ (x,y)\in \mathbb R_+^2,
\end{equation*}
where 
\begin{equation}
\mathcal G _0 ^-(x,y\,;\lambda) = \begin{cases} 2\pi \Ai(e^{i\alpha} w_x) \Ai(e^{-i\alpha} w_y)  \quad (x<y), \cr
2\pi \Ai(e^{-i\alpha} w_x) \Ai(e^{i\alpha} w_y)  \quad (x>y), \end{cases}
\end{equation}
and
\begin{equation}
\begin{split}
\mathcal G ^{-,R}_1(x,y\,;\kappa, \lambda) & = - 2\pi \frac{i e^{i\alpha} \Ai'(e^{i\alpha} \lambda) - \kappa \Ai(e^{i\alpha}\lambda)}
{i e^{-i\alpha} \Ai'(e^{-i\alpha} \lambda) - \kappa \Ai(e^{-i\alpha} \lambda)}  \\
& \times  \Ai\bigl(e^{-i\alpha} (ix+\lambda)\bigr)~ \Ai\bigl(e^{-i\alpha} (iy+\lambda)\bigr)\,.   
\end{split}
\end{equation}
The kernel $ \mathcal G _0 ^-(x,y\,;\lambda)$ is holomorphic in
$\lambda$ and independent of $\kappa$.  Setting $\kappa = 0\,$, one
retrieves the resolvent for the Neumann case.  Its poles are
determined as (complex-valued) solutions of the equation
\begin{equation}\label{eigenvalueR1}
f^R(\kappa, \lambda):=i e^{-i\alpha} \Ai'(e^{-i\alpha} \lambda) - \kappa \Ai(e^{-i\alpha} \lambda) = 0\,.
\end{equation} 
For $\kappa=0$, we recover the equation determining the poles of the
Neumann problem: 
\begin{equation*}
f^N(\lambda):=i e^{-i\alpha} \Ai'(e^{-i\alpha} \lambda)=0\,.
\end{equation*}
We look at the first pole and observe that
\begin{equation}\label{eigenvalueR2}
(\partial_\lambda f^R)(0, \lambda^{R,(1)}(0)) = (\partial_\lambda f^R)(0, \lambda^{N,(1)}) = (f^N)'(\lambda^{N,(1)}) \neq 0\,.
\end{equation}

This evidently remains true for $\kappa$ small enough:
\begin{equation}
(\partial_\lambda f^R)(\kappa, \lambda^{R,(1)}(\kappa))\neq 0\,.
\end{equation}
As done in \cite{GHH}, we can compute the distribution kernel of the
projector associated with
\begin{equation*}
\lambda_0 (\kappa) := \lambda^{R,(1)}(\kappa)\,.
\end{equation*}
We get
\begin{equation}
\begin{split}
\Pi^R_1 (x,y;\kappa)  & = - 2\pi \frac{i e^{i\alpha} \Ai'(e^{i\alpha} \lambda_0(\kappa)) 
- \kappa \Ai(e^{i\alpha}\lambda_0(\kappa) )} {(\partial_\lambda f^R)(\kappa, \lambda_0(\kappa))}  \\
       & \times  \Ai\bigl(e^{-i\alpha} (ix+\lambda_0(\kappa))\bigr)~ \Ai\bigl(e^{-i\alpha} (iy+\lambda_0(\kappa))\bigr)\,.   
\end{split}
\end{equation}
This kernel is regular with respect to $\kappa$.\\
The distribution kernel of the regularized resolvent at $\lambda_0
(\kappa)$ is obtained as
\begin{equation*}
\begin{split}
\mathcal G^{R,reg} (x,y;\kappa,  \lambda_0 (\kappa)) &:=  \mathcal G _0 ^-(x,y;\kappa, \lambda_0(\kappa) ) \\
& + \lim_{\lambda \rightarrow \lambda_0} \left(  \mathcal G _1 ^{-,R} (x,y\,;\kappa, \lambda) - (\lambda_0-\lambda)^{-1} \Pi^R_1 (x,y;\kappa)\right) \,. \\
\end{split}
\end{equation*}
It remains to compute the second term of the right hand side.  Writing
$\mathcal G_1^{-,R}(x,y;\kappa,\lambda)$ in the form
\begin{equation*}
\mathcal G_1^{-,R}(x,y,\kappa,\lambda)= \frac{\Phi(x,y;\kappa,\lambda)}{\lambda -\lambda_0 (\kappa)} ,
\end{equation*}
we observe that $\Phi(x,y;\kappa,\lambda)$ is regular in $\kappa,
\lambda$ and we get
\begin{equation*}
\mathcal G^{R,reg} (x,y;\kappa, \lambda_0 (\kappa)):=  \mathcal G _0 ^-(x,y;\kappa, \lambda_0(\kappa) ) + \partial_\lambda \Phi (x,y;\kappa, \lambda_0(\kappa)) \,.
\end{equation*}
It is regular in $\kappa$ and we recover for $\kappa=0$ the
regularized resolvent of the Neumann problem at $\lambda=\lambda^{N,
(1)}$.\\

With this regularity with respect to $\kappa$, we can control all the
constructions for $j=0,\dots, 4$ (and actually any $j$) and in
particular solve \eqref{5.52} for $\kappa$ small and similarly
\eqref{5.59}, with a complete expansion in powers of $\kappa$ at the
origin. \\
 
\begin{remark}
Similarly, one can treat the transmission case.
\end{remark}

\section{WKB construction}\label{sec.WKB}

In this section, we propose an alternative analysis based on the  WKB
method.  This construction is restricted to quasimodes with $k = 1$ in
\eqref{eq:phi0_k} but it gives a quasimode state that is closer to the
eigenfunction than that obtained by the earlier perturbative approach.
Here we follow the constructions of \cite{HK,HKR} developed for a
Robin problem.

We start from
\begin{equation}\label{a1ba}
 \A_h=- h^2 a^{-2}\partial_s^2  +  h^2 a^{-3}(\partial_s a)\, \partial_s  - h^2 \partial_\t^2 
- h^2a^{-1}(\partial_\t a) \,\partial_\t  + i \,  \widetilde  V(s,\t) \,.
\end{equation}
Here, instead of what was done in \eqref{dil1}, we only dilate in the
$\rho$ variable:
\begin{equation*}
\rho =h^\frac 23 \tau\,.
\end{equation*}
In the $(s,\tau)$ coordinates, we get
\begin{equation}\label{a6a}
 \widehat{\A}_h=- h^2 \check a_h^{-2}\partial_s^2  +  h^2 \check a_h^{-3} {(\partial_s \check a_h)} \, \partial_s
-h^\frac 23 \partial_\tau^2 - h^\frac 43  \check a_h^{-1} \check{(\partial_\t  a)} \,\partial_\tau  + i\,  \check V_h(s, \tau)\,,
\end{equation}
with
\begin{equation}
\begin{split}
\check V_h(s, \tau) &=\widetilde  V(s,h^\frac 23 \tau)\, ,\\ \check a_h (s,\tau) & = 1 - \tau h^\frac 23  \curv\, (s)\,,\\
\partial_s  \check a_h (s,\tau) & =  - \tau h^\frac 23  \curv'\, (s)\, , \\
  \check{\partial_\t  a} & = -\curv (s)\,,\\
 \check a_h (s,\tau)^2 & = 1 - 2 \tau h^\frac 23  \curv\, (s)+ \tau^2  h^\frac 43  \curv\, (s)^2 \,, \\
 \check a_h (s,\tau)^{-2}& =  1 + 2 \tau h^\frac 23  \curv\, (s) + 3  \tau^2  h^\frac 43  \curv\, (s)^2 + \mathcal O (h^2)\,. \\
\end{split}
\end{equation}

We consider the Taylor expansion of $\check V_h$:
\begin{equation}
\check V_h(s, \tau) \sim  \sum_{j\in \mathbb N} v_j(s) h^\frac{ 2j}{3}\tau^j \, ,
\end{equation}
with
\begin{equation}  \label{eq:vj_WKB}
v_j(s)= \frac{1}{j!} (\partial_\t^j \widetilde V) (s,0)\,.
\end{equation}
We look for a trial state in the form 
\begin{equation} \label{defwkb}  
u^{\#,wkb}_h:=  d(h) b_h(s,\tau) \exp \left(-\frac{\theta(s,h)}{h}\right)\,,
\end{equation}
with
\begin{equation}\label{deftheta}
\theta (s,h)= \theta_0 (s) + h^\frac23 \theta_1(s) \,,
\end{equation}
and
\begin{equation}\label{defb}
b_h  (s,\tau) \sim \sum_{j\in \mathbb N} b_j (s,\tau) h^\frac j3\,.
\end{equation}
Here $d(h)$ is a normalizing constant such that, when coming back to
the initial coordinates, the $L^2$ norm of $u^{\#,wkb}_h$ is $1$.  In
the initial coordinates, we should actually consider $u^{\#,wkb}_h (s,
h^{-\frac 23} \rho)$ multiplied by a suitable cut-off function in the
neighborhood of the point $x^0$ of $\partial \Omega^\perp$. \\
This gives an operator acting on $b_h$
\begin{equation}\label{a6ab}
\begin{split}
\widehat{\A}_{h,\theta}&:= \exp \left(\frac{\theta(s,h)}{h}\right) \, \widehat{\A}_h \, \exp \left(-\frac{\theta(s,h)}{h}\right)  \\
& =- \check a_h^{-2}(h\partial_s-\theta'(s,h)) ^2  +   h \check a_h^{-3} { (\partial_s \check a_h)} \, (h\partial_s-\theta'(s,h)) \\
&  -h^\frac 23 \partial_\tau^2   - h^\frac 43  \check a_h^{-1} \check{(\partial_\t  a)} \,\partial_\tau  + i\,  \check V_h(s, \tau)\,. \\
\end{split}
\end{equation}
We rewrite this operator in the form 
\begin{equation}\label{eq:exp} 
\widehat { \mathcal A}_{h,\theta} \sim \sum_{j\geq 0} \Lambda_j h^\frac j 3\,,
\end{equation}
with
\begin{equation}
\begin{split}
\Lambda_0 &:= i v_0(s) -\theta_0'(s)^2\,,\\
\Lambda_1 &:=  0 \,, \\
\Lambda_2 &:= -\partial_\tau^2 + (i v_1(s)  -2\mathfrak c (s) \theta_0'(s)^2)  \tau  - 2 \theta'_0 (s) \theta'_1 (s) \,, \\
\Lambda_3 &:= 2 \theta'_0(s) \partial_s  + \theta''_0(s) \,, \\
\Lambda_4 &:= \mathfrak c (s) \partial_\tau   + \bigl( iv_2(s)  -  3 \mathfrak c (s)^2 \theta_0'(s)^2  \bigr)\tau^2 
+ 4  \mathfrak c (s)^2 \theta_0'(s)\theta'_1(s) \tau    - \theta_1'(s)^2    \,. \\
\end{split}
\end{equation}
We recall that  $v_0'(0)=0$, $v_1(0) \neq 0\,$.\\
We look for a quasimode in the form
\begin{equation}
\lambda_h^{\#,wkb} \sim i v_0(0) + h^\frac 23 \sum_{j\in \mathbb N}  \mu_j h^\frac j3 \,.
\end{equation}
The construction should be local in the $s$-variable near $0$ and
global in the $\tau$ variable in $\mathbb R^\#$.\\

Expanding $(\widehat{\mathcal A}_{h,\theta} -\lambda_h) b_h$ in powers of
$h^\frac 13$ and looking at the coefficient in front of $h^0$, we get
\begin{equation*}
 (\Lambda_0 - i v_0 (0)) b_0 =0\,
\end{equation*}
as a necessary condition.  Hence we choose $\theta_0$ as a solution of
\begin{equation}\label{E0}
i (v_0(s)-v_0(0)) - \theta_0'(s) ^2  =0\,,
\end{equation}
which is usually called the (first) eikonal equation.\\
We take the solution such that 
\begin{equation} \label{E1}
\Re\, \theta_0 (s) \geq 0\,,\quad \theta_0(0) =0\,,
\end{equation}
and we note that
\begin{equation}
\theta'_0(0)=0 \quad \textrm{and} \quad \theta''_0(0) \ne 0\,.
\end{equation}
With this choice of $\theta_0$, we note that
\begin{equation}
\Lambda_2 = -\partial_\tau^2 + i \bigl(v_1(s)  -2\mathfrak c (s) [v_0(s)-v_0(0)]\bigr)  \tau  - 2 \theta'_0(s) \theta'_1(s) \,, \\
\end{equation}
with 
\begin{equation}  \label{eq:v1hat_WKB}
\hat v_1(s) := v_1(s)  -2\mathfrak c (s) [v_0(s)-v_0(0)]
\end{equation}
being real.\\
As operator on $L^2_\#$, with the corresponding boundary or
transmission condition $\#\in \{D,N,R,T\}$, it satisfies
\begin{equation*}
\Lambda_2^{\#,*}= \overline{\Lambda_2^\#}\,.
\end{equation*}
The coefficient in front of $h^\frac 13$ vanishes and we continue with
imposing the cancellation of the coefficient in front of $h^\frac 23$
which reads
\begin{equation*}
 (\Lambda_0 - i v_0 (0)) b_2 + \Lambda_2 b_0 = \mu_0 b_0\,,
\end{equation*}
or, taking account of our choice of $\theta_0$,
\begin{equation}\label{E3}
- 2 \theta'_0 (s) \theta'_1(s) b_0(s,\tau) + (- \partial_\tau^2 + i \hat v_1  (s) \tau) b_0 (s,\tau)   - \mu_0 b_0(s,\tau) =0\,.
\end{equation}
Considering this equation at $s=0$, we get as a necessary condition
\begin{equation}\label{E4}
(- \partial_\tau^2 + i v_1  (0) \tau) \, b_0(0,\tau) = \mu_0 \, b_0(0,\tau)\,.
\end{equation}
If we impose a choice such that $b_0(0,\tau)$ is not identically $0$,
we get that $\mu_0$ should be an eigenvalue of (the suitable
realization of) $- \partial_\tau^2 + i v_1 (0) \tau$, i.e. $\mathcal
L_0^\#$.  We take some simple eigenvalue $\mu_0$ and define $\mu_0(s)$
as the eigenvalue of the operator
\begin{equation}  \label{eq:WKB_mu0s}
 - \partial_\tau^2 + i \hat v_1 (s) \tau
\end{equation}
such that $\mu_0(0) =\mu_0$.  If $f_0(s,\tau)$ denotes the
corresponding eigenfunction normalized as
\begin{equation}\label{E4a} 
\int f_0(s,\tau)^2 d\tau =1\,,
\end{equation} 
we can look for 
\begin{equation}\label{E5} 
b_0(s,\tau) = c_0 (s) f_0(s,\tau)\,.  
\end{equation} 
We now come back to \eqref{E3}, which reads, assuming $c_0(s) \neq 0$,
\begin{equation}\label{E6} 
- 2 \theta'_0 (s) \theta'_1(s) + (\mu_0 (s) - \mu_0) =0\,.
\end{equation}
This equation can be seen as the second eikonal equation.  It has a
unique regular solution $\theta_1$ if we add the condition
\begin{equation}\label{E7}
\theta_1(0) =0\,.
\end{equation}

The first transport equation is obtained when looking at the
coefficient in front of $h$ which reads
\begin{equation*}
 (\Lambda_0 - i v_0 (0)) b_3  +  (\Lambda_2-\mu_0)  b_1 + \Lambda_3 b_0= \mu_1 b_0\,,
\end{equation*}
or 
\begin{equation}\label{E8}
 (- \partial_\tau^2 + i \hat v_1  (s) \, \tau -\mu_0(s)) b_1(s,\tau)   + 2   \theta_0'(s)  \partial_s  b_0 (s,\tau)  
 + \theta_0''(s) b_0 (s,\tau)  - \mu_1 b_0 (s,\tau) =0 \,.
\end{equation}
We assume 
\begin{equation} \label{E9}
b_1 (s,\tau) = c_1(s) f_0(s,\tau) + \hat b_1(s,\tau) \,,\, \mbox{ with } \int  f_0 (s,\tau)   \hat b_1 (s,\tau) d\tau =0\,.
\end{equation}
Multiplying \eqref{E8} by $f_0(s,\tau)$ and integrating with respect
to $\tau$, we get
\begin{equation}\label{E10}
  2   \theta_0'(s)  \int  \partial_s   b_0 (s,\tau)  f_0 (s,\tau)  d\tau  +  \theta_0''(s) c_0(s)   = \mu_1 c_0(s) \,,
\end{equation}
which leads to 
\begin{equation}\label{E11}
  2   \theta_0'(s) c'_0(s)   +  \theta_0'' (s) c_0(s)   = \mu_1 c_0(s)  \,,
\end{equation}
where we have used in the last line \eqref{E4a}.  Taking $s=0$ and
assuming $c_0(0)\neq 0$, one gets
\begin{equation}\label{E12}
 \theta_0''(0)   = \mu_1 \,,
\end{equation}
which is also sufficient for solving \eqref{E11}.  We have determined
at this stage $c_0(s)$ assuming for normalization
\begin{equation}
c_0(0) =1\,.
\end{equation}
Coming back to \eqref{E8}, we have to solve, for each $s$ in a
neighborhood of $0$
\begin{equation}
 \bigl(- \partial_\tau^2 + i \hat v_1  (s) \tau -\mu_0(s) \bigr) \hat b_1(s,\tau)  = g_1(s,\tau)\,,
\end{equation}
with $g_1(s,\tau)$ satisfying $\int f_0(s,\tau) g_1(s,\tau) d\tau  =0\,$.\\
At this stage,  the function $c_1$ is free.\\
We continue, one step more, in order to see if the proposed approach
is general.\\
The second transport equation is obtained when looking at the
coefficient in front of $h^\frac 43$, which reads
\begin{equation*}
 (\Lambda_0 - i v_0 (0)) b_4  +  (\Lambda_2-\mu_0)  b_2 + ( \Lambda_3 -\mu_1) b_1 + \Lambda_4 b_0= \mu_2 b_0\,,
\end{equation*}
or
\begin{equation}\label{E8a}
\begin{split}
& (- \partial_\tau^2 + i \hat v_1  (s)   \tau -\mu_0(s)) b_2(s,\tau)  + 2 \theta_0'(s)  \partial_s  b_1 (s,\tau) + \theta_0''(s) b_1 (s,\tau) \\
& -\mu_1 b_1 (s,\tau) - \mu_2  b_0 (s,\tau) + (iv_2(s )\tau^2  - \theta'_1(s)^2)  b_0 (s,\tau)  -  3 \tau ^2 \mathfrak c (s)^2 \theta_0'(s)^2  \\
& + 4 \tau \mathfrak c (s)^2 \theta_0'(s)\theta'_1(s) b_0 +  \mathfrak c (s) \partial_\tau b_0  =0 \,.  \\
\end{split}
\end{equation}
We look for $b_2$ in the form
\begin{equation} \label{E92}
b_2 (s,\tau) = c_2(s) f_0(s,\tau) + \hat b_2(s,\tau) \,,\, \mbox{ with } \int  f_0 (s,\tau)   \hat b_1 (s,\tau) d\tau =0\,.
\end{equation}
We then proceed as before.  If we write
\begin{equation*}
\begin{split}
g_2 (s,\tau) & =  - 2   \theta_0'(s)  \partial_s  b_1 (s,\tau)   \\
 & \quad  - \theta_0''(s) b_1 (s,\tau) + \mu_1 b_1(s,\tau)  +  \mu_2 b_0(s,\tau) (s)  + ( \theta'_1(s)^2-iv_2\tau^2 )  b_0 (s,\tau) \\
 &\quad -  3 \tau ^2 \mathfrak c (s)^2 \theta_0'(s)^2 b_0  + 4 \tau \mathfrak c (s)^2 \theta_0'(s)\theta'_1(s) b_0  
  +  \mathfrak c (s) \partial_\tau b_0 \,, \\
\end{split}
\end{equation*}
the orthogonality condition reads
\begin{equation*}
\begin{split}
0 &= \int g_2(s,\tau) f_0(s,\tau) \, d\tau \\
& =  - 2 \theta'_0 (s) c'_1(s) + (\mu_1 -  \theta_0''(s)) c_1(s)  - 2   \theta_0'(s) \int   \partial_s  \hat b_1 (s,\tau) f_0(s,\tau) d\tau \\
& \quad   +  \left(\mu_2  + \theta'_1(s)^2-iv_2 \int \tau^2 f_0(s,\tau)^2 d\tau  \right) c_0(s)\\
& \quad  +  \int \left( -  3 \tau ^2 \mathfrak c (s)^2 \theta_0'(s)^2 b_0  
+ 4 \tau \mathfrak c (s)^2 \theta_0'(s)\theta'_1(s) b_0 \, f_0 (s,\tau) d\tau\right)  \\
& \quad  +\int   \mathfrak c (s) \partial_\tau b_0\,  f_0(s,\tau) \, d\tau \,. \\
\end{split}
\end{equation*}
Observing that
\begin{equation*}
\begin{split}
&   \int \bigl( -  3 \tau ^2 \mathfrak c (s)^2 \theta_0'(s)^2 b_0(s,\tau)  + 4 \tau \mathfrak c (s)^2 \theta_0'(s)\theta'_1(s) b_0(s,\tau)  \, f_0 (s,\tau) 
+  \mathfrak c (s) \partial_\tau b_0\,  f_0(s,\tau)\bigr)   d\tau\\
& \quad  =  \mathfrak c (0) \left(\int \partial_\tau f_0(0,\tau) f_0(0,\tau) \, d\tau\right)  \, c_0(0)  \,, \\
\end{split}
\end{equation*}
for $s=0$, this determines $\mu_2$ as a necessary condition at $s=0$
which reads
\begin{equation}\label{detmu2}
\mu_2 = i v_2(0) \int \tau^2 f_0(0,\tau)^2 d\tau - \theta'_1(0)^2  -  \mathfrak c (0) \int \partial_\tau f_0(0,\tau) f_0(0,\tau) \, d\tau   \, .
\end{equation}
Note that in the case when $\# \in \{D,N,R\}$, we get
\begin{equation*}
\int \partial_\tau f_0(0,\tau) f_0(0,\tau) \, d\tau = \frac 12  f_0(0,0)^2\,.
\end{equation*}

We can then determine $c_1$ if we add the condition $c_1(0)=0\,$.\\
Since $g_2$ is orthogonal to $\bar f_0$, we can find $\hat b_2$, while
$c_2$ remains free for the next step.\\

Hence, we have obtained the following theorem
\begin{theorem} 
Under the assumptions of Theorem \ref{th:app}, if $\mu_0^\#$ is a
simple eigenvalue of the realization ``$\#$'' of the complex Airy
operator $-\frac{d^2}{dx^2} + ix$ in $L^{2}_{\#}\,$, and $\tilde\mu_1$
is the eigenvalue of the Davies operator $-\frac{d^2}{dy^2} + i y^2$
on $L^2(\mathbb R)$ with the smallest real part, then there exists an
approximate pair $(\lambda_h^{\#,wkb},u_h^{\#,wkb})$ with
$u_h^{\#,wkb}$ in the domain of $\mathcal A_h^\#$, such that
\eqref{defwkb}, \eqref{deftheta} and \eqref{defb} are satisfied and
\begin{equation} 
\exp \left(\frac \theta h\right)\, (\A_h^\# -\lambda_h^\#)\, u_h^{\#,wkb} = \mathcal O(h^\infty) \, 
\mbox{ in } L^{2}_{\#}(\Omega)\,,\, ||u_h^{\#,wkb}||_{L^2} \sim  1\,,
\end{equation}
where
\begin{equation}
\lambda_0^\# = \mu_0^\# \,|\, v_{01}|^\frac 23  \exp \left(i\frac{\pi}{3} \sign \, v_{01}\right)  \,,\quad 
\lambda_2 =\tilde\mu_1  |v_{20}|^\frac 12 \exp \left(i \frac{\pi}{4} \sign v_{20} \right)\,,
\end{equation}
with $v_{01}:= \nu \cdot \nabla V (x^0)\,$.
\end{theorem}

\begin{remark}
In this approach, we understand more directly why no odd power of
$h^\frac 16$ appears for $\lambda_h$.  Note that $\mu_j=\lambda_{2j}$.
\end{remark}

\section{Examples}
\label{sec:examples}

In this Section, we illustrate the above general results for the
potential $V(x) = x_1$ and some simple domains.

\subsection{Disk}
\label{sec:disk_example}

Let $\Omega = \{ (x_1,x_2)\in \R^2 ~:~ |x| < R_0\}$ be the disk of
radius $R_0$.  In this case, $\Omega_\perp = \{(R_0,0), (-R_0,0)\}$.
The local parameterization around the point $(R_0,0)$ reads in polar
coordinates $(r,\theta)$ as $\t = R_0 - r$, $s = R_0
\theta$, so that 
\begin{equation} \label{eq:Vdisk}
V(x) = x_1(s,\t) = (R_0-\rho) \cos(s/R_0),
\end{equation}
$\curv(0) = 1/R_0$, and we get
\begin{equation}
v_{00} = R_0, \quad \, v_{01} = - 1,  \quad v_{20} = - \frac{1}{2R_0} , \quad v_{11} = v_{02} = 0\, .
\end{equation}
Using Eqs. \eqref{eq:lambda4N}, \eqref{eq:lambda4D},
\eqref{eq:lambda4R} or \eqref{eq:lambda4T} for $\lambda_4^{\#,(n)}$,
one can write explicitly the four-term expansion for four types of
boundary condition: \\
$\bullet$ Dirichlet case, 
\begin{equation} 
\label{eq:lambda_app2D}
\lambda_{h}^{D, (n,k)} = i R_0 - h^{\frac 23} a_n e^{-i\pi/3} + h (2k-1) \frac{e^{-i\pi/4}}{\sqrt{2 R_0}} + \mathcal O(h^{\frac 53})\, .
\end{equation} \\
$\bullet$ Neumann case
\begin{equation}  
\label{eq:lambda_app2N}
\lambda_{h}^{N, (n,k)}= i R_0 - h^{\frac 23} a'_n e^{-i\pi/3} + h (2k-1)\frac{e^{-i\pi/4}}{\sqrt{2 R_0}} + h^{\frac 43} 
\frac{e^{-\pi i/6}}{2R_0\, a'_n} + \mathcal O(h^{\frac 53})\, .
\end{equation} \\
$\bullet$ Robin case
\begin{equation}
\label{eq:lambda_app2R}
\begin{split}
\lambda_{h}^{R, (n,k)} & = i R_0 - h^{\frac 23} a_n^R(\kappa) e^{-i\pi/3}  + h (2k-1) \frac{e^{-i\pi/4}}{\sqrt{2 R_0}} \\
& + h^{\frac 43} \frac{i}{2R_0 (\kappa^2 - a_n^R(\kappa) e^{-\pi i/3})} + \mathcal O(h^{\frac 53})\, . \\
\end{split}
\end{equation}
When $\kappa = 0$, $a_n^R(0) = a'_n$, and one retrieves the expansion
(\ref{eq:lambda_app2N}) for Neumann case. \\
$\bullet$ Transmission case, 
\begin{equation}  
\label{eq:lambda_app2T}
\lambda_{h}^{T, (n,k)} =  i R_0 - h^{\frac 23} a_n^+(\kappa) e^{-i\pi/3}
 + h (2k-1)\frac{e^{-i\pi/4}}{\sqrt{2 R_0}} + h^{\frac 43} \frac{e^{-i\pi/6}}{2R_0\, a_n^+(\kappa) }
  + \mathcal O(h^{\frac 53})\, .
\end{equation}
When $\kappa = 0$, one has $a_n^+(0) = a'_n$ and thus retrieves
the expansion (\ref{eq:lambda_app2N}) for Neumann case.

We recall that the indices $n = 1,2,\ldots$ and $k = 1,2,\ldots$
enumerate eigenvalues of the operators $\L_0^\#$ and $\L_2^\#$ that
were used in the asymptotic expansion.  The approximate eigenvalue
with the smallest real part corresponds to $n = k = 1$.

The three-terms version of the Neumann expansion
(\ref{eq:lambda_app2N}) was first derived by de Swiet and Sen
\cite{deSwiet94} (note that we consider the eigenvalues of the
operator $-h^2 \Delta + ix_1$ while de Swiet and Sen looked at the
complex conjugate operator).

\begin{remark} \label{lem:complex_conj}
At the other point $(-R_0,0)$, the parameterization is simply 
\begin{equation*}
V(x) = -(R_0-\rho) \cos(s/R_0)
\end{equation*}
that alters the signs of the all involved coefficients $v_{jk}$.  As a
consequence, the asymptotics is obtained as the complex conjugate of
$\lambda_h^{\#,(n,k)}$.
\end{remark}

In the WKB approach, one needs to compute the functions $\theta_0(s)$
and $\theta_1(s)$ that determine the asymptotic decay of the quasimode
state in the tangential direction.  We only consider the Neumann
boundary condition while the computation for other cases is similar.
From \eqref{eq:vj_WKB} and \eqref{eq:v1hat_WKB}, we have for the
potential in (\ref{eq:Vdisk}):
\begin{equation*}
v_0(s) = R_0 \cos (s/R_0),  \qquad v_1(s) = - \cos (s/R_0),  \qquad \hat v_1(s) = 2-3\cos (s/R_0).
\end{equation*}
In what follows, we consider $s > 0$ though the results will be the
same for $s < 0$ due to the symmetry.  From Eqs. (\ref{E0}, \ref{E1}),
we first obtain
\begin{equation}
\label{eq:theta0_disk}
\theta_0(s) = \int\limits_0^s \sqrt{-iR_0(1-\cos(s'/R_0))} \, ds' = e^{-\pi i/4}\, (2R_0)^{\frac 32} \bigl(1 - \cos(s/(2R_0))\bigr) \, .
\end{equation}
For Neumann boundary condition, $\mu_0 = -a'_1 e^{-\pi i/3}$ (here
$v_1(0)=-1$) and the eigenvalue of the operator in (\ref{eq:WKB_mu0s})
reads
\begin{equation*}
\mu_0(s) = -a'_1 \, |2-3\cos(s/R_0)|^{\frac 23} \exp\left(\frac{\pi i}{3} \sign (2-3\cos(s/R_0))\right) .
\end{equation*}
Since $\hat v_1(s)$ was assumed to be nonzero, we restrict the
analysis to $|s/R_0| < \arccos(2/3)$ for which $2-3\cos x$ does not
vanish (and remains negative) so that
\begin{equation}
\mu_0(s) = -a'_1 \, \bigl(3\cos(s/R_0)-2\bigr)^{\frac 23} \exp\left(-\frac{\pi i}{3}\right) .
\end{equation}
From (\ref{E6}), one gets then
\begin{equation}
\label{eq:theta1_disk}
\begin{split}
\theta_1(s) & = \int\limits_0^s \frac{-a'_1 \, e^{-\pi i/3} \left[(3\cos(s'/R_0)-2)^{\frac 23}  - 1\right]}
{2e^{-\pi i/4} R_0^{\frac 12} \sqrt{1-\cos(s'/R_0)}}\, ds'\\
& = \frac12 |a'_1| \, e^{-\pi i/12} R_0^{\frac 12} \int\limits_0^{s/R_0} \frac{(3\cos x -2)^{\frac 23} - 1}{\sqrt{1-\cos x}} \, dx \,.  \\
\end{split}
\end{equation}

\subsection{Annulus}
\label{sec:annulus_example}

For an annulus $\Omega = \{ (x_1,x_2)\in \R^2 ~:~ R_1 < |x| < R_2\}$
between two circles of radii $R_1$ and $R_2$, there are four points in
$\Omega_\perp$: $(\pm R_1,0)$ and $(\pm R_2,0)$.  In order to
determine the candidate for an eigenvalue with the smallest real part
(in short the ``first eigenvalue''), one needs to compare the
asymptotics of the quasimodes associated with these points and
identify those with the smallest real part.  Of course, this analysis
depends on the imposed boundary conditions.  We consider four
combinations: NN (Neumann condition on both circles), ND (Neumann
condition on the inner circle and Dirichlet on the outer circle), DN
(Dirichlet condition on the inner circle and Neumann on the outer
circle), and DD (Dirichlet condition on both circles).  Since the
leading contribution is proportional $|a_1|\approx 2.3381$ for the
Dirichlet case and to $|a'_1| \approx 1.0188$ for the Neumann case,
the asymptotics for the circle with Neumann boundary condition always
contributes to the first eigenvalue.  In turn, when the same boundary
condition is imposed on the two circles, the first eigenvalue
expansion corresponds to the outer circle of larger radius because the
real part of the next-order term (of order $h$) is always positive and
scales as $1/\sqrt{R_0}$.  As a consequence, the first eigenvalue
asymptotics is given by (\ref{eq:lambda_app2N}) with $R_0 = R_2$ for
cases NN and DN, and by (\ref{eq:lambda_app2D}) with $R_0 = R_2$ for
the case DD.  Only in the case ND, the first eigenvalue asymptotics is
determined by the points $(\pm R_1,0)$ on the inner circle.  In this
case, the potential reads in local coordinates around $(R_1,0)$ as
$V(s,\t) = (R_1 + \t) \cos(s/R_1)$ so that the only change with
respect to the above results is $\, v_{01} = 1$ (instead of $\, v_{01}
= -1$) and $\curv(0) = -1/R_1$ (instead of $\curv(0) = 1/R_1$) so that
Eq. (\ref{eq:lambda_app2N}) becomes
\begin{equation}
\label{eq:lambda_app3}
\lambda_{app}^{ND,(n,k)} = i R_1 + h^{\frac 23} |a'_n| e^{i\pi/3} + h (2k-1) \frac{e^{-i\pi/4}}{\sqrt{2 R_1}} 
+ h^{\frac 43} \frac{e^{\pi i/6}}{2|a'_n| R_1} + \mathcal O(h^{\frac 53})\, .
\end{equation}

\begin{remark}
When the outer radius $R_2$ of an annulus goes to infinity, the above
problem should progressively\footnote{We do not have a mathematical
proof, the statement remains conjectural.} become an exterior problem
in the complement of a disk: $\Omega = \{ (x_1,x_2)\in \R^2 ~:~ |x| >
R_1\}$.  Due to the local character of the asymptotic analysis, the
expansion (\ref{eq:lambda_app3}) is independent of the outer radius
$R_2$ and holds even for the unbounded case.  This argument suggests
the non-emptiness of the spectrum for unbounded domains.  This
conjecture is confirmed by numerical results in
Sec. \ref{sec:numerics}.
\end{remark}

\subsection{Domain with transmission condition}
\label{sec:twolayers_example}

Finally, we consider the union of two subdomains, the disk $\Omega_- =
\{ (x_1,x_2)\in \R^2~:~ |x| < R_1\}$ and the annulus $\Omega_+ = \{
(x_1,x_2)\in \R^2~:~ R_1 < |x| < R_2\}$ separated by a circle on which
the transmission boundary condition is imposed.  A Dirichlet, Neumann
or Robin boundary condition can be imposed at the outer boundary
(circle of radius $R_2$).  As for the annulus, there are four points
in $\Omega_\perp$: $(\pm R_1,0)$ and $(\pm R_2,0)$.  Here we focus
only on the asymptotic behavior at points $(\pm R_1,0)$ for the
transmission boundary condition (the behavior at the points $(\pm
R_2,0)$ was described in Sec. \ref{sec:disk_example}).  We consider
the case described in Theorem \ref{thm:scaling} when the transmission
parameter $\kappa$ scales with $h$ according to
(\ref{eq:kappa_newscaling}).  As discussed in
Sec. \ref{sec:kappa_scaling}, this situation is relevant for diffusion
MRI applications.  The case with fixed $\kappa$ can be treated
similarly.

As stated in Theorem \ref{thm:scaling}, the asymptotic expansion is
obtained by starting from the ``basic'' expansion (with $\kappa =0$)
of either of two problems with Neumann boundary condition
corresponding to the two subdomains $\Omega_-$ and $\Omega_+$.\\
If we start from the expansion for the disk, one has $V(x) =
(R_1-\rho) \cos(s/R_1),$ and the asymptotic expansion
(\ref{newscaling}) at the point $(R_1,0)$ reads
\begin{equation}
\label{eq:lambda_app2Th}
\begin{split}
\lambda^{\#,(n,k)}_h & = i \, R_1 - h^{\frac 23} a'_n  e^{-\pi i/3} + h\, (2k-1)\, \frac{e^{-\pi i/4}}{\sqrt{2R_1}}   \\
& \qquad - h^{\frac 43} \frac{e^{-\pi i/6}}{a'_n} \left(\hat \kappa - \frac{1}{2R_1}\right) + O(h^{\frac 53})\,.  \\
\end{split}
\end{equation}
In turn, if we start from the expansion for the inner boundary of the
annulus, one has $V(x) = (R_1 + \rho) \cos(s/R_1)$, and the asymptotic
expansion (\ref{newscaling}) at the point $(R_1,0)$ reads 
\begin{equation}
\label{eq:lambda_app2Th2}
\begin{split}
\lambda^{\#,(n,k)}_h & = i \, R_1 - h^{\frac 23} a'_n  e^{\pi i/3} + h\, (2k-1)\, \frac{e^{-\pi i/4}}{\sqrt{2R_1}}   \\
& \qquad - h^{\frac 43} \frac{e^{\pi i/6}}{a'_n} \left(\hat \kappa + \frac{1}{2R_1}\right) + O(h^{\frac 53})\,.  \\
\end{split}
\end{equation}
These two expressions are different, in particular, their imaginary
parts differ already in the order $h^{\frac 23}$.  In turn, the real
parts differ at the term of order $h^{\frac 43}$ that contains two
contributions: from the curvature of the boundary, and from the
transmission.  While the curvature changes its sign on both sides of
the boundary, the contribution due to the transmission remains the
same.  As a consequence, the real part of (\ref{eq:lambda_app2Th2}) is
larger than the real part of (\ref{eq:lambda_app2Th}).  One can thus
expect the existence of two distinct eigenstates living on both sides
of the boundary, as confirmed numerically in the next section.  For
$k=1$, the eigenstate associated with the eigenvalue with the smallest
real part is mainly localized in the disk side of the boundary.

\section{Numerical results}
\label{sec:numerics}

This section presents some numerical results to illustrate our
analysis.  The claims of this section are supported by numerical
evidence but should not be considered as rigorous statements, in
contrast to previous sections.

The numerical analysis will be limited to bounded domains, for which
the BT-operator has compact resolvent and hence  discrete spectrum (see
Sec.~\ref{sec:definition}).  In order to compute numerically its
eigenvalues and eigenfunctions, one needs to approximate the
BT-operator in a matrix form.  For this purpose, one can either (i)
discretize the domain by a square lattice and replace the Laplace
operator by finite differences (finite difference method); (ii)
discretize the domain by a mesh and use a weak formulation of the
eigenvalue problem (finite elements method); or (iii) project the
BT-operator onto an appropriate complete basis of functions.  We
choose the last option and project the BT-operator onto the Laplacian
eigenfunctions which for rotation-invariant domains (such as disk,
annuli, circular layers) are known explicitly \cite{Grebenkov13}.  In
this basis, the Laplace operator $-\Delta$ is represented by a
diagonal matrix $\Lambda$.  Moreover, the matrix representation of the
potential $V(x) = x_1$ was computed analytically, i.e., the elements
of the corresponding matrix $\B$ are known explicitly
\cite{Grebenkov07,Grebenkov08,Grebenkov10}.  As a consequence, the
computation is reduced to finding the Laplacian eigenvalues for these
rotation-invariant domains, constructing the matrices $\Lambda$ and
$\B$ through explicit formulas, and then diagonalizing numerically the
truncated matrix $h^2 \Lambda + i\B$ which is an approximate
representation of the BT-operator $-h^2\Delta + ix_1$.  This numerical
procedure yields the eigenvalues $\lambda^{(m)}_h$ of the truncated
matrix $h^2 \Lambda + i\B$, while the associated eigenvectors allow
one to construct the eigenfunctions $u^{(m)}_h$.  All eigenvalues are
ordered according to their increasing real parts:
\begin{equation}
\label{eq:eigen_ordering}
\Re\{\lambda^{(1)}_h\} \leq \Re\{\lambda^{(2)}_h\} \leq \ldots 
\end{equation}

Note that, for a bounded domain, the potential $ix$ is a bounded
perturbation of the unbounded Laplace operator $-h^2\Delta$, if $h \ne
0$.  To preserve this property after truncation of the matrix $h^2
\Lambda + i\B$, the truncation size should be chosen such that $h^2
\mu^{(M)} \gg 1$, where $\mu^{(M)}$ is the largest element of the matrix
$\Lambda$.  Due to the Weyl's law, $M\sim \frac{|\Omega|}{4\pi}
\mu^{(M)}$ so that the truncation size $M$ should satisfy:
\begin{equation}
h^2 M \gg \frac{|\Omega|}{4\pi},
\end{equation}
where $|\Omega|$ is the surface area of $\Omega$.  For larger domains,
either larger truncation sizes are needed (that can be computationally
limiting), or $h$ should be limited to larger values.  In practice, we
use $M$ around $3000$ to access $h$ up to $0.01$.  We have checked
that the truncation size does not affect the computed eigenvalues.

\subsection{Eigenvalues}

For large $h$, one can divide the BT operator by $h^2$, $-\Delta +
ix_1/h^2$, to get a small bounded perturbation of the Laplace
operator.  In particular, the eigenvalues of the operator $-h^2 \Delta
+ ix_1$ behave asymptotically as $h^2 \mu^{(m)}$, where $\mu^{(m)}$
are the eigenvalues of the Laplace operator.  In this Section, we
focus on the more complicated semi-classical limit $h\to 0$ which is
the main topic of the paper.

\subsubsection{Disk}

In order to check the accuracy of the asymptotic expansion of
eigenvalues, we first consider the BT-operator in the unit disk:
$\Omega = \{ (x_1,x_2)\in\R^2~:~ |x| < R_0\}$, with $R_0 = 1$.  We
will present {\it rescaled} eigenvalues, $(\lambda^{(m)}_h -
iR)/h^{\frac 23}$, for which the constant imaginary offset $iR$ is
subtracted and the difference $\lambda^{(m)}_h - iR$ is divided by
$h^{\frac 23}$ in order to emphasize the asymptotic behavior.  Note
also that, according to Remark \ref{lem:complex_conj}, the asymptotic
expansions for the approximate eigenvalues corresponding to the points
$(-R,0)$ and $(R,0)$ are the complex conjugates to each other.  In
order to facilitate their comparison and check this property for
numerically computed eigenvalues, we will plot the absolute value of
the imaginary part.

\begin{figure}
\begin{center}
\includegraphics[width=62mm]{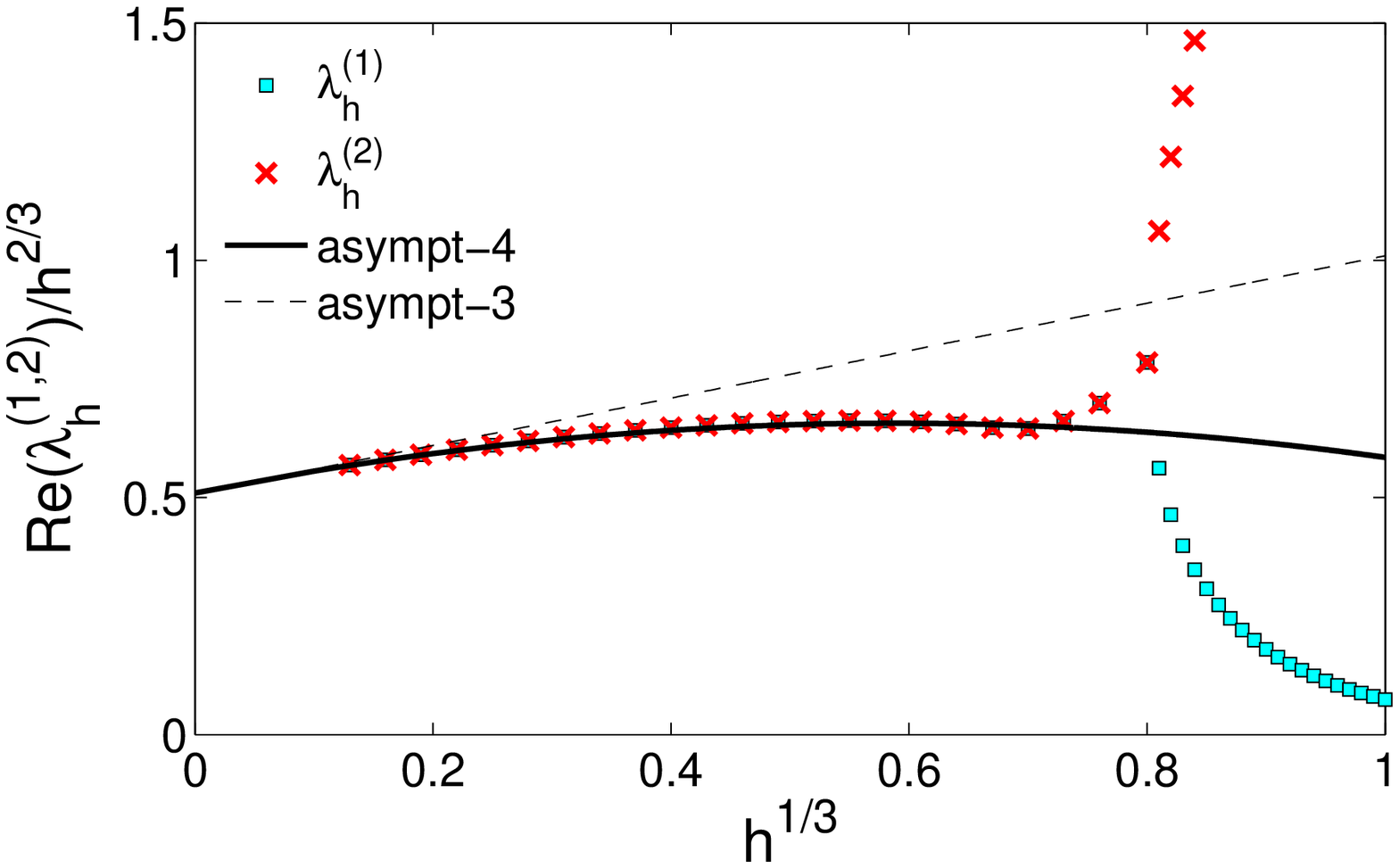}  
\includegraphics[width=62mm]{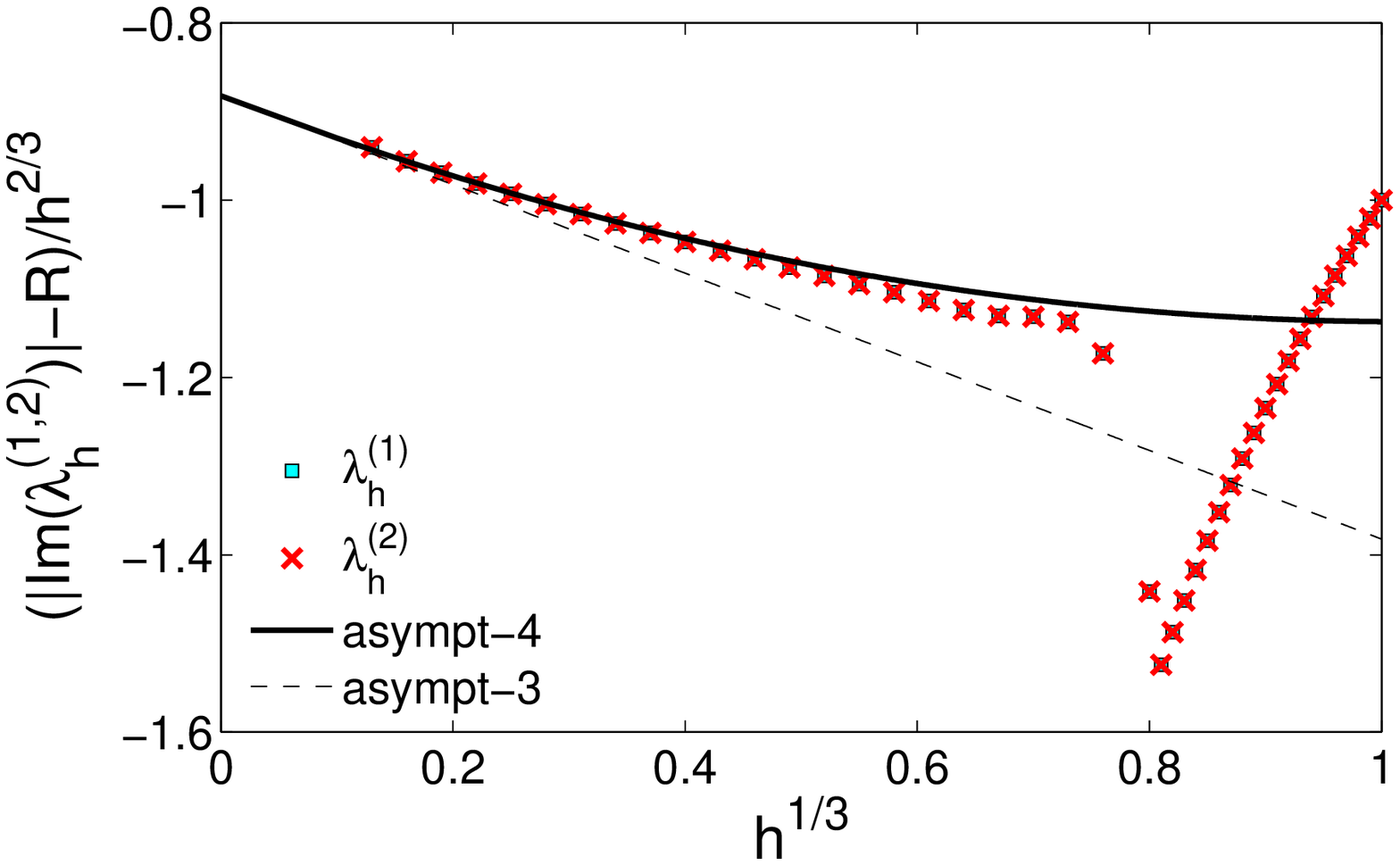}  
\end{center}
\caption{
The rescaled eigenvalues $\lambda^{(1)}_h$ and $\lambda^{(2)}_h$ of
the BT-operator in the unit disk with Neumann boundary condition.
Symbols (squares and crosses) show the numerical results of the
diagonalization of the matrix $h^2 \Lambda + i\B$ (truncated to the
size $2803 \times 2803$), solid line presents the four-terms
asymptotics (\ref{eq:lambda_app2N}) for $\lambda_h^{N,(1,1)}$ while
the dashed line shows its three-terms versions (without $h^{\frac 43}$
term).  }
\label{fig:diskN_mu12}
\end{figure}

Figure \ref{fig:diskN_mu12} shows the first two eigenvalues
$\lambda^{(1)}_h$ and $\lambda^{(2)}_h$.  For $h^{\frac 13} \lesssim
0.8$, these eigenvalues turn out to be the complex conjugate to each
other, as expected from their asymptotic expansions (the difference
$\lambda^{(1)}_h - \bar{\lambda}^{(2)}_h$ being negligible within
numerical precision).  In turn, the eigenvalues $\lambda^{(1)}_h$ and
$\lambda^{(2)}_h$ become real and split for $h^{\frac 13} \gtrsim
0.8$.  The splitting is expected because these eigenvalues behave
differently in the large $h$ limit.  This numerical observation
suggests the existence of branch points in the spectrum (similar
features were earlier reported for the complex Airy operator on the
one-dimensional interval with Neumann boundary condition, see
\cite{Stoller91}).  For comparison, the four-terms asymptotics
(\ref{eq:lambda_app2N}) for $\lambda_h^{N,(1,1)}$ and its three-terms
version (without term $h^{\frac 43}$) are shown by solid and dashed
lines, respectively.  These expansions start to be applicable for
$h^{\frac 13} \lesssim 0.7$, while their accuracy increases as $h$
decreases.\\

\begin{figure}
\begin{center}
\includegraphics[width=62mm]{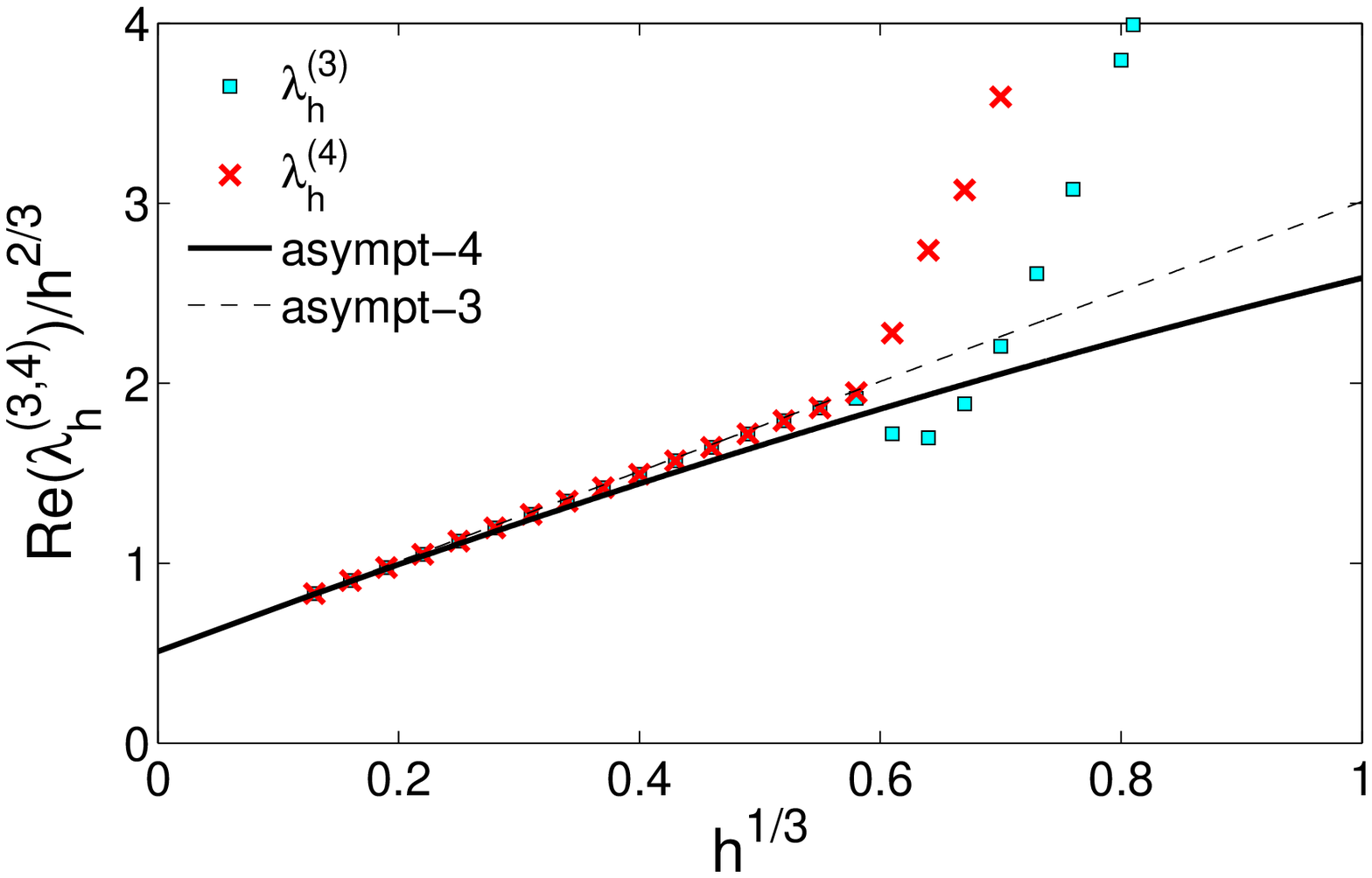}   
\includegraphics[width=62mm]{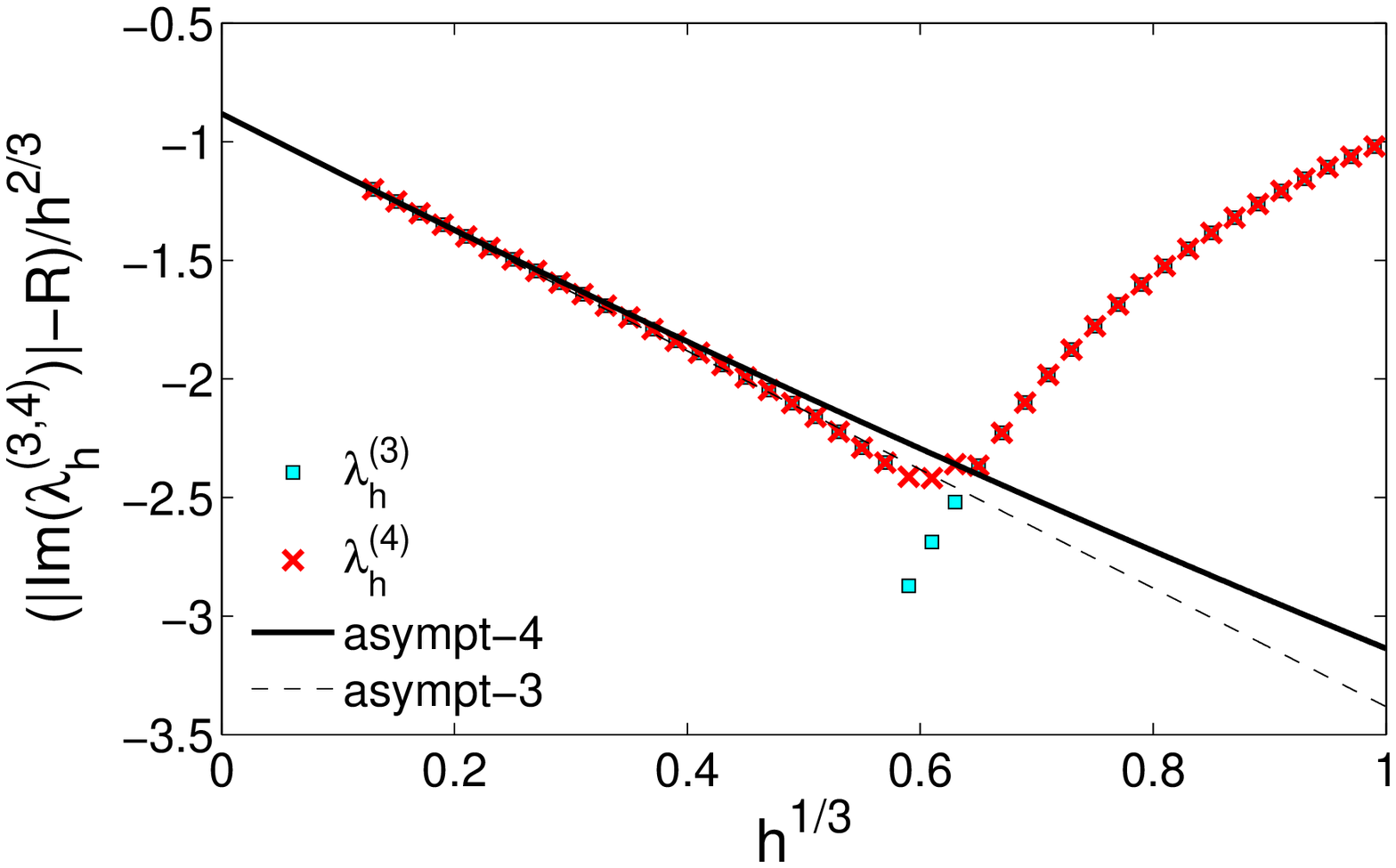}   
\end{center}
\caption{
The rescaled eigenvalues $\lambda^{(3)}_h$ and $\lambda^{(4)}_h$ of
the BT-operator in the unit disk with Neumann boundary condition.
Symbols (squares and crosses) show the numerical results of the
diagonalization of the matrix $h^2 \Lambda + i\B$ (truncated to the
size $2803 \times 2803$), solid line presents the four-terms
asymptotics (\ref{eq:lambda_app2N}) for $\lambda_h^{N,(1,3)}$ while
the dashed line shows its three-terms versions (without $h^{\frac 43}$
term).  }
\label{fig:diskN_mu34}
\end{figure}

Figure \ref{fig:diskN_mu34} shows the next eigenvalues
$\lambda^{(3)}_h$ and $\lambda^{(4)}_h$, the four-terms asymptotics
(\ref{eq:lambda_app2N}) for $\lambda_h^{N,(1,3)}$ and its three-terms
version.  These eigenvalues are the complex conjugates to each other
for $h^{\frac 13} \lesssim 0.57$ while become real and split for
larger $h$.  One can see that the four-terms asymptotics is less
accurate for these eigenvalues than for those from
Fig. \ref{fig:diskN_mu12}.  A small deviation can probably be
attributed to higher-order terms (it is worth noting that
contributions from the $h^{\frac 43}$ and $h^{\frac 53}$ terms can be
comparable for the considered values of $h$).\\

\begin{figure}
\begin{center}
\includegraphics[width=62mm]{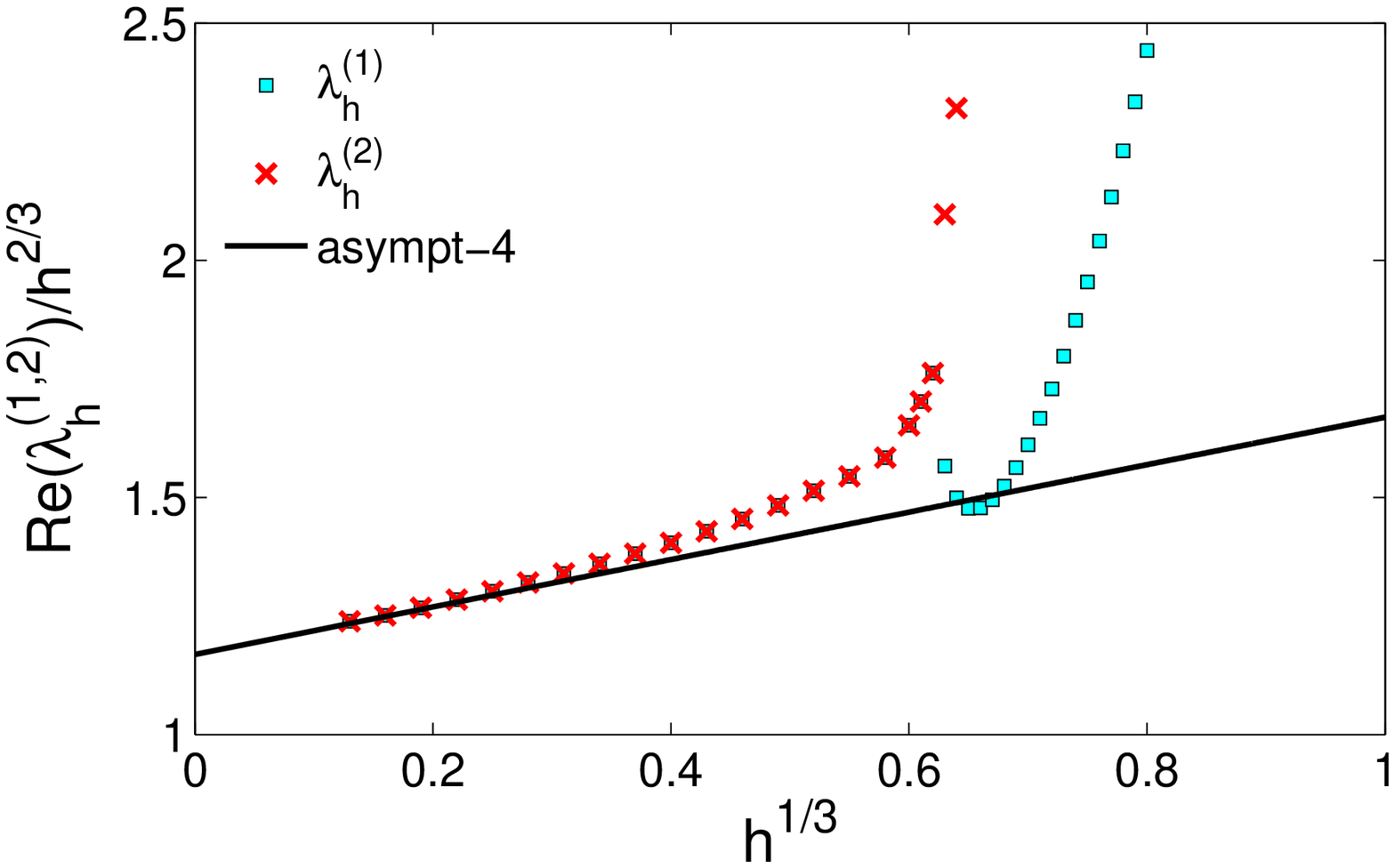} 
\includegraphics[width=62mm]{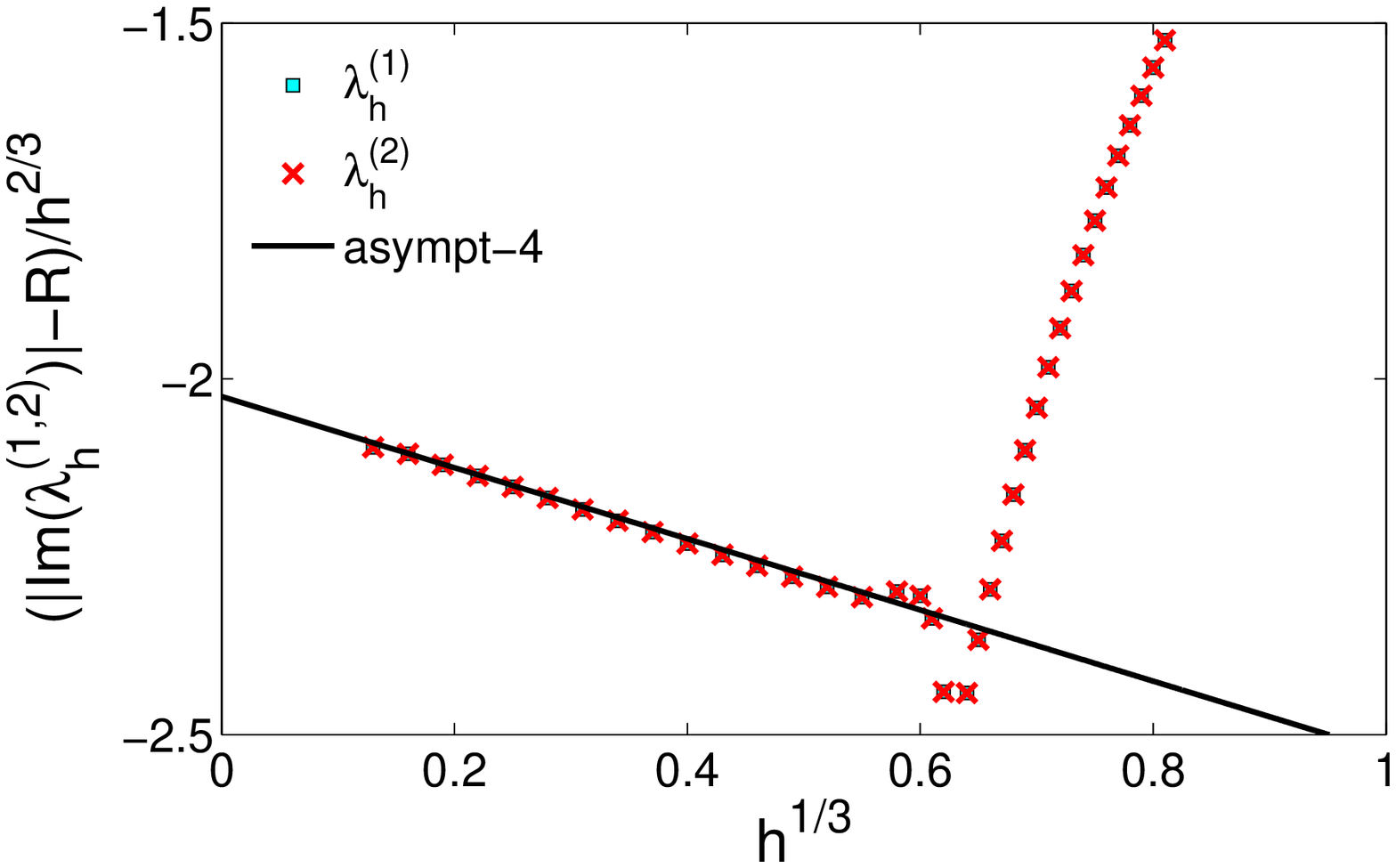} 
\end{center}
\caption{
The rescaled eigenvalues $\lambda^{(1)}_h$ and $\lambda^{(2)}_h$ of
the BT-operator in the unit disk with Dirichlet boundary condition.
Symbols (squares and crosses) show the numerical results of the
diagonalization of the matrix $h^2 \Lambda + i\B$ (truncated to the
size $2731 \times 2731$), while solid line shows the four-terms
asymptotic expansion (\ref{eq:lambda_app2D}) for
$\lambda_h^{D,(1,1)}$. }
\label{fig:diskD_mu}
\end{figure}

For comparison, Figure \ref{fig:diskD_mu} shows the first rescaled
eigenvalues $\lambda^{(1)}_h$ and $\lambda^{(2)}_h$ of the BT-operator
in the unit disk with Dirichlet boundary condition.  As earlier for
the Neumann case, these eigenvalues are complex conjugate to each
other for $h^{\frac 13} \lesssim 0.6$ while become real and split for
larger $h$.  One can see that the asymptotics (\ref{eq:lambda_app2D})
for $\lambda_h^{D,(1,1)}$ captures the behavior for the imaginary part
very accurately.  In turn, the behavior of the real part is less
accurate, probably due to higher-order corrections. \\
 
\begin{figure}
\begin{center}
\includegraphics[width=62mm]{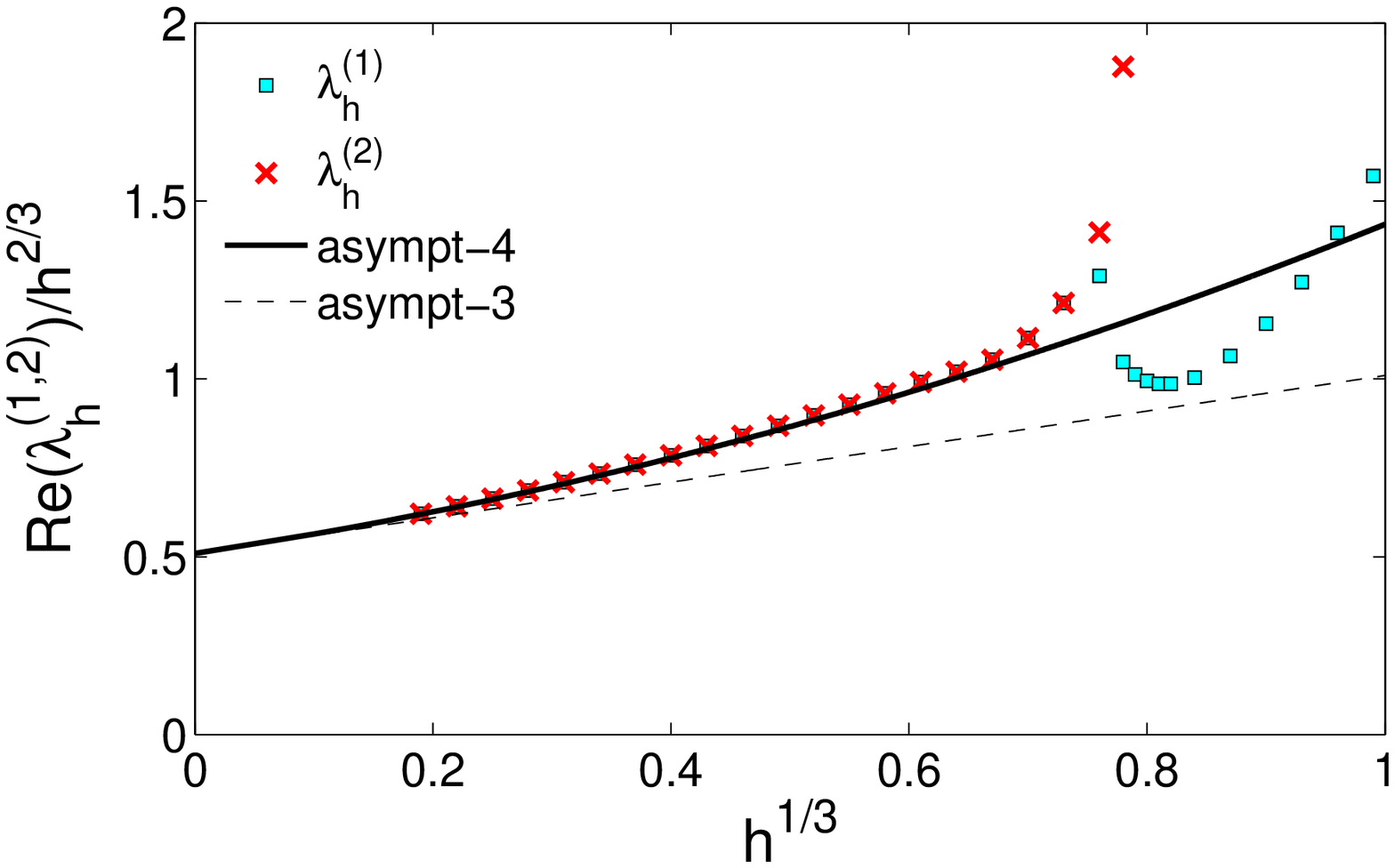}  
\includegraphics[width=62mm]{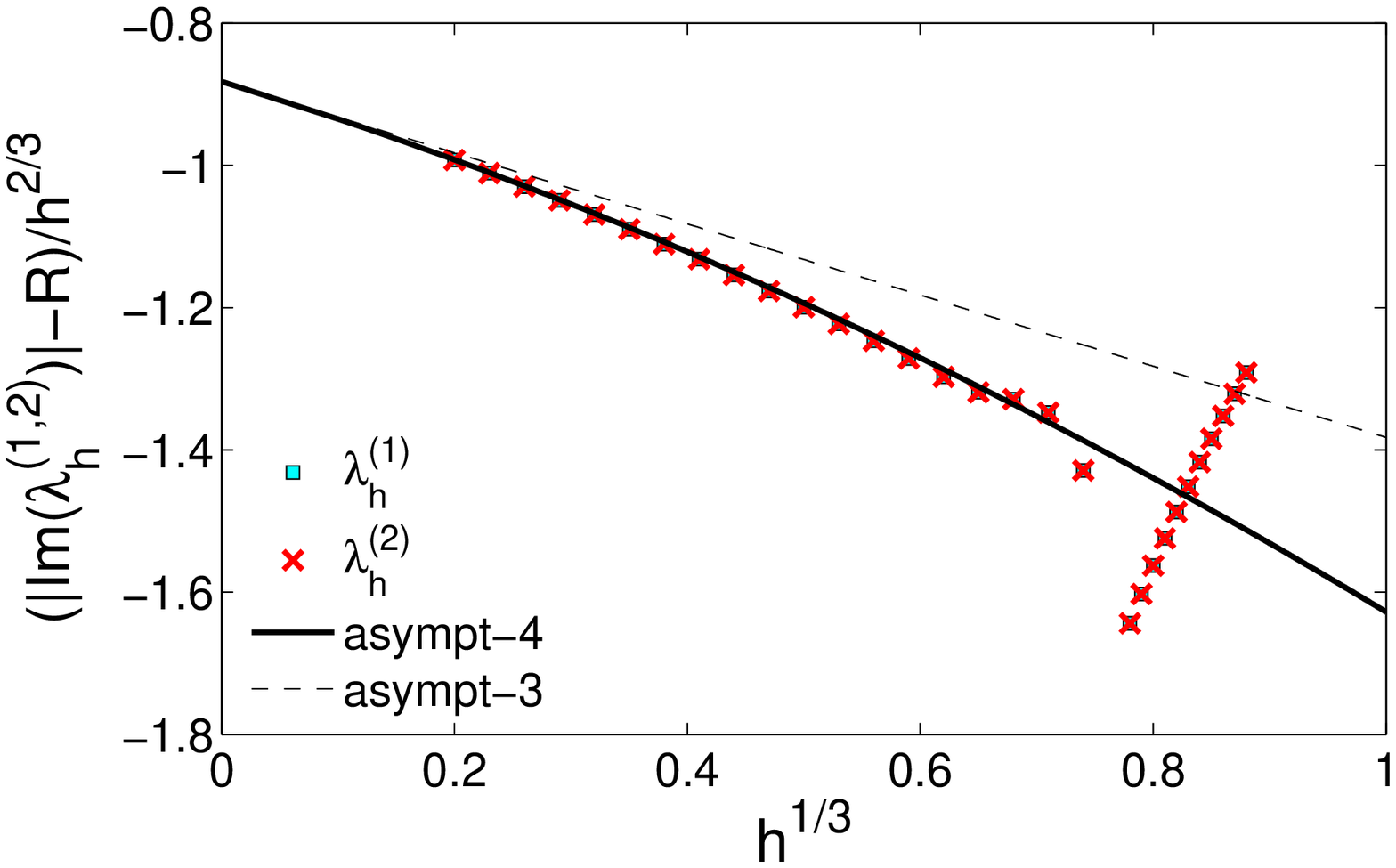}  
\end{center}
\caption{ 
The rescaled eigenvalues $\lambda^{(1)}_h$ and $\lambda^{(2)}_h$ of
the BT-operator in the unit disk with Robin boundary condition (with
$\hat\kappa = 1$ and $\kappa = \hat\kappa h^{\frac23}$).  Symbols
(squares and crosses) show the numerical results of the
diagonalization of the matrix $h^2 \Lambda + i\B$ (truncated to the
size $2803 \times 2803$), while solid and dashed lines show the
four-terms asymptotic expansion (\ref{eq:lambda_app2D}) for
$\lambda_h^{R,(1,1)}$ and its three-term version (without term
$h^{\frac 43}$). }
\label{fig:diskR_mu}
\end{figure}

Finally, Figure \ref{fig:diskR_mu} illustrates the case with Robin
boundary condition, with $\hat\kappa = 1$ while $\kappa$ scaling as
$\hat\kappa h^{\frac23}$.  The four-term expansion (\ref{newscaling})
accurately captures their asymptotic behavior.

\subsubsection{Annulus}

Due to its local character, the quasimodes construction is expected to
be applicable to the exterior problem, i.e., in the complement of a
disk of radius $R_1$, $\Omega = \{ (x_1,x_2)\in \R^2 ~:~ |x| >R_1\}$.
Since we cannot numerically solve this problem for unbounded domains,
we consider a circular annulus $\Omega = \{ x\in \R^2~:~ R_1 < |x|
<R_2\}$ with a fixed inner radius $R_1 = 1$ and then increase the
outer radius $R_2$.  In the limit $h\to 0$, the eigenfunctions are
expected to be localized around the four points $(\pm R_1,0)$, $(\pm
R_2,0)$ from the set $\Omega_\perp$, with corresponding asymptotic
expansions for eigenvalues.

\begin{figure}
\begin{center}
\includegraphics[width=62mm]{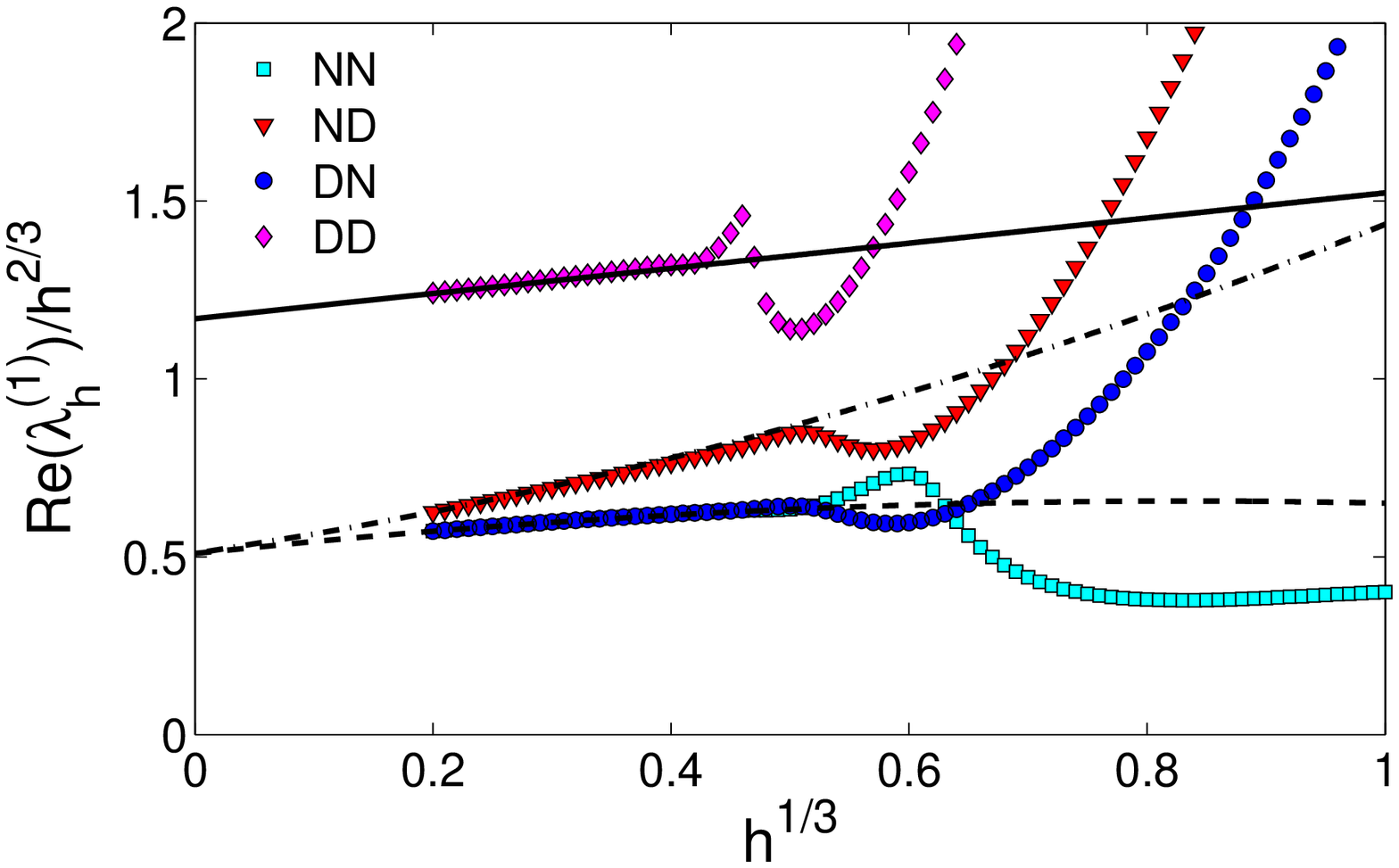} 
\includegraphics[width=62mm]{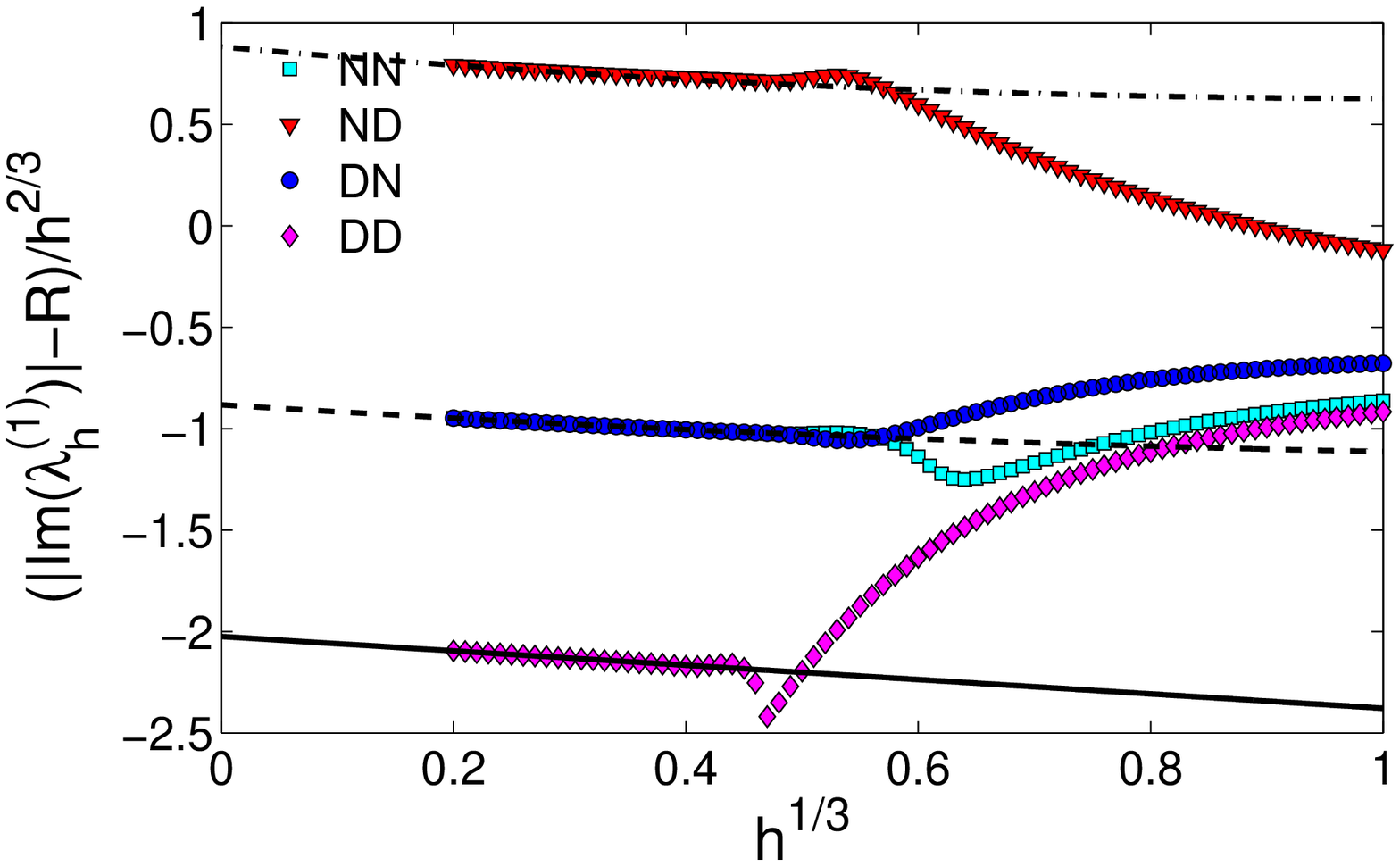} 
\end{center}
\caption{
The rescaled eigenvalue $\lambda^{(1)}_h$ of the BT-operator in the
annulus with four combinations of Neumann/Dirichlet boundary
conditions at the inner and outer circles of radii $R_1 = 1$ and $R_2
= 2\,$: NN (squares), ND (triangles), DN (circles), and DD (diamonds),
obtained by the diagonalization of the truncated matrix $h^2 \Lambda +
i\B$.  The solid line presents the expansion (\ref{eq:lambda_app2D}) with
$R = R_2$ for Dirichlet condition, the dashed line shows the expansion
(\ref{eq:lambda_app2N}) with $R = R_2$ for Neumann condition, and
the dash-dotted line shows the expansion (\ref{eq:lambda_app3}) with $R =
R_1$ for Neumann condition. }
\label{fig:annulus_combined_mu}
\end{figure}

Figure \ref{fig:annulus_combined_mu} illustrates the discussion in
Sec. \ref{sec:annulus_example} about different asymptotics of the
first eigenvalue $\lambda^{(1)}_h$ for four combinations of
Neumann/Dirichlet boundary conditions on inner and outer circles.  In
particular, one observes the same asymptotic expansion
(\ref{eq:lambda_app2N}) with $R = R_2$ for NN and DN cases because the
first eigenvalue is determined by the local behavior near the point
$(R_2,0)$ which is independent of the boundary condition on the inner
circle as $h\to 0\,$.  The expansion (\ref{eq:lambda_app2D}) with $R =
R_2$ for the Dirichlet condition appears only for the case DD.
Finally, the case ND is described by the local behavior at the inner
circle by the expansion (\ref{eq:lambda_app3}) with $R = R_1$.  In
what follows, we focus on this case in order to illustrate that the
local behavior at the inner boundary is not affected by the position
of the outer circle as $h\to 0\,$.

\begin{figure}
\begin{center}
\includegraphics[width=62mm]{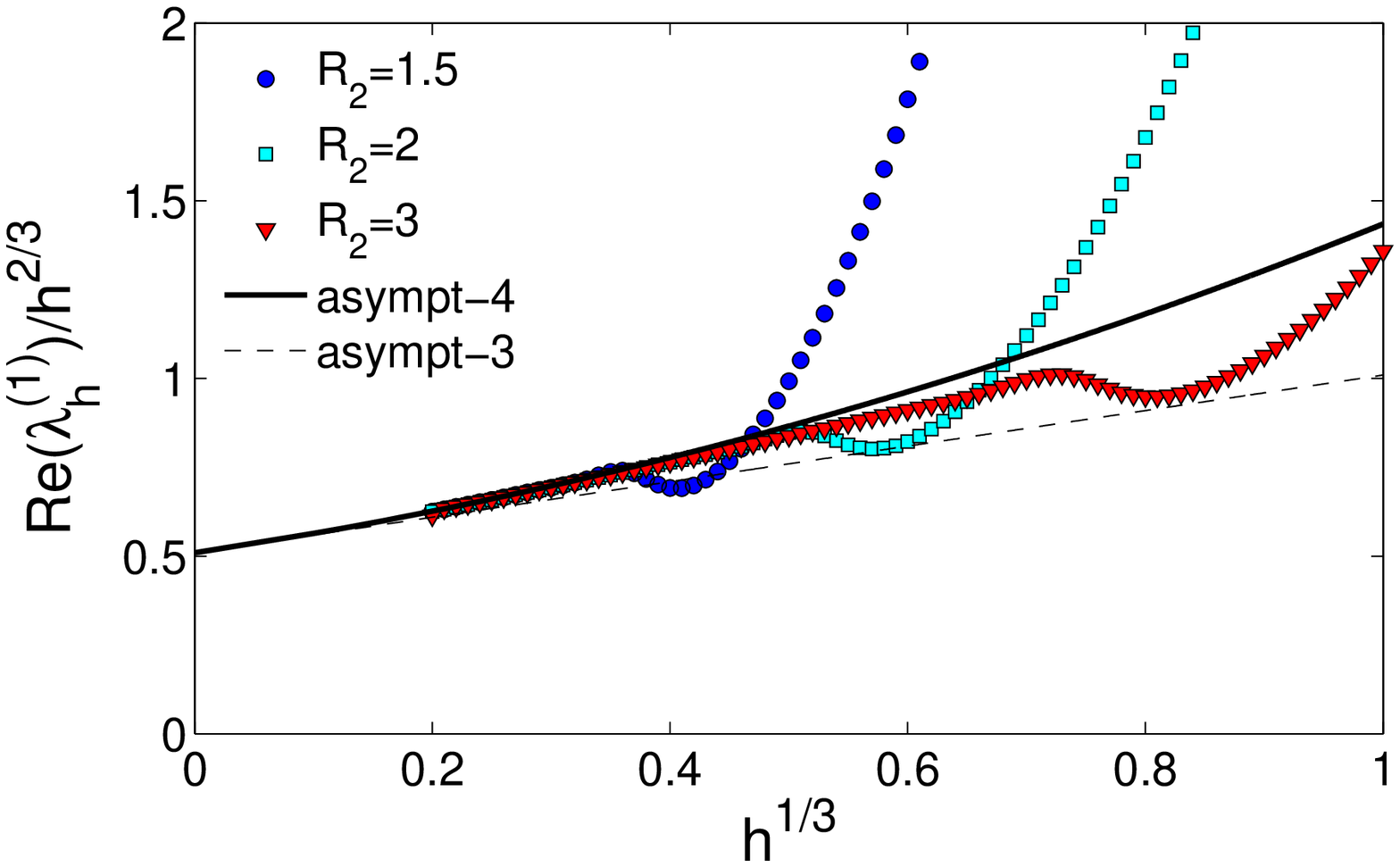} 
\includegraphics[width=62mm]{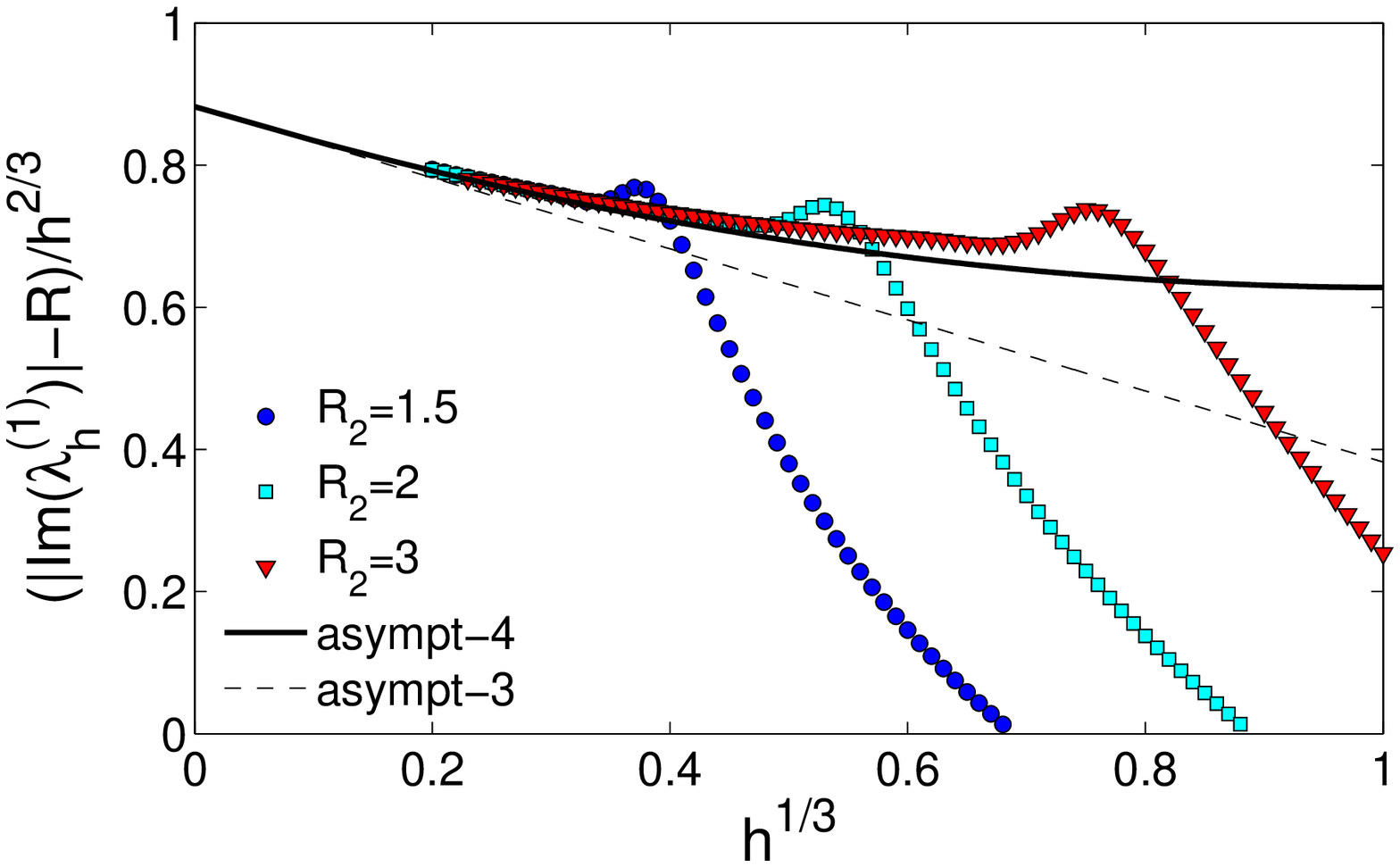} 
\end{center}
\caption{
The rescaled eigenvalue $\lambda^{(1)}_h$ of the BT-operator in the
annulus with Neumann boundary condition at the inner circle of radius
$R_1 = 1$ and Dirichlet boundary condition at the outer circle of
radius $R_2$, with $R_2 = 1.5$ (circles), $R_2 = 2$ (squares) and $R_3
= 3$ (triangles), obtained by the diagonalization of the matrix $h^2
\Lambda + i\B$ (truncated to sizes $1531\times 1531$ for $R_2 = 1.5$,
$2334 \times 2334$ for $R_2 = 2$, and $2391\times 2391$ for $R_2 =
3$).  Solid line presents the four-terms expansion
(\ref{eq:lambda_app3}) for $\lambda_h^{ND,(1,1)}$, while dashed line
shows its reduced three-terms version (without $h^{\frac 43}$ term). }
\label{fig:annulusND_mu2}
\end{figure}

For the case ND, Fig. \ref{fig:annulusND_mu2} shows the first rescaled
eigenvalue $\lambda^{(1)}_h$ that corresponds to an eigenfunction
which, for small $h$, is localized near the inner circle.  As a
consequence, the asymptotic behavior of $\lambda^{(1)}_h$ as $h\to 0$
is expected to be independent of the outer boundary.  This is indeed
confirmed because the numerical results for three annuli with $R_2 =
1.5\,$, $R_2 = 2$ and $R_3 = 3$ are indistinguishable for $h^{\frac
13}$ smaller than $0.5$.  For comparison, we also plot the four-terms
asymptotics (\ref{eq:lambda_app3}) that we derived for the exterior of
the disk of radius $R_1 = 1$.  One can see that the inclusion of the
term $h^{\frac 43}$ improves the quality of the expansion (as compared
to its reduced three-terms version without $h^{\frac 43}$ term).

\subsubsection{Domain with transmission condition}

Finally, we consider the BT-operator in the union of two subdomains,
the disk $\Omega_- = \{ (x_1,x_2)\in \R^2~:~ |x| < R_1\}$ and the
annulus $\Omega_+ = \{ (x_1,x_2)\in \R^2~:~ R_1 < |x| < R_2\}$
separated by the circle of radius $R_1$ on which the transmission
boundary condition is imposed.  We impose the Dirichlet boundary
condition at the outer boundary of the domain (at the circle of radius
$R_2$) to ensure that first eigenfunctions are localized near points
$(\pm R_1,0)$ with transmission boundary condition.\\

\begin{figure}
\begin{center}
\includegraphics[width=62mm]{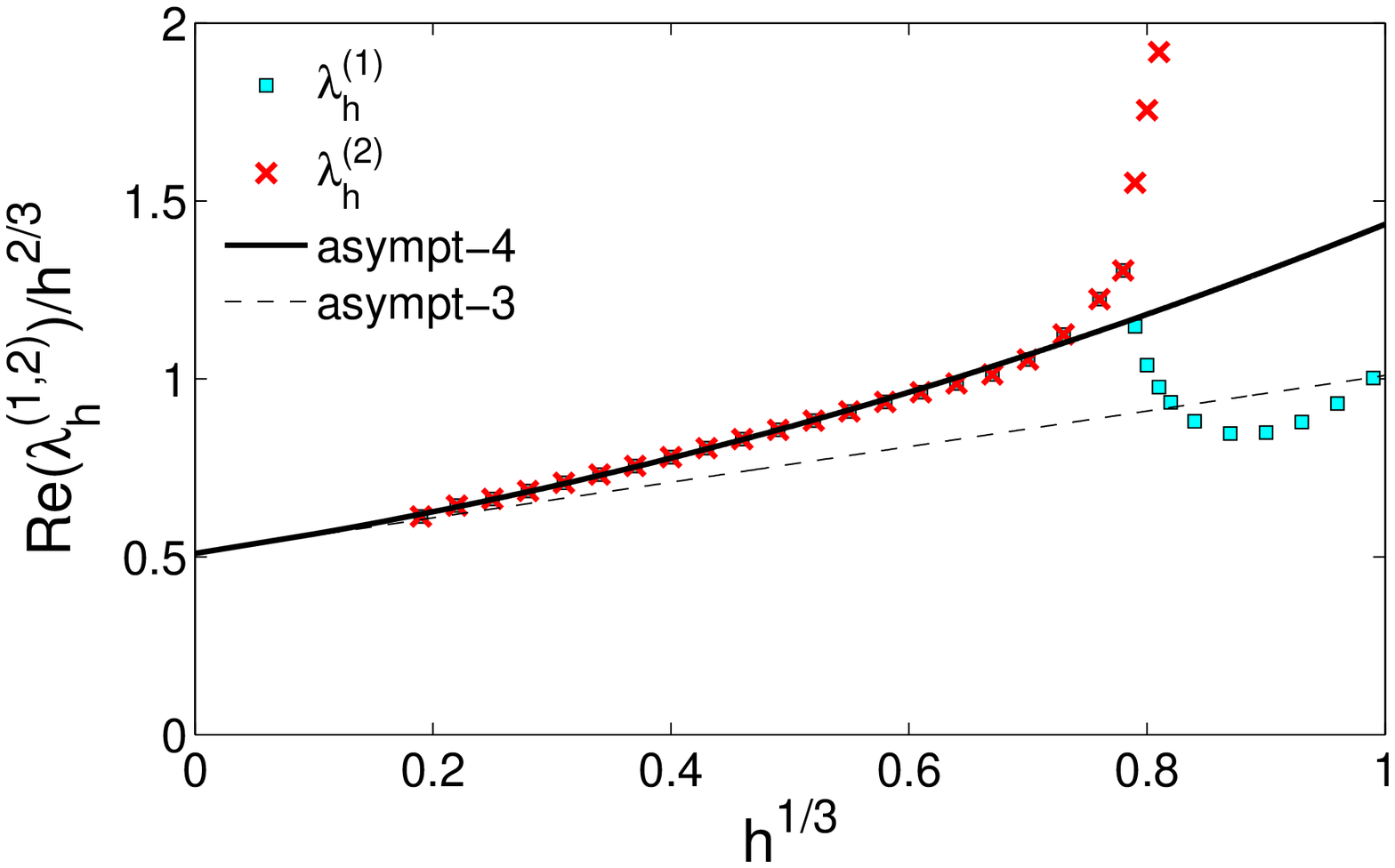}  
\includegraphics[width=62mm]{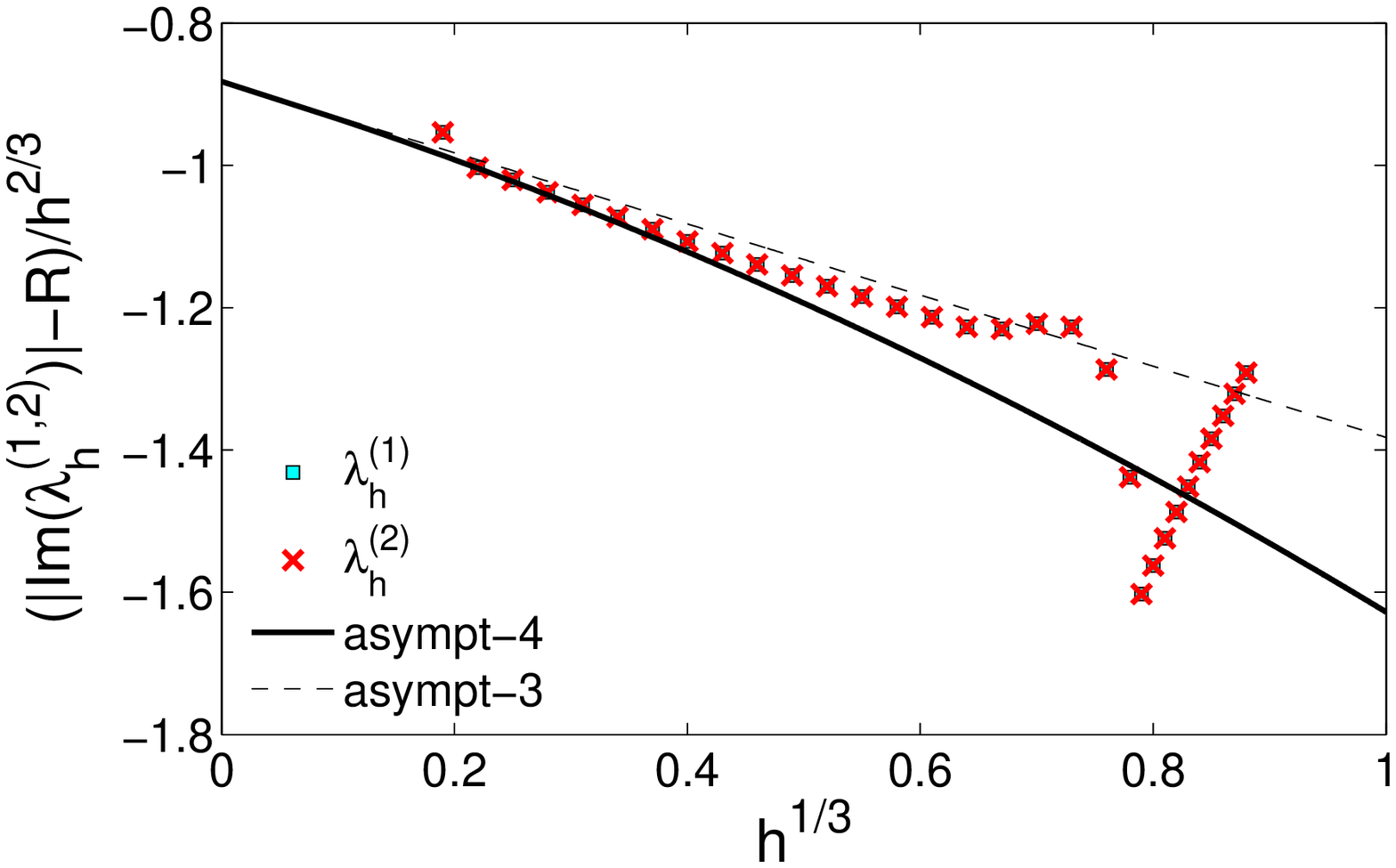}  
\end{center}
\caption{
The rescaled eigenvalues $\lambda^{(1)}_h$ and $\lambda^{(2)}_h$ of
the BT-operator in the union of the disk and annulus with transmission
condition at the inner boundary of radius $R_1 = 1$ (with $\hat \kappa
= 1$) and Dirichlet condition at the outer boundary of radius $R_2 =
2$.  Symbols (squares and crosses) show the numerical results of the
diagonalization of the matrix $h^2 \Lambda + i\B$ (truncated to the
size $3197\times 3197$), solid line presents the four-terms expansion
(\ref{eq:lambda_app2Th}) for $\lambda_h^{T,(1,1)}$, while dashed line
shows its reduced three-terms version (without $h^{\frac 43}$ term). }
\label{fig:twolayers_mu12}
\end{figure}

Figure \ref{fig:twolayers_mu12} shows the rescaled eigenvalues
$\lambda^{(1)}_h$ and $\lambda^{(2)}_h$ of the BT-operator with a
fixed $\hat \kappa = 1$ and $\kappa$ scaling as $\hat\kappa h^{\frac
23}$.  As in earlier examples, the first two eigenvalues are complex
conjugate to each other for small $h$ but they split at larger $h$.
One can see that the asymptotic relation (\ref{eq:lambda_app2Th}) with
$n=k=1$ accurately describes the behavior of these eigenvalues for
small $h$.  \\

\begin{figure}
\begin{center}
\includegraphics[width=62mm]{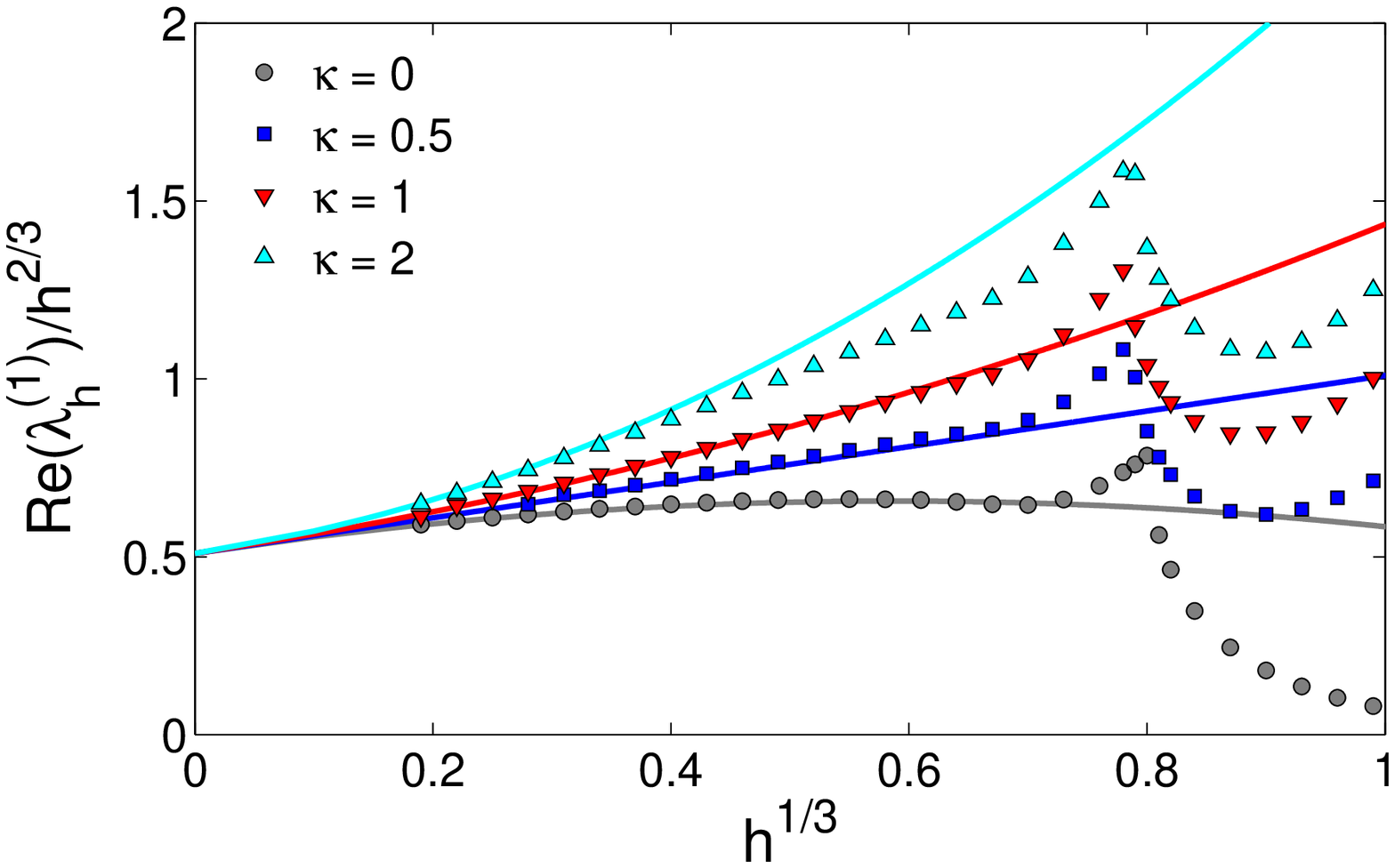} 
\includegraphics[width=62mm]{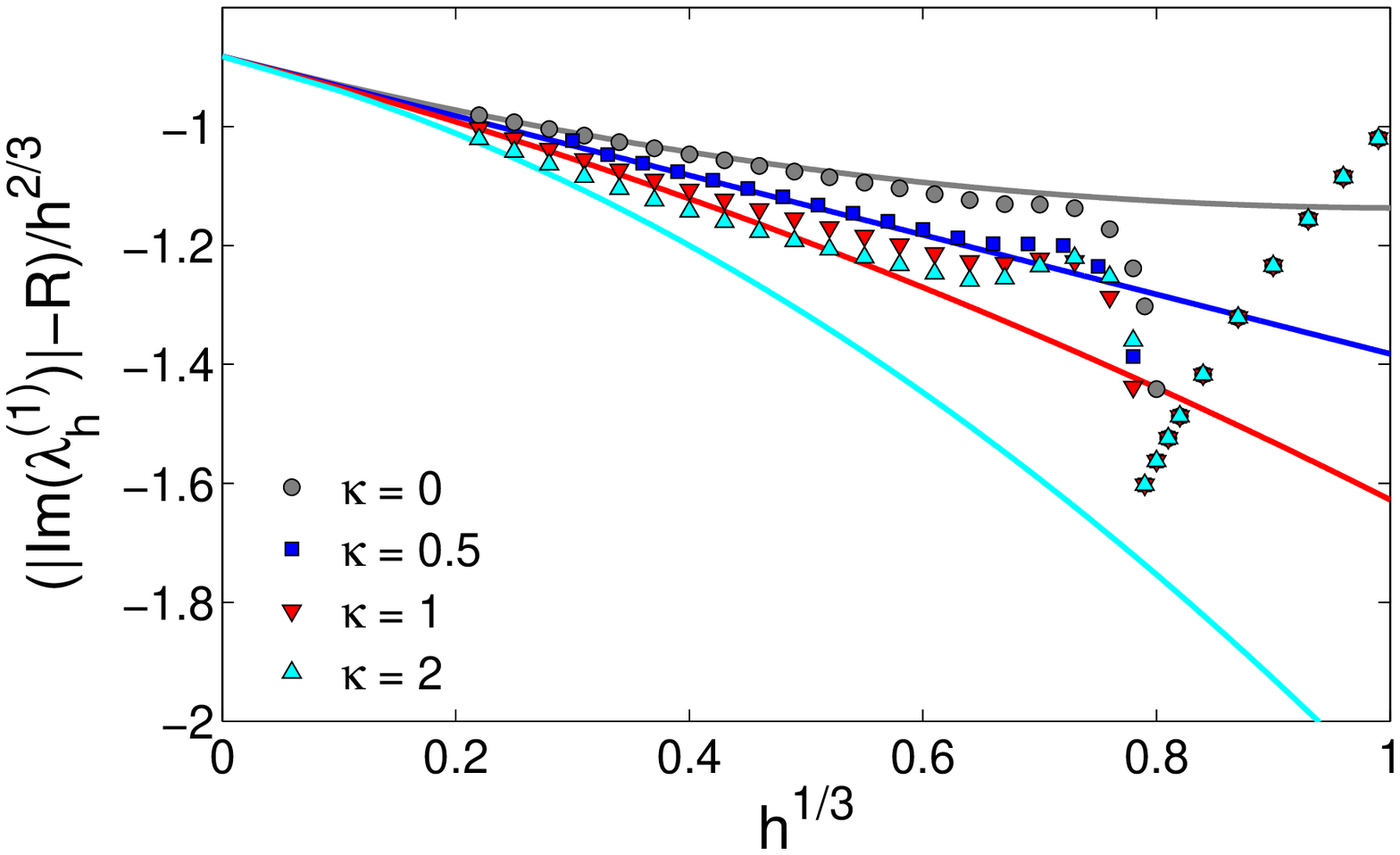} 
\end{center}
\caption{
The rescaled eigenvalue $\lambda^{(1)}_h$ of the BT-operator in the
union of the disk and annulus with transmission condition at the inner
boundary of radius $R_1 = 1$ (with several values of $\hat \kappa$:
$0\,$, $0.5\,$, $1\,$, $2$) and Dirichlet condition at the outer boundary of
radius $R_2 = 2$.  Symbols (circles, squares, triangles) show the
numerical results of the diagonalization of the truncated matrix $h^2
\Lambda + i\B$, solid lines present the four-terms expansion
(\ref{eq:lambda_app2Th}) for $\lambda_h^{T,(1,1)}$. }
\label{fig:twolayers_mu1_kappa}
\end{figure}

Figure \ref{fig:twolayers_mu1_kappa} shows the first rescaled
eigenvalue $\lambda^{(1)}_h$ for several values of $\hat \kappa$ (with
$\kappa$ scaling as $\hat\kappa h^{\frac23}$).  In the special case
$\hat \kappa = 0$, the two subdomains are separated from each other by
Neumann boundary condition, and the spectrum of the BT operator is
obtained from its spectra for each subdomain.  As a consequence, we
plot in this case the first rescaled eigenvalue for the BT operator in
the unit disk with Neumann boundary condition (as in Fig.
\ref{fig:diskN_mu12}).  One can see that the expansion
(\ref{eq:lambda_app2Th}) accurately captures the asymptotic behavior.
We recall that the transmission parameter $\hat\kappa$ appears only in
the fourth term of order $h^{\frac 43}$.  Note also that this term
vanishes in the case $\hat \kappa = 1/2$ as two contributions in
(\ref{eq:lambda_app2Th}) compensate each other.

\subsection{Eigenfunctions}

\begin{figure}
\begin{center}
\includegraphics[width=62mm]{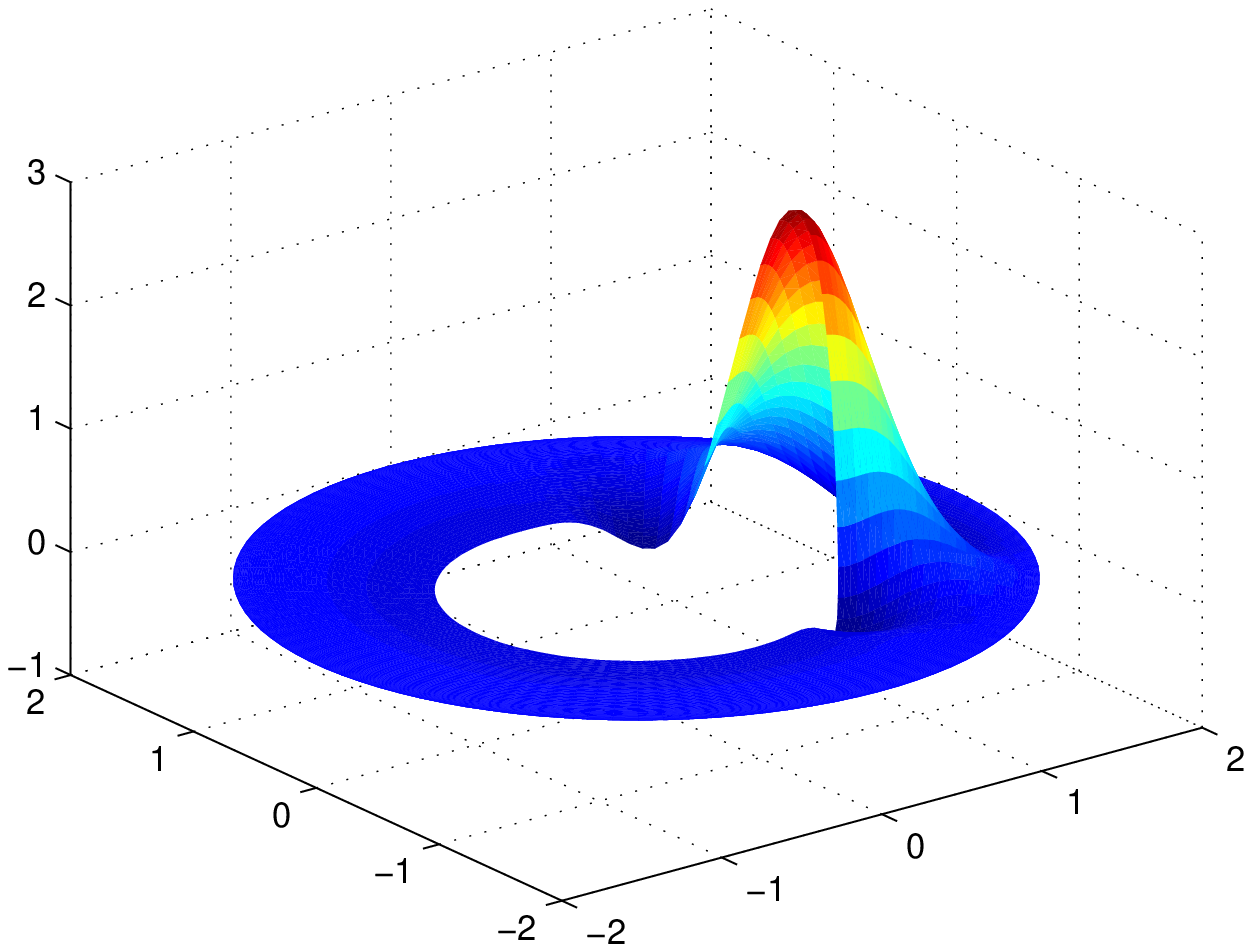} 
\includegraphics[width=62mm]{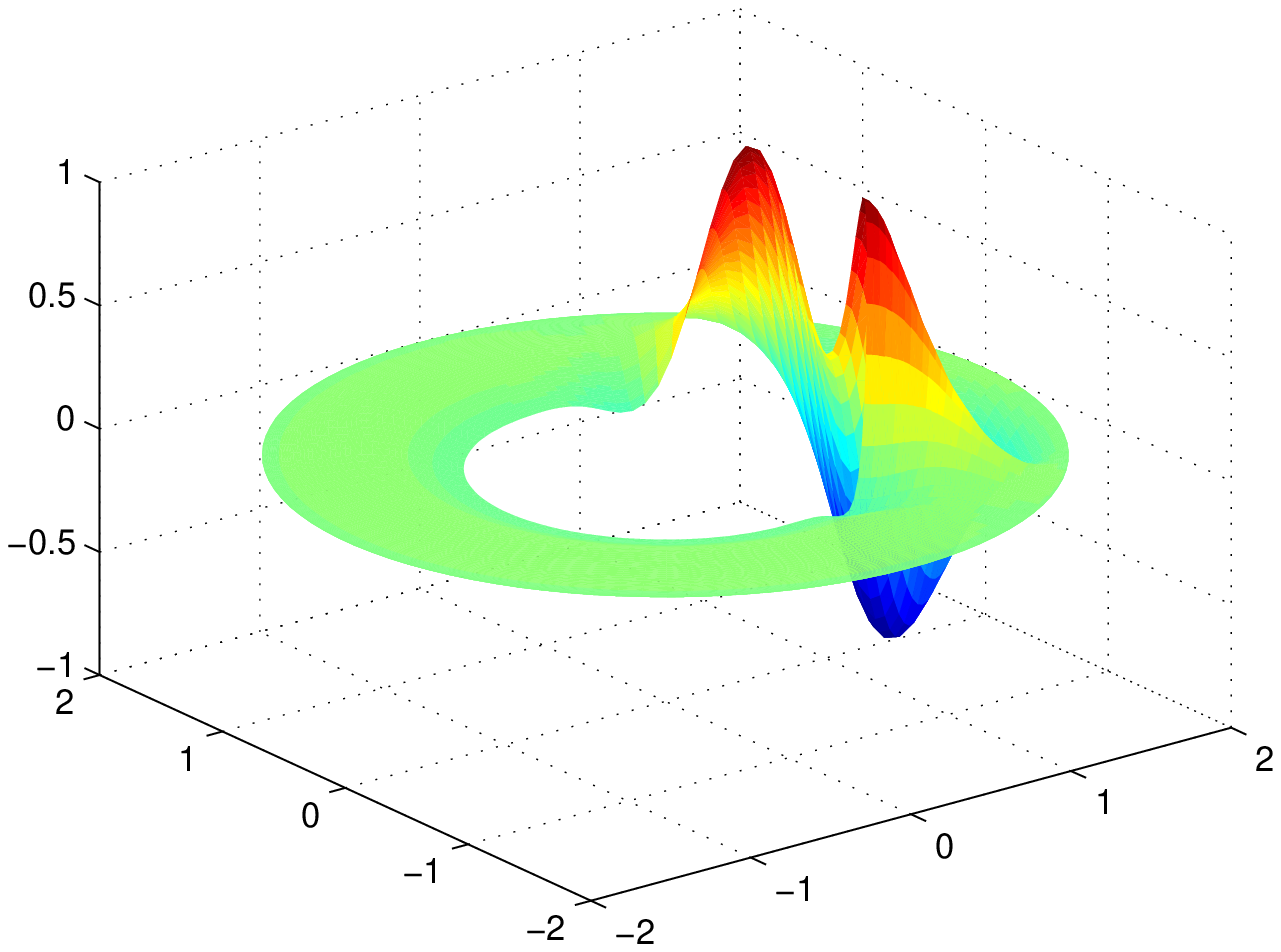} 
\includegraphics[width=62mm]{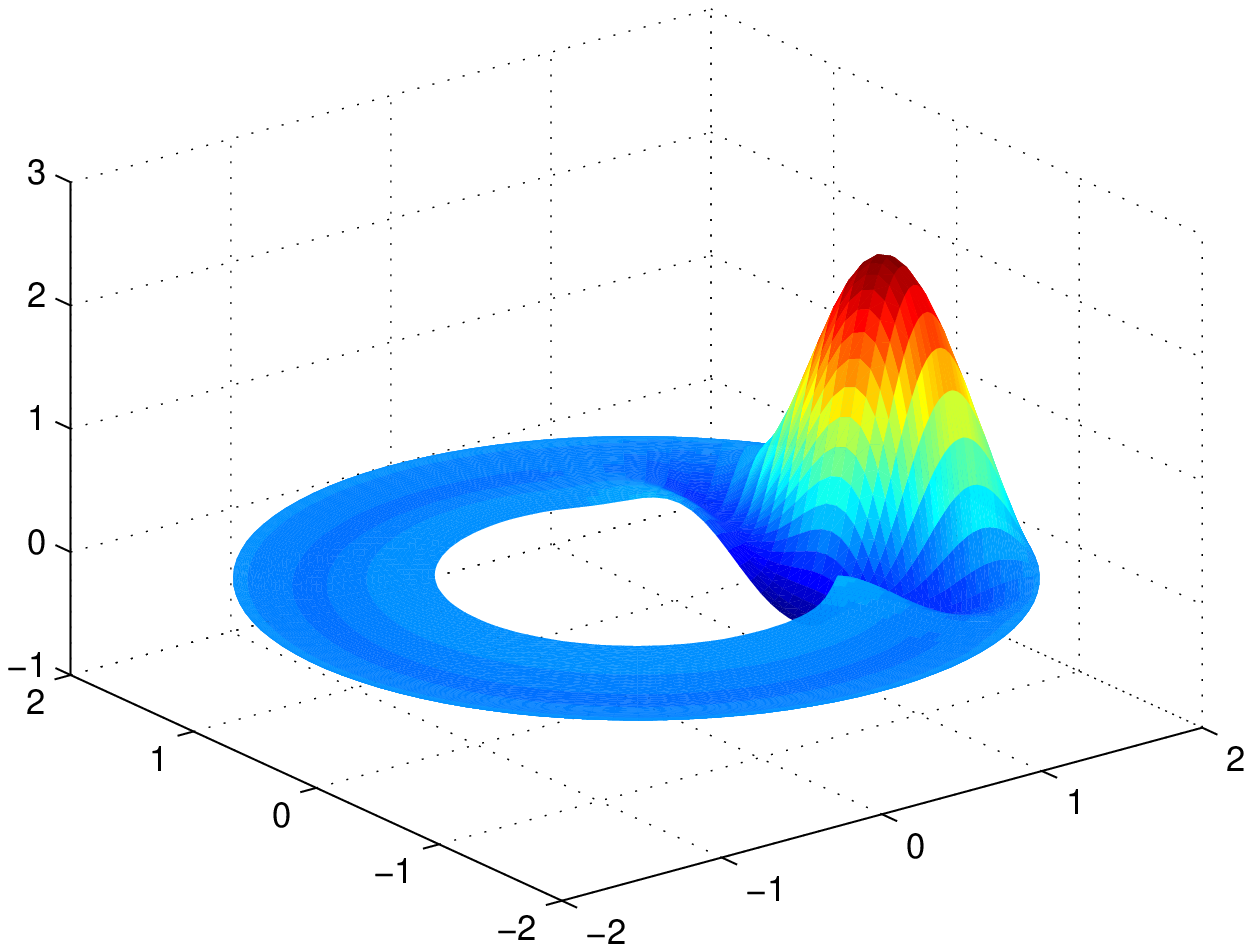} 
\includegraphics[width=62mm]{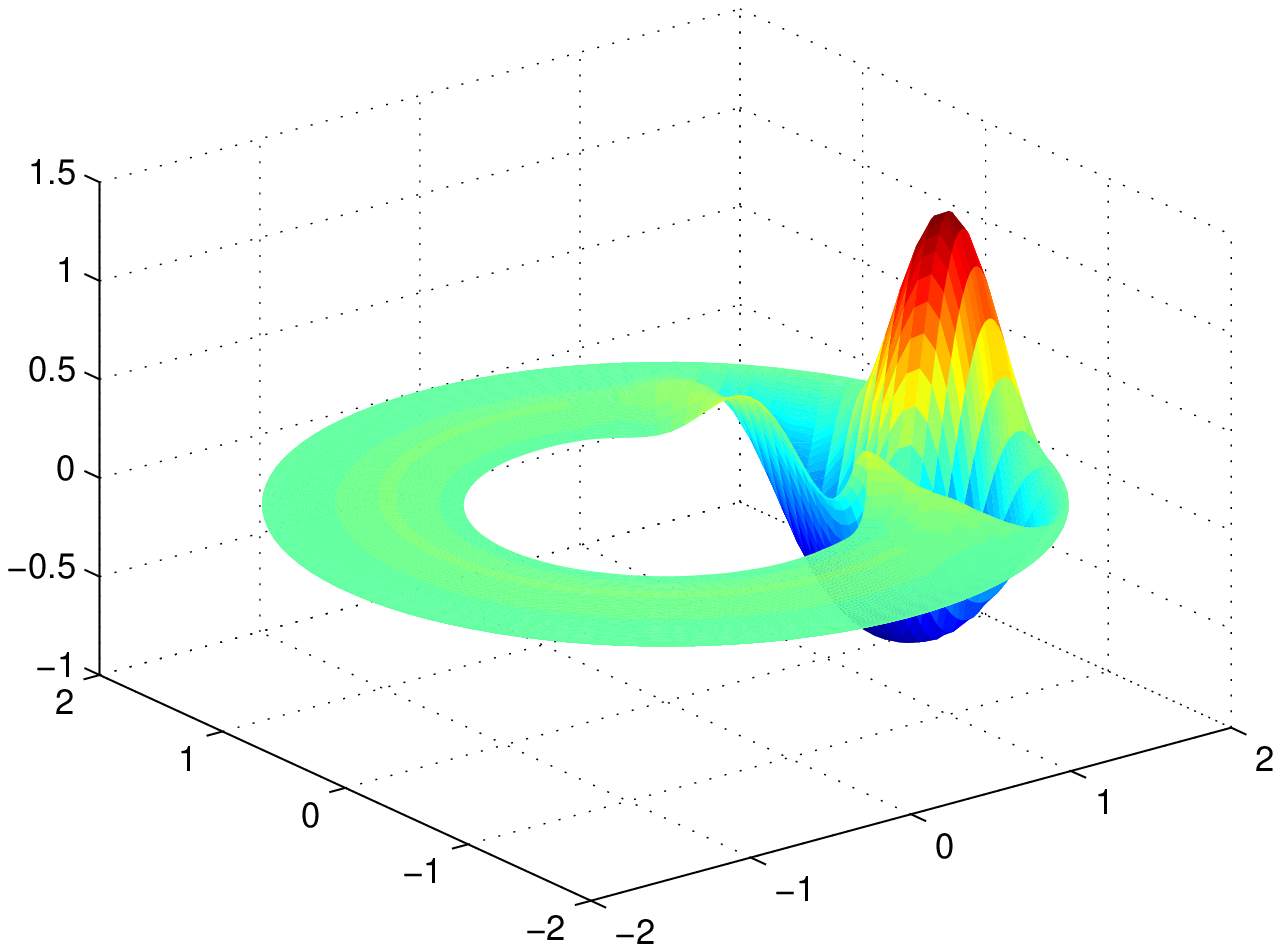} 
\noindent\makebox[\linewidth]{\rule{124mm}{0.4pt}}
\includegraphics[width=62mm]{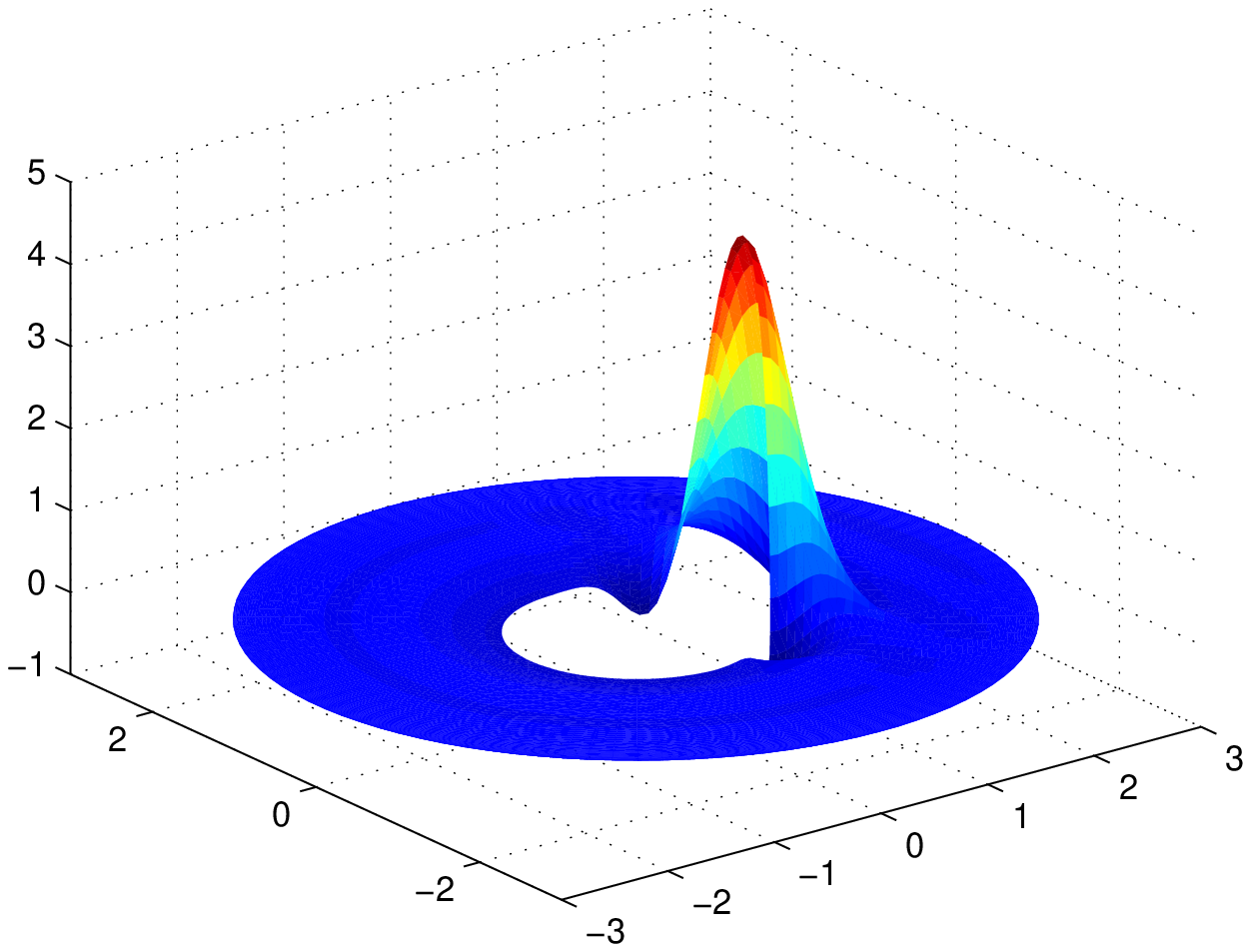} 
\includegraphics[width=62mm]{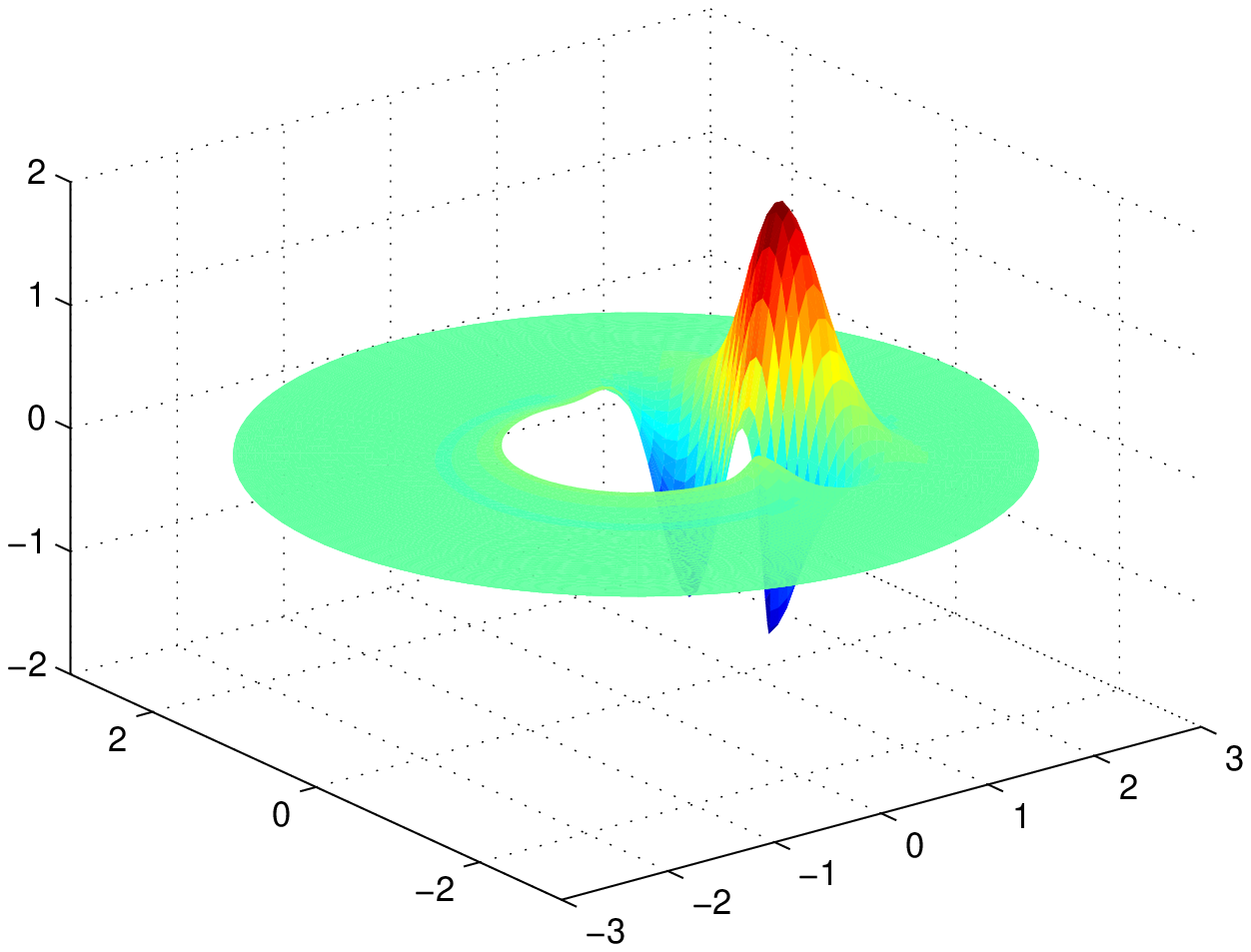} 
\includegraphics[width=62mm]{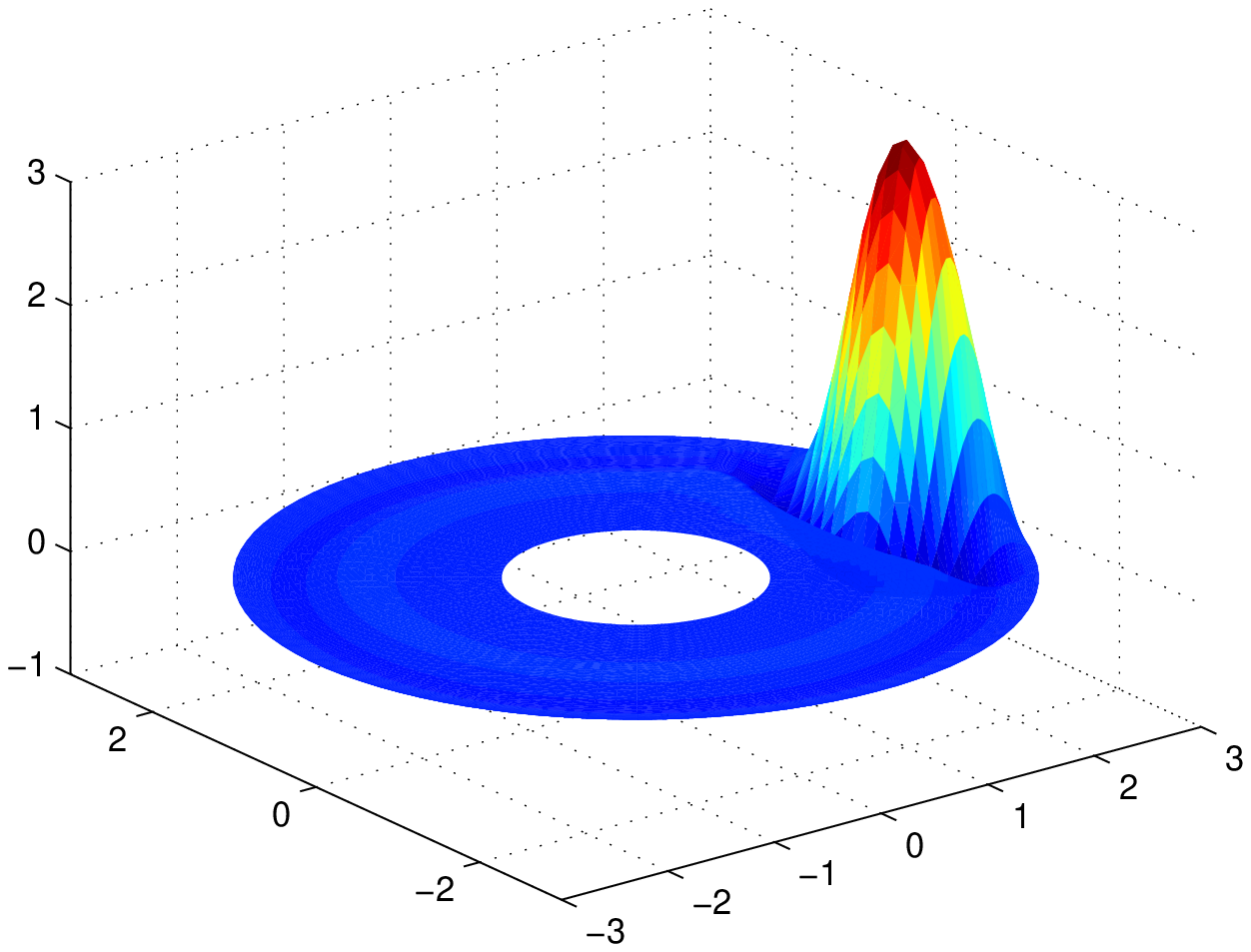} 
\includegraphics[width=62mm]{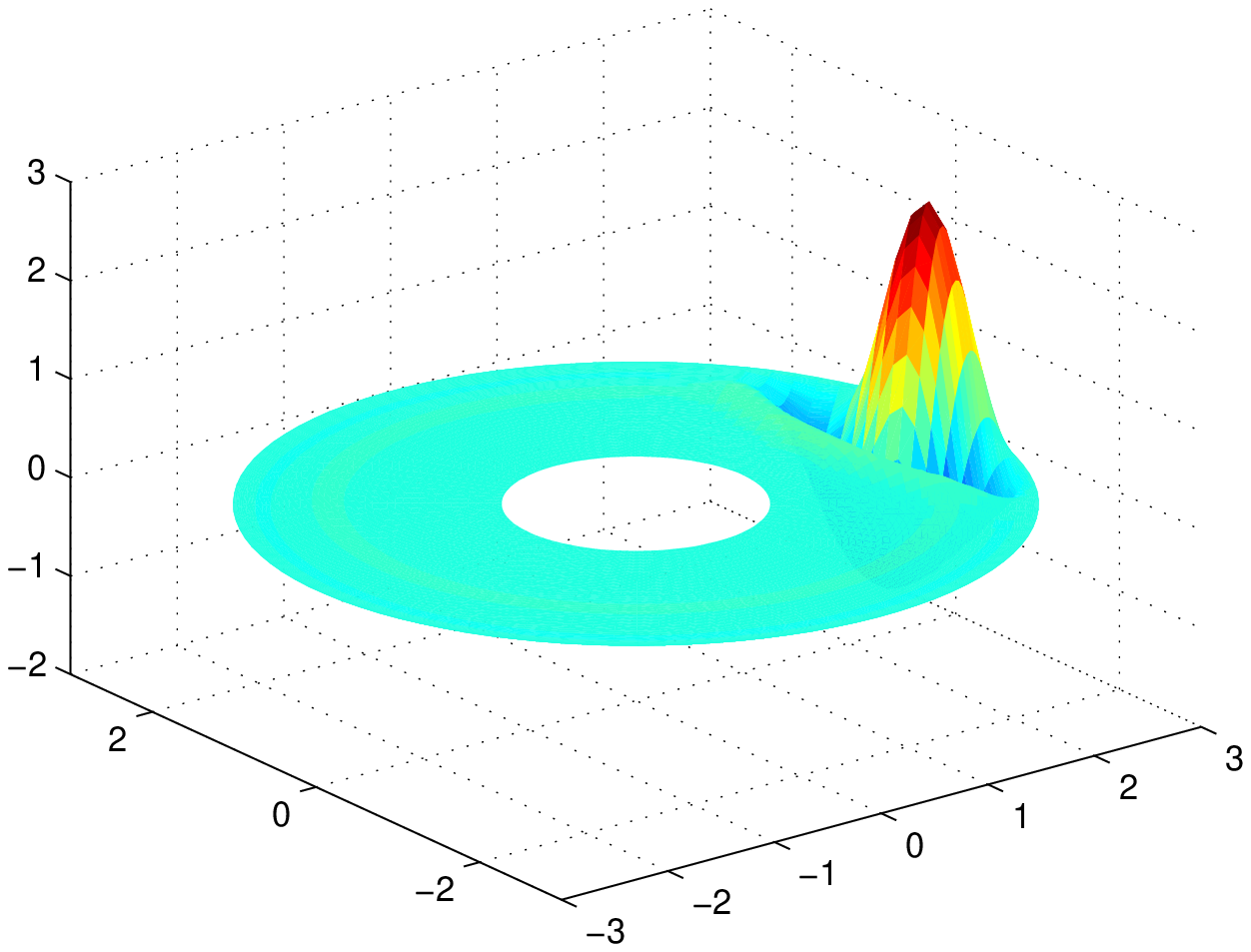} 
\end{center}
\caption{
Real (left) and imaginary (right) parts of the eigenfunctions
$u^{(1)}_h$ (top) and $u^{(3)}_h$ (bottom) at $h = 0.1$ for the
annulus with Neumann boundary condition at the inner circle of radius
$R_1 = 1$ and Dirichlet boundary condition at the outer circle of
radius $R_2 = 2$ (four plots above horizontal line) or $R_2 = 3$ (four
plots below horizontal line).  Numerical computation is based on the
truncated matrix representation of sizes $2334\times 2334$ and
$2391\times 2391$, respectively. }
\label{fig:annulus_eigenfunctions}
\end{figure}

For the annulus with Neumann boundary condition at the inner circle of
radius $R_1 = 1$ and Dirichlet boundary condition at the outer circle
of radius $R_2 = 2\,$, Fig. \ref{fig:annulus_eigenfunctions}(top)
shows two eigenfunctions of the BT operator with $h = 0.1$
(corresponding to $h^{\frac 13} \approx 0.4642$).  One can already
recognize the localization of the first eigenfunction $u^{(1)}_h$ at
the inner boundary, while the eigenfunction $u^{(3)}_h$ tends to
localize near the outer boundary.  Their pairs $u^{(2)}_h$ and
$u^{(4)}_h$ (not shown) exhibit the same behavior near the opposite
points $(-R_1,0)$ and $(-R_2,0)$, respectively.  Since $h = 0.1$ is
not small enough, the localization becomes less and less marked for
other eigenfunctions which progressively spread over the whole annulus
(not shown).  For comparison, we also plot in
Fig. \ref{fig:annulus_eigenfunctions}(bottom) the eigenfunctions
$u^{(1)}_h$ and $u^{(3)}_h$ for a thicker annulus of outer radius $R_2
= 3$.  One can see that these eigenfunctions look very similar to that
of the annulus with $R_2 = 2$.

\begin{figure}
\begin{center}
\includegraphics[width=62mm]{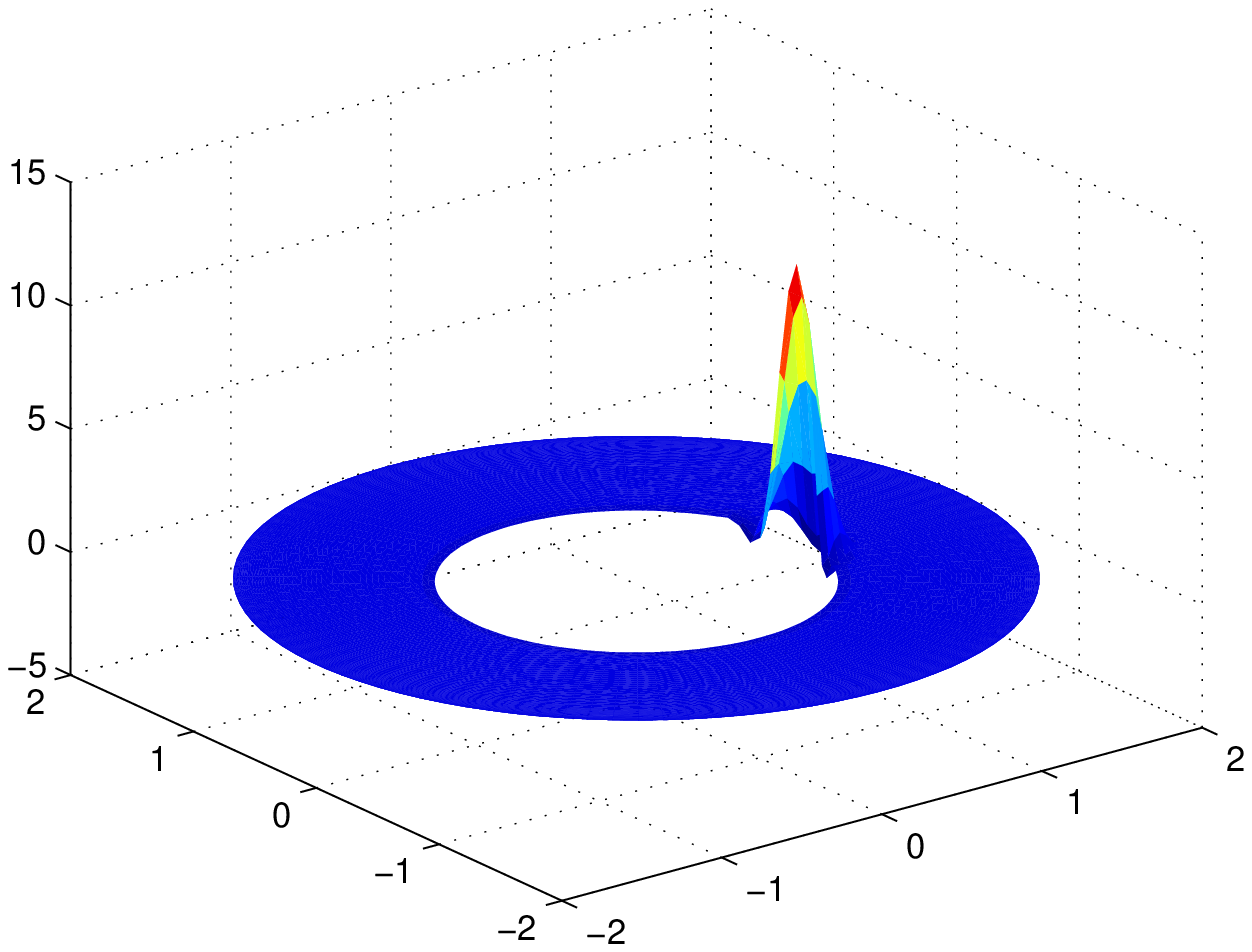} 
\includegraphics[width=62mm]{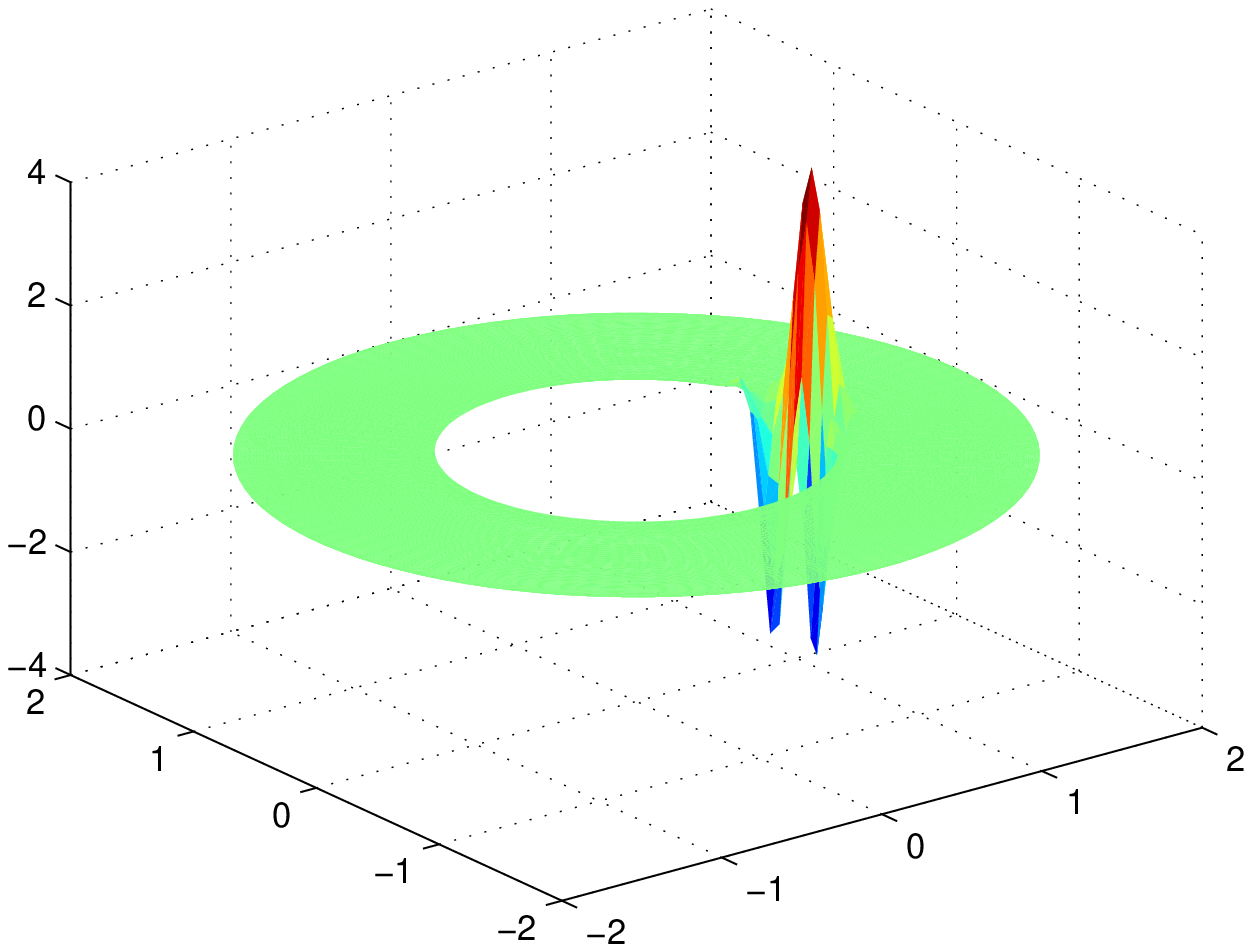} 
\includegraphics[width=62mm]{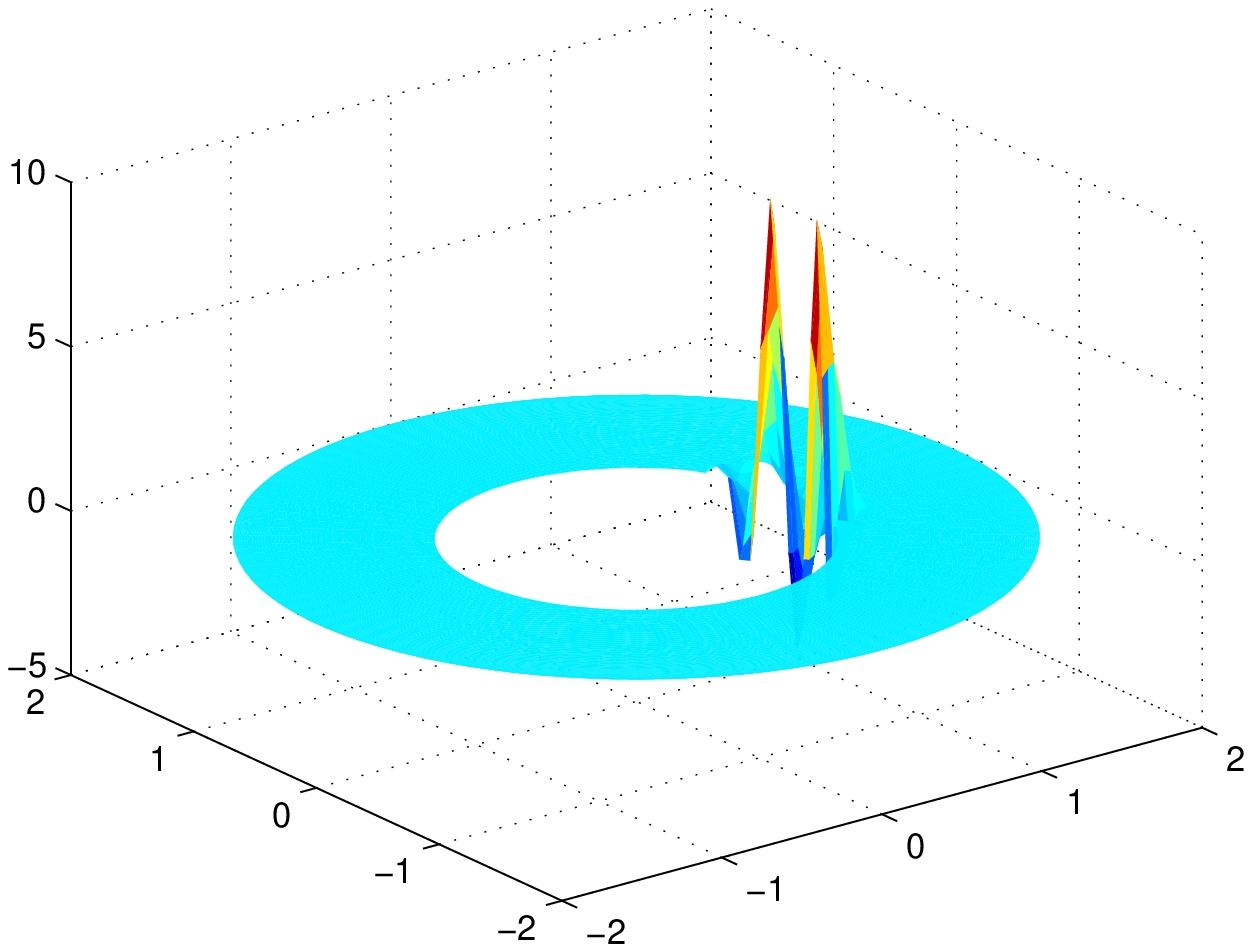} 
\includegraphics[width=62mm]{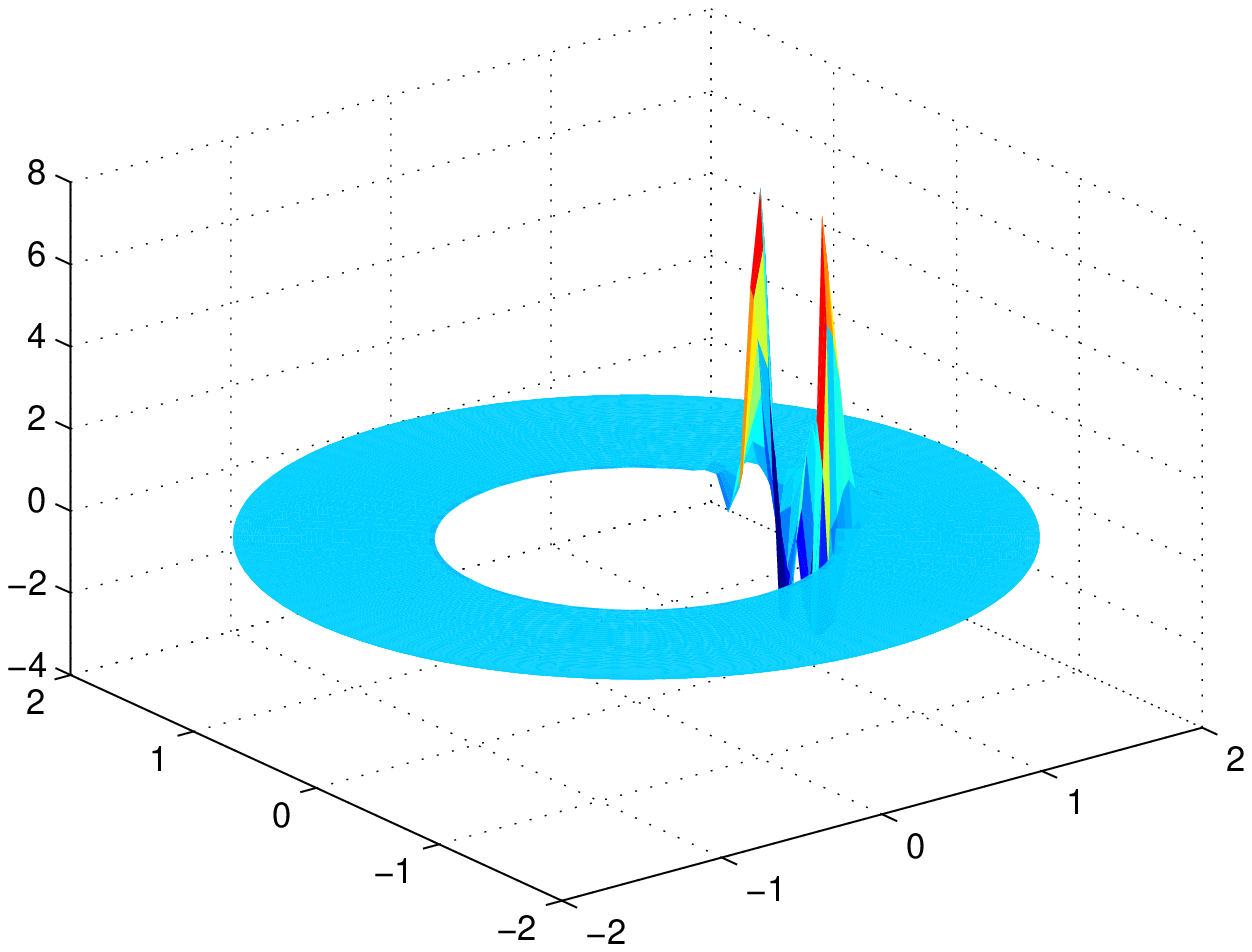} 
\includegraphics[width=62mm]{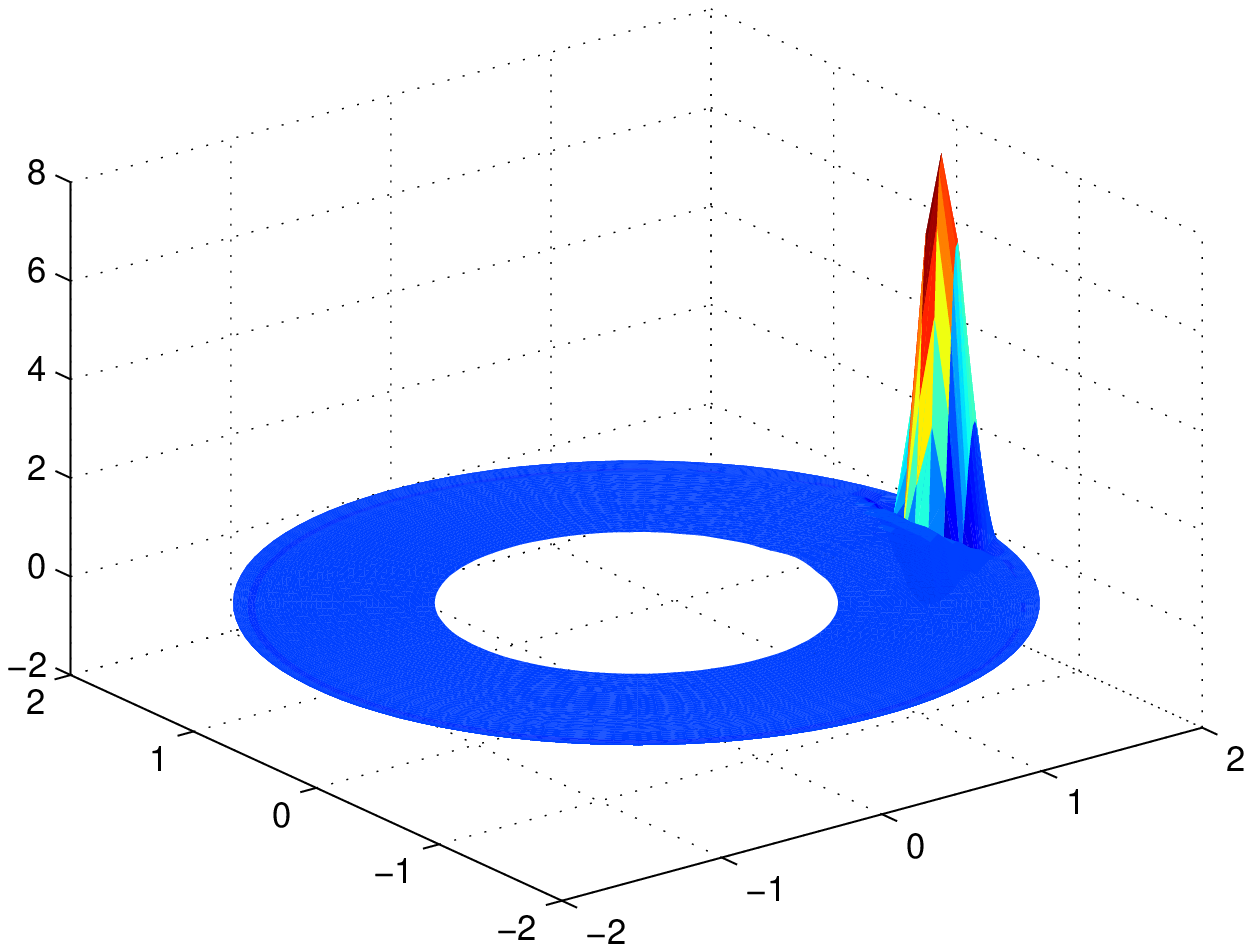} 
\includegraphics[width=62mm]{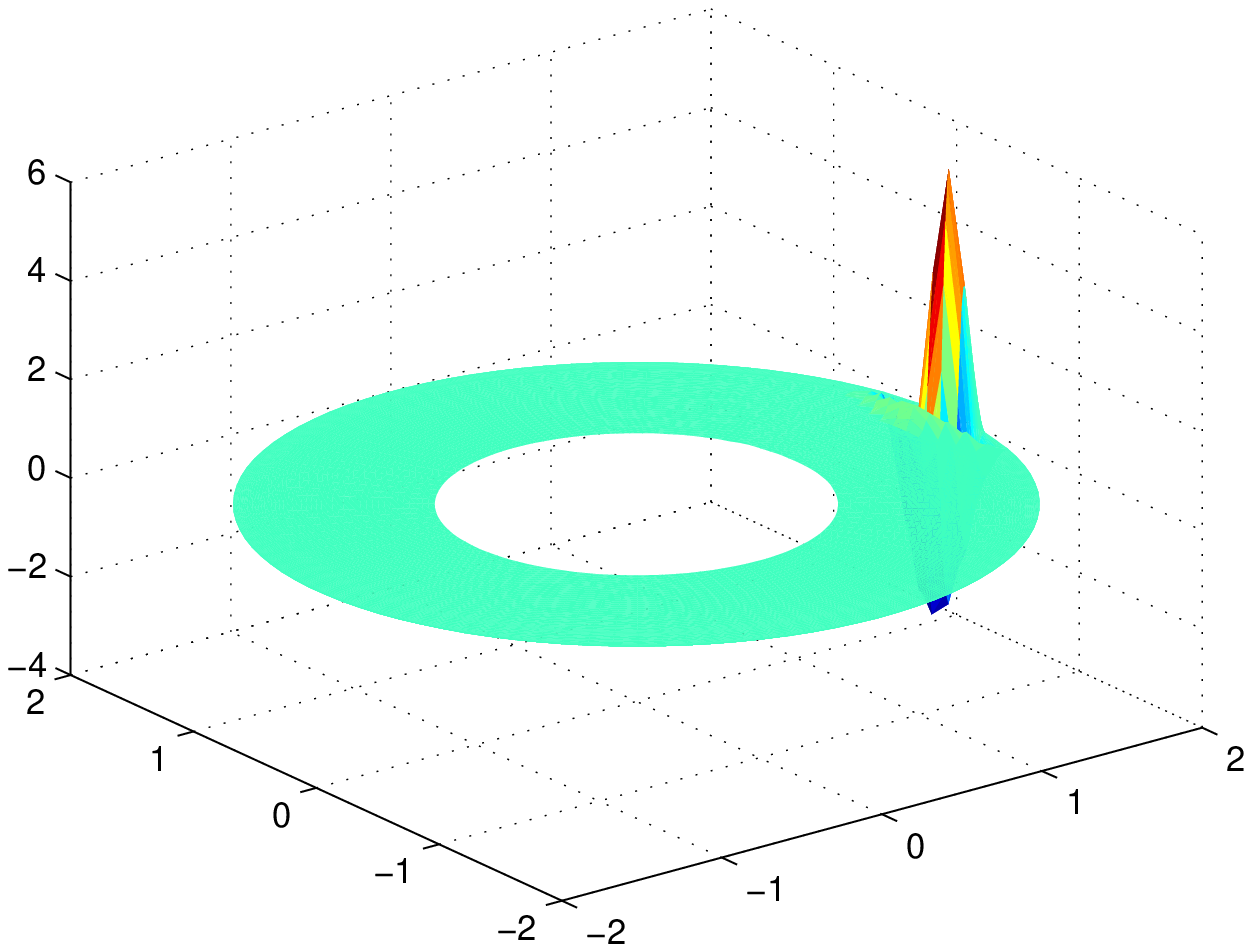} 
\includegraphics[width=62mm]{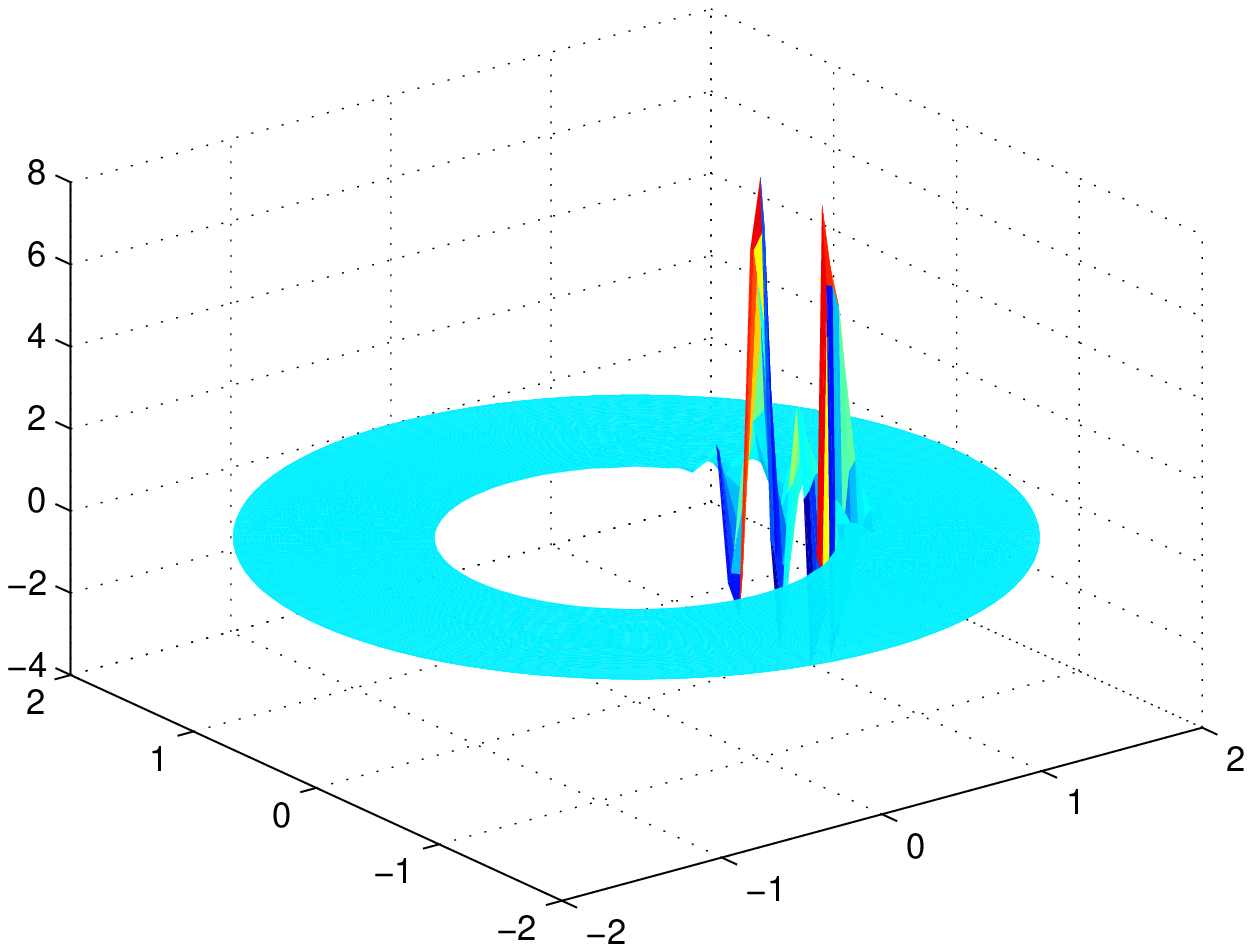} 
\includegraphics[width=62mm]{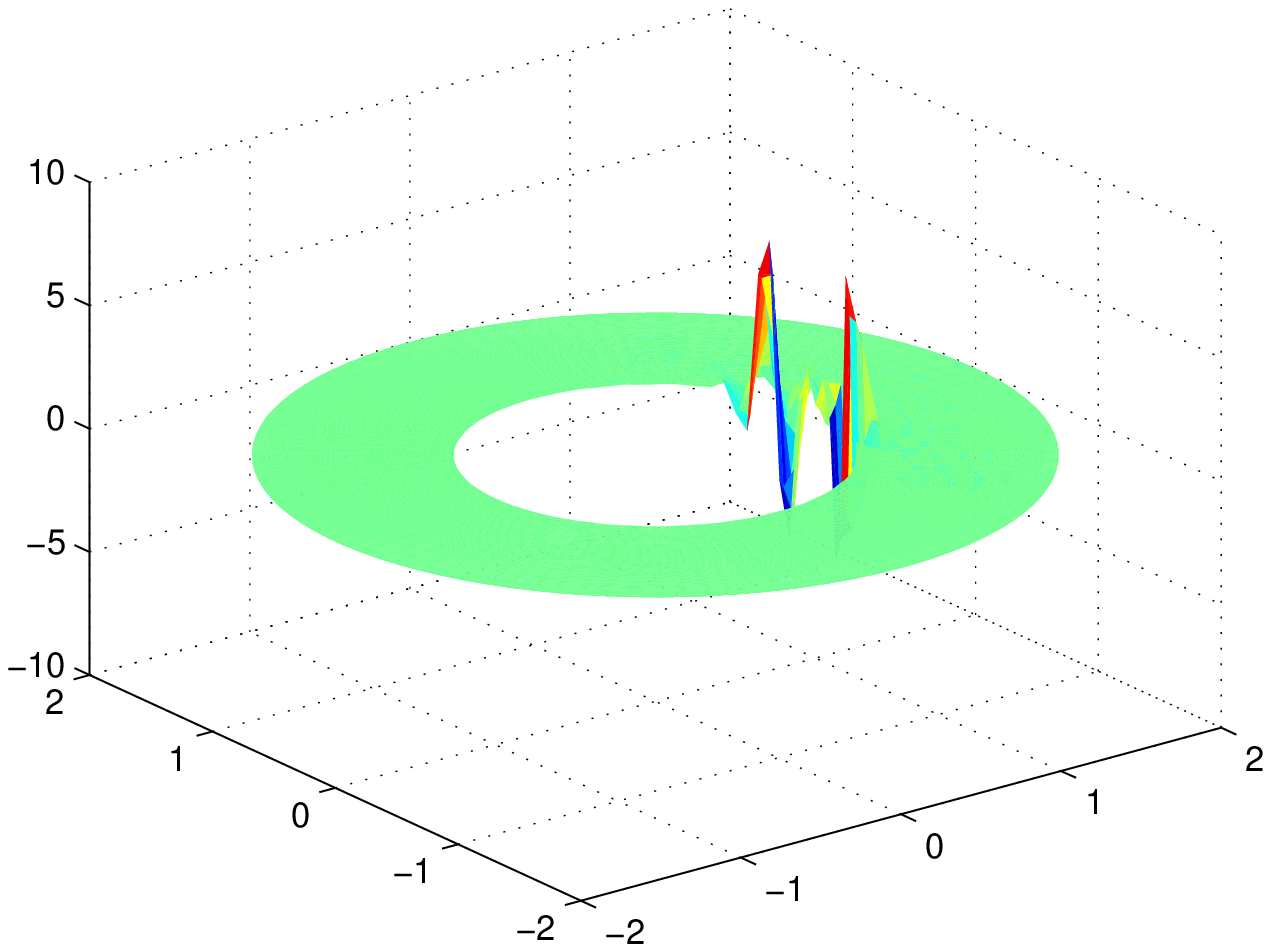} 
\end{center}
\caption{
Real (left) and imaginary (right) parts of the eigenfunctions
$u^{(1)}_h$ (top), $u^{(3)}_h$, $u^{(5)}_h$ and $u^{(7)}_h$ (bottom)
at $h = 0.01$ for the annulus with Neumann boundary condition on the
inner circle of radius $R_1 = 1$ and Dirichlet boundary condition on
the outer circle of radius $R_2 = 2$ (numerical computation based on
the truncated matrix representation of size $2334\times 2334$). }
\label{fig:annulus_eigenfunctions2}
\end{figure}

For smaller $h = 0.01$ (corresponding to $h^{\frac 13} \approx
0.2154$), the localization of eigenfunctions is much more pronounced.
Figure \ref{fig:annulus_eigenfunctions2} shows four eigenfunctions for
the annulus of radii $R_1 = 1$ (Neumann condition) and $R_2 = 2$
(Dirichlet condition).  One can see that the eigenfunctions
$u^{(1)}_h$, $u^{(3)}_h$, and $u^{(7)}_h$ are localized near the inner
circle while $u^{(5)}_h$ is localized near the outer circle.  When the
outer circle is moved away, the former eigenfunctions remain almost
unchanged, suggesting that they would exist even in the limiting
domain with $R_2 = \infty$, i.e., in the complement of the unit disk.
In turn, the eigenfunctions that are localized near the outer boundary
(such as $u^{(5)}_h$) will be eliminated.  In spite of this numerical
evidence, the existence of eigenfunctions of the BT operator for
unbounded domains remains conjectural.

\begin{figure}
\begin{center}
\includegraphics[width=62mm]{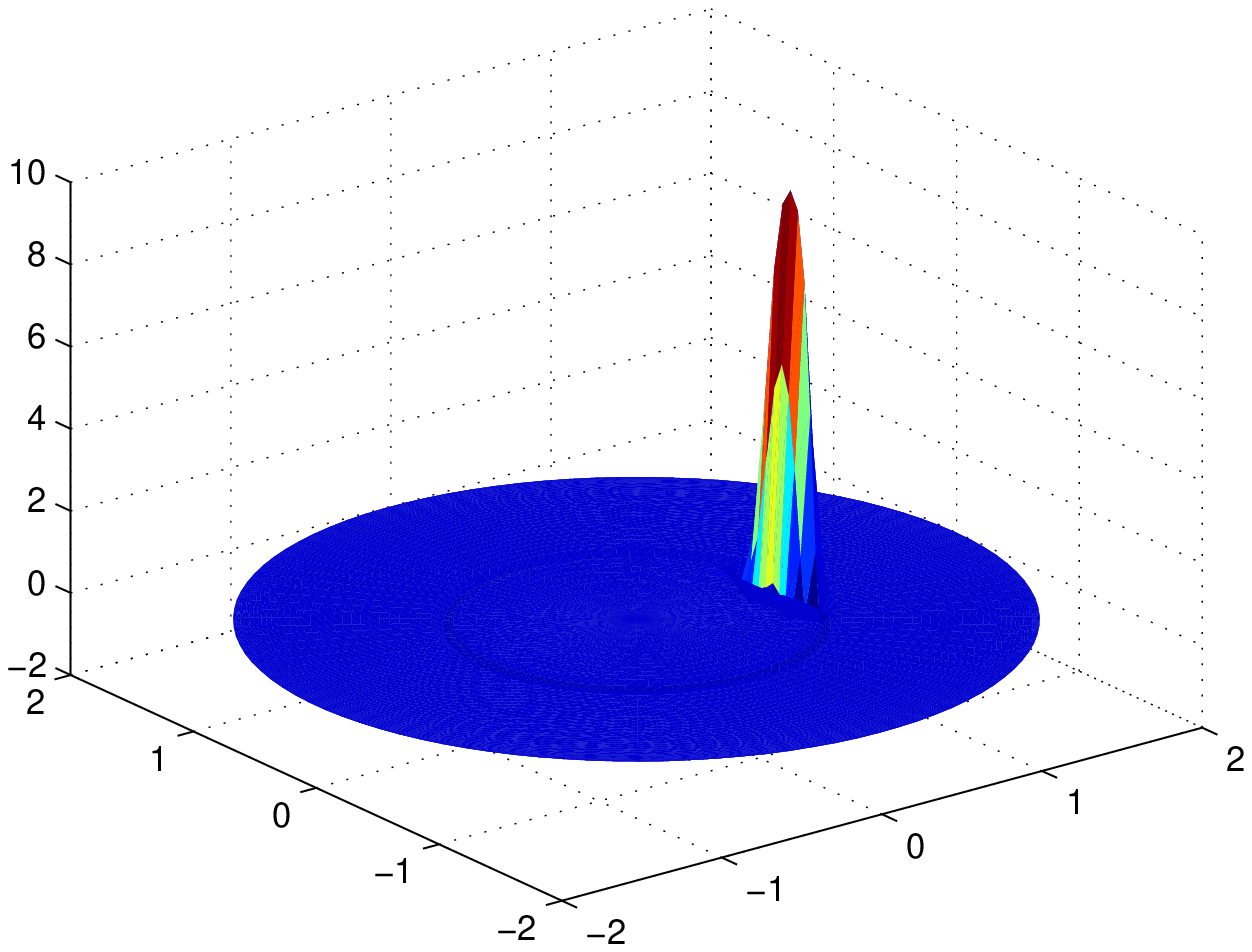} 
\includegraphics[width=62mm]{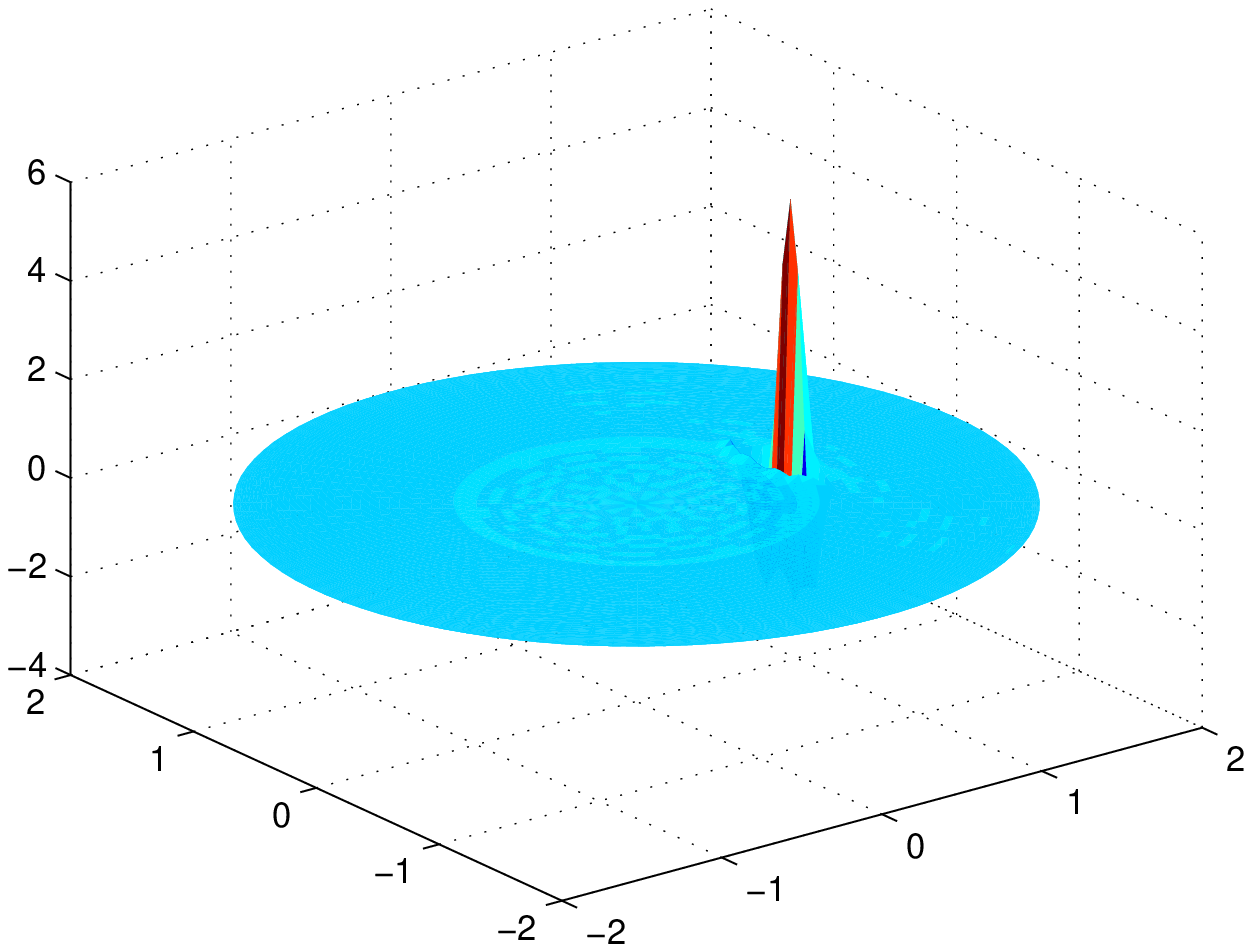} 
\includegraphics[width=62mm]{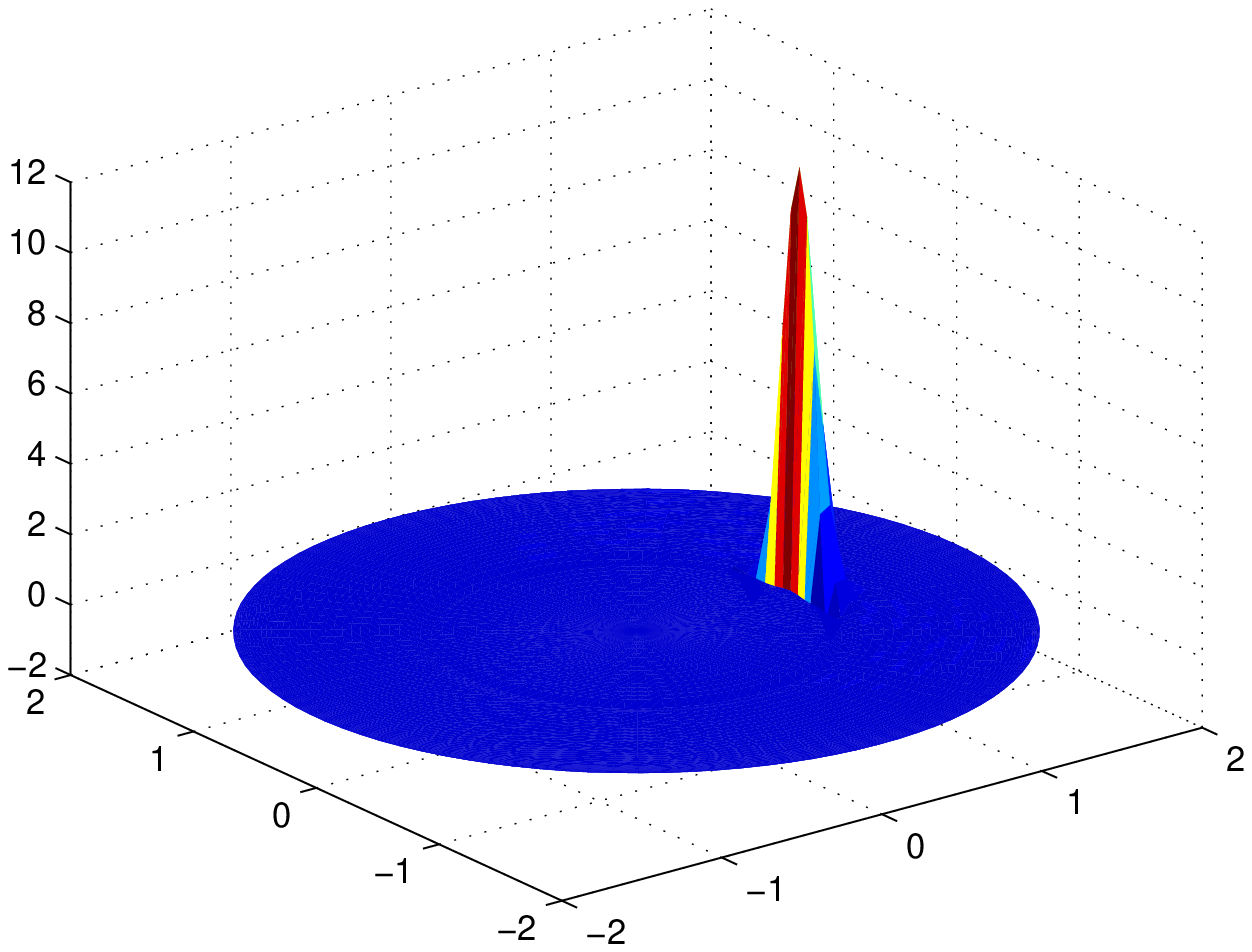} 
\includegraphics[width=62mm]{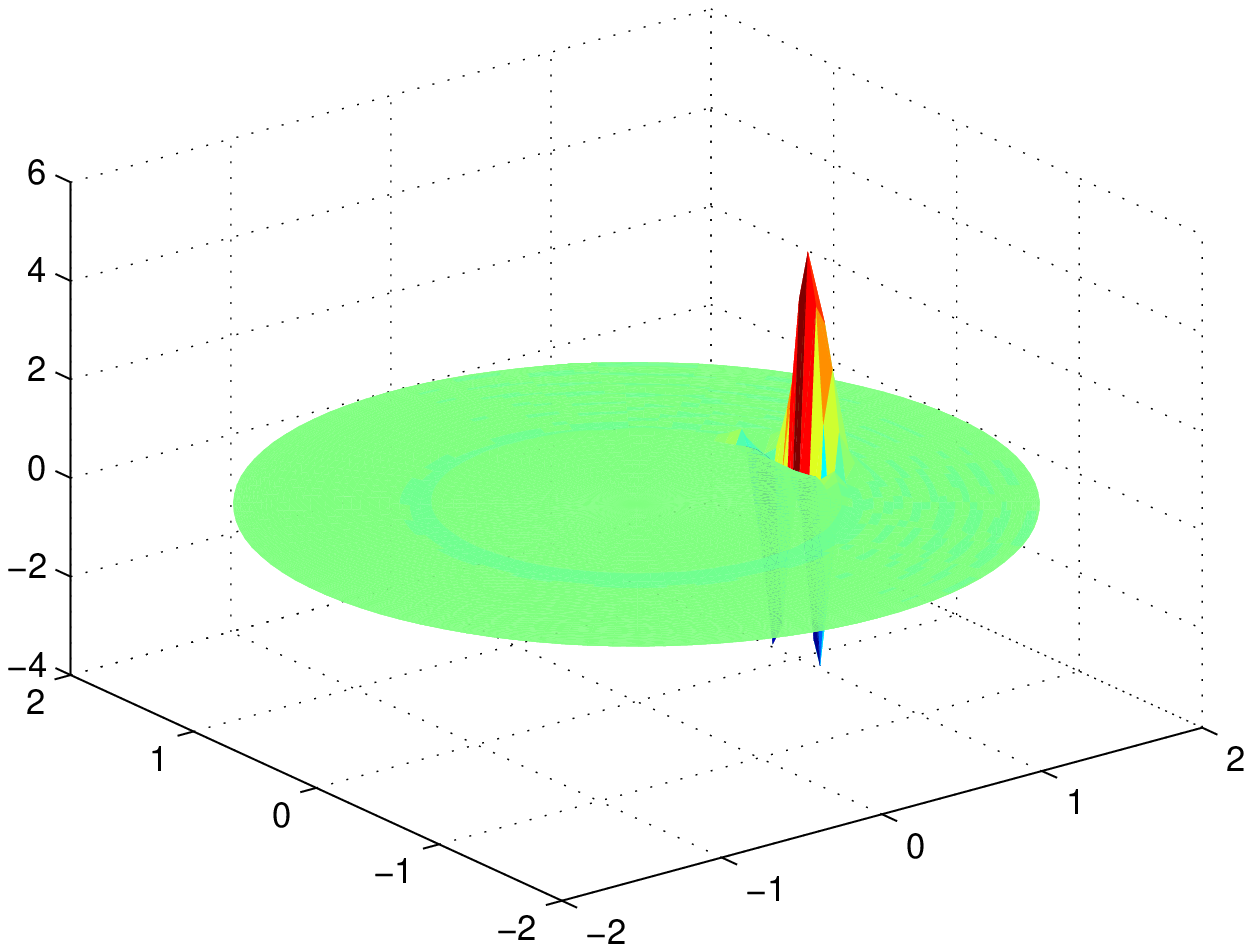} 
\end{center}
\caption{ 
Real (left) and imaginary (right) parts of the eigenfunctions
$u^{(1)}_h$ (top) and $u^{(3)}_h$ (bottom) at $h = 0.01$ for the union
of the disk and annulus with a transmission boundary condition (with
$\hat\kappa = 1$ and $\kappa = \hat\kappa h^{\frac23}$) at the inner
circle of radius $R_1 = 1$ and Dirichlet boundary condition at the
outer circle of radius $R_2 = 2$ (numerical computation based on the
truncated matrix representation of size $3197\times 3197$). }
\label{fig:twolayersTD_1_2_eigenfunctions}
\end{figure}

Figure \ref{fig:twolayersTD_1_2_eigenfunctions} shows the
eigenfunctions $u^{(1)}_h$ and $u^{(3)}_h$ at $h = 0.01$ for the union
of the disk and annulus with transmission condition at the inner
boundary of radius $R_1 = 1$ (with $\hat \kappa = 1$ and $\kappa =
\hat\kappa h^{\frac23}$) and Dirichlet condition at the outer boundary
of radius $R_2 = 2$.  Both eigenfunctions are localized near the inner
boundary.  Moreover, a careful inspection of this figure shows that
$u^{(1)}_h$ is mainly supported by the disk and vanishes rapidly on
the other side of the inner circle (i.e., in the annulus side), while
$u^{(3)}_h$ exhibits the opposite (i.e., it is localized in the
annulus).  This is a new feature of localization as compared to the
one-dimensional case studied in \cite{Gr1,GHH} because the curvature
has the opposite signs on two sides of the boundary.

\begin{figure}
\begin{center}
\includegraphics[width=62mm]{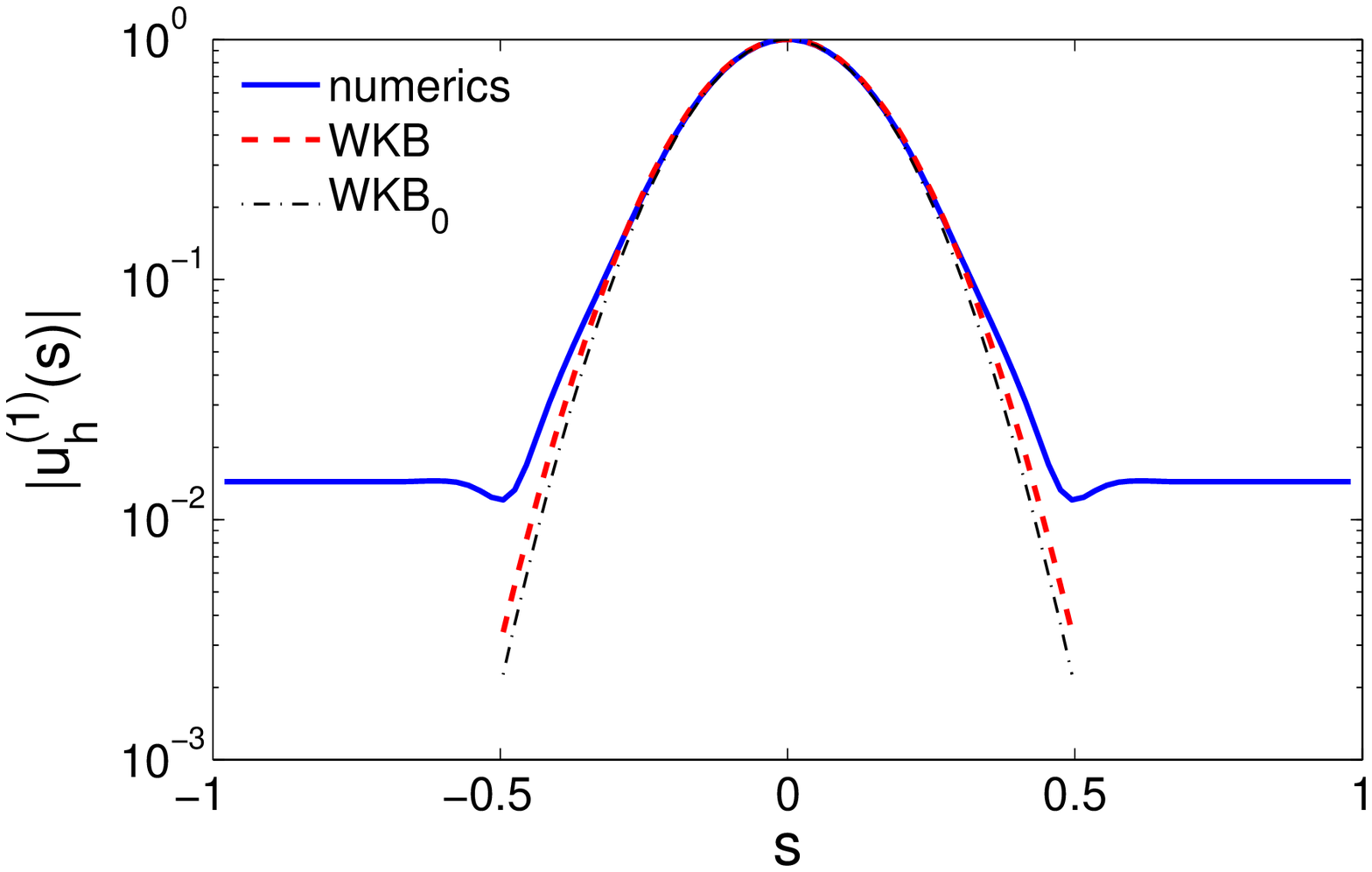} 
\end{center}
\caption{ 
The absolute value of the first eigenfunction $u^{(1)}_h(r,s)$ (solid
line) at $h = 0.01$ and $r = 1$ for the unit disk with Neumann
boundary condition, near the boundary point $s = 0$.  For convenience,
$u^{(1)}_h(r,s)$ is normalized by its maximum at $s=0$.  For
comparison, the absolute value of the WKB approximation,
$\exp(-(\theta_0(s) + h^{\frac 23}\theta_1(s))/h)$ and of its reduced
version, $\exp(-\theta_0(s)/h)$, are shown by dashed and dash-dotted
lines, respectively. }
\label{fig:diskN_WKB}
\end{figure}

Finally, we check the accuracy of the WKB approximation of the first
eigenfunction $u^{(1)}_h$ for the unit disk with Neumann boundary
condition.  To make the illustration easier, we plot in Figure
\ref{fig:diskN_WKB} the absolute value of $u^{(1)}_h$ at $h = 0.01$,
normalized by its maximum, along the boundary (on the circle of radius
$R_0=1$), near the localization point $s = 0$.  One can see that the
WKB approximation, $\exp(-(\theta_0(s) + h^{\frac 23}\theta_1(s))/h)$,
obtained with $\theta_0(s)$ and $\theta_1(s)$ given by
\eqref{eq:theta0_disk} and \eqref{eq:theta1_disk}, accurately captures
the behavior over the range of $s$ between $-0.3$ and $0.3\,$.  Note
that its reduced version, $\exp(-\theta_0(s)/h)$, is also accurate.

\section{Application to diffusion NMR}
\label{sec:NMR}

In this section, we briefly discuss (with no pretention to
mathematical rigor) a possible application of the proposed spectral
analysis of the Bloch-Torrey operator to diffusion NMR
\cite{Grebenkov07}.  In this field, the BT-operator governs the
evolution of the transverse nuclear magnetization which satisfies the
Bloch-Torrey equation
\begin{equation}
\frac{\partial}{\partial t} m(x,t) = \bigl[D \Delta - i \gamma g x_1\bigr] m(x,t) ,
\end{equation}
subject to the uniform initial condition $m(x,0) = 1$.  Here $D$ is
the diffusion coefficient, $g$ the magnetic field gradient, $\gamma$
the gyromagnetic ratio, and the gradient is considered to be constant
in time.  For a bounded domain, the long-time asymptotic behavior of
the solution is determined by the first eigenvalue $\lambda^{(1)}$ of
the BT-operator (with the smallest real part):
\begin{equation}
m(x,t) \simeq C u^{(1)}(x) \exp(-\omega t)  \quad (t\to\infty),
\end{equation}
where 
\begin{equation}
\omega = \gamma g \lambda^{(1)}_h  , \qquad h^2 = D/(\gamma g).  
\end{equation}
Admitting\footnote{This has not be proven mathematically.}  that the
formal asymptotic expansion (\ref{eq:lambda2}) with $n=k=1$ is the
asymptotics of the eigenvalue $\lambda^{(1)}_h$ with the smallest real
part, we obtain in the limit of large $g$
\begin{equation}
\begin{split}
\omega & = i \, \gamma g v_{00} + D^{\frac 13} (\gamma g)^{\frac 23} \mu_0^\# |v_{01}|^{\frac23}  \exp\left(\frac{i\pi}{3} \sign \, v_{01}\right) \\
& + D^{\frac 12} (\gamma g)^{\frac 12} |v_{20}|^{\frac12} \exp\left(\frac{i\pi}{4} \sign v_{20}\right) 
+ D^{\frac 23} (\gamma g)^{\frac 13} \lambda_4^{\#,(1)} + \mathcal O(g^\frac 16)\,,  \\
\end{split}
\end{equation}
where the coefficients $v_{jk}$ are defined by the local
parameterization $V(x) = x_1$ of the boundary near a point from
$\Omega_\perp$.  The real part of $\omega$ determines the decay rate
of the transverse magnetization and the related macroscopic signal.

The leading term of order $(\gamma g)^{\frac 23}$ was predicted for
impermeable one-dimensional domains (with Neumann boundary condition)
by Stoller {\it et al.} \cite{Stoller91} and experimentally confirmed
by H\"urlimann {\it et al.} \cite{Hurlimann95}.  The next-order
correction was obtained by de Swiet and Sen \cite{deSwiet94} for an
impermeable disk.  In the present paper, we generalized these results
to arbitrary planar domains with smooth boundary and to various
boundary conditions (Neumann, Dirichlet, Robin, transmission) and
provided a general technique for getting higher-order corrections (in
particular, we derived the last term).  Moreover, we argued (without
rigorous proof) that these asymptotic relations should also hold for
unbounded domains.

\appendix
\section{Explicit computation of $\lambda_4$}
\label{sec:lambda4}

This Appendix presents the explicit computation of the coefficient
$\lambda_4$ in front of the $h^{4/3}$ term of the four-term
asymptotics (\ref{eq:lambda_asympt}).  Although this is not the
leading term, it is sensitive to the type of boundary condition.  This
is particularly clear for the physically relevant case when the
parameter $\kappa$ of the Robin or transmission boundary condition
scales as $h^{2/3}$.  In this case, the boundary condition for the
rescaled problem is getting closer and closer to the Neumann one, and
the information about the boundary properties appears only in the
$h^{4/3}$ term (e.g., compare Eqs. (\ref{eq:lambda_app2Th}) and
(\ref{eq:lambda_app2Th2})).  The related information on the membrane
permeability or the surface relaxivity of a sample can potentially be
extracted from diffusion NMR experiments.

\subsection{Evaluation of the integral with $\phi_1$}

In order to compute $\lambda_4$ from (\ref{eq:lambda4_def}), we first
evaluate the integral
\begin{equation}
\label{eq:eta_auxil}
\eta = \int\limits_{-\infty}^\infty \sigma \, \phi_1(\sigma) \, \phi_0(\sigma) \, d\sigma .
\end{equation}
We recall that $\phi_1(\sigma)$ satisfies 
\begin{equation}\label{a18bis}
  (\L_2-\lambda_2)\, \phi_1 =  c_{11} \, \sigma\, \phi_0 \,,
\end{equation}
with
\begin{equation}
c_{11} := -  i \, v_{11} \, \int  \tau  \psi_0^\#(\tau)^2 d\tau  \, .
\end{equation}

As a solution of (\ref{a18bis}), we search for some eigenpair
$\{\lambda_2, \phi_0\} = \{ \lambda_2^{(k)}, \phi_0^{(k)}\}$, with
some fixed $k\geq 1$, where $\lambda_2^{(k)}$ and $\phi_0^{(k)}$ are
the eigenvalues and eigenfunctions of the quantum harmonic oscillator
given explicitly in (\ref{eq:phi0_k}).  Since $\phi_0^{(k)}$ are
expressed through the Hermite polynomials $H_k$, one can use their
recurrence relation, $H_{k+1}(x) = 2x H_k(x) - 2 k H_{k-1}(x)$, to
express
\begin{equation}
\sigma \, \phi_0^{(k)} = \frac{\sqrt{k} \, \phi_0^{(k+1)} + \sqrt{k-1} \, \phi_0^{(k-1)}}{(2\gamma)^{\frac12}}  .
\end{equation}
It is therefore natural to search for the solution of (\ref{a18bis})
in the form
\begin{equation}
\phi_1(\sigma) = C_1 \, \phi_0^{(k+1)}(\sigma) + C_2 \, \phi_0^{(k-1)}(\sigma) \,.
\end{equation}
The coefficients $C_1$ and $C_2$ are determined by substituting this
expression into (\ref{a18bis}):
\begin{equation}
\begin{split}
  (\L_2-\lambda_2)\, \phi_1 & = C_1 \bigl(\lambda_2^{(k+1)} - \lambda_2^{(k)}\bigr) \phi_0^{(k+1)} 
+ C_2 \bigl(\lambda_2^{(k-1)} - \lambda_2^{(k)}\bigr) \phi_0^{(k-1)} \\
& = c_{11} \, \frac{\sqrt{k} \, \phi_0^{(k+1)} + \sqrt{k-1} \, \phi_0^{(k-1)}}{(2\gamma)^{\frac12}}  , \\
\end{split}
\end{equation}
from which $C_1 = c_{11} \sqrt{k}/(2\gamma)^{\frac 32}$ and $C_2 = -
c_{11} \sqrt{k-1}/(2\gamma)^{\frac 32}$, where we used
$\lambda_2^{(k)} = \gamma (2k-1)$, with $\gamma = |v_{20}|^{\frac12}
\exp\left(\frac{\pi i}{4}\, \sign v_{20}\right)$.  We get then
\begin{equation}
\phi_1(\sigma) = \frac{c_{11}}{(2\gamma)^{\frac 32}} \biggl( \sqrt{k}\, \phi_0^{(k+1)}(\sigma) - \sqrt{k-1} \, \phi_0^{(k-1)}(\sigma) \biggr) \, .
\end{equation} 
Substituting this expression into (\ref{eq:eta_auxil}), one gets
\begin{equation}
\label{eq:eta}
\eta = \frac{c_{11}}{4\gamma^2} = -  \frac{v_{11}}{4 v_{20}} \, \int  \tau  \psi_0^\#(\tau)^2 d\tau  \, ,
\end{equation}
independently of $n$.  We conclude from (\ref{eq:lambda4_def}) that
\begin{equation}
\label{eq:lambda4_aux1}
\lambda_4^\# = -i \frac{v_{11}^2 [I_1^\#]^2}{4v_{20}} + \frac{\curv(0)}{2} \int \partial_\tau [\psi_0^\#(\tau)]^2 + i v_{02} I_2^\# ,
\end{equation}
where
\begin{equation}
I_1^\# = \int \tau \, \psi_0^\#(\tau)^2 \, d\tau\, , \qquad  I_2^\# = \int \tau^2 \, \psi_0^\#(\tau)^2 \, d\tau\, . 
\end{equation}

\subsection{Evaluation of the integrals with $\psi_0^\#$}

In order to compute these integrals, we consider the function $\Psi(x)
= \Ai(\alpha + \beta x)$ that satisfies the Airy equation
\begin{equation}
(-\partial_x^2 + \beta^3 x + \beta^2 \alpha) \Psi(x) = 0\, .
\end{equation}
Multiplying this equation by $\Psi'(x)$, $\Psi(x)$, $x \Psi'(x)$, $x
\Psi(x)$, or $x^2 \Psi'(x)$ and integrating from $0$ to infinity, one
gets the following five relations:
\begin{enumerate}
\item  
\begin{equation*}
- \int\limits_0^\infty \Psi'' (x) \Psi' (x)\, dx   + \int\limits_0^\infty  ( \beta^3 x  + \beta^2 \alpha) \Psi(x) \Psi' (x)\, dx  = 0\, ,
\end{equation*}
which leads to the determination of $\int_0^{+\infty} \Psi (x)^2 dx$
by the formula
\begin{equation}
 \Psi'(0)^2   -  \beta^2 \alpha \Psi(0)^2 -    \beta^3  \int\limits_0^\infty    \Psi(x)^2\, dx  = 0\, .
\end{equation}
\item 
\begin{equation*}
- \int\limits_0^\infty \Psi'' (x) \Psi (x)\, dx   + \int\limits_0^\infty  ( \beta^3 x  + \beta^2 \alpha) \Psi(x)^2 dx  = 0\, .
\end{equation*}
Here we remark that
\begin{equation*}
\int\limits_0^\infty \Psi'' (x) \Psi (x)\, dx = \Psi'(0) \Psi (0) - \int\limits_0^\infty \Psi'(x)^2 \, dx
\end{equation*}
and get
\begin{equation}
 - \Psi'(0) \Psi (0) + \int\limits_0^\infty \Psi'(x)^2 \, dx   + \int\limits_0^\infty  ( \beta^3 x  + \beta^2 \alpha) \Psi(x)^2 dx  = 0\, .
\end{equation}
\item 
\begin{equation*}
- \int\limits_0^\infty \Psi'' (x) x \Psi' (x)\, dx   + \int\limits_0^\infty  ( \beta^3 x  + \beta^2 \alpha)x \Psi(x)\Psi' dx  = 0\,  \quad \Longrightarrow
\end{equation*}
\begin{equation}
\frac 12  \int\limits_0^\infty  \Psi' (x) ^2 \, dx   - \frac 12  \int\limits_0^\infty  (2 \beta^3x   + \beta^2 \alpha) \Psi(x)^2 dx  = 0\, .
\end{equation}
\item 
\begin{equation*}
- \int\limits_0^\infty \Psi'' (x) x \Psi (x)\, dx   + \int\limits_0^\infty  ( \beta^3 x  + \beta^2 \alpha) x \Psi(x)^2 dx  = 0\,  \quad \Longrightarrow
\end{equation*}
\begin{equation*}
 \int\limits_0^\infty \Psi'(x) (x \Psi (x))' \, dx   + \int\limits_0^\infty  ( \beta^3 x  + \beta^2 \alpha) x \Psi(x)^2 dx  = 0\,  \quad \Longrightarrow
\end{equation*}
\begin{equation}
 \int\limits_0^\infty x  \Psi'(x)^2 \, dx  - \frac 12 \Psi (0)^2    + \int\limits_0^\infty  ( \beta^3 x  + \beta^2 \alpha) x \Psi(x)^2 dx  = 0\, .
\end{equation}
\item 
\begin{equation*}
- \int\limits_0^\infty \Psi'' (x) x^2 \Psi' (x)\, dx   + \int\limits_0^\infty  ( \beta^3 x  + \beta^2 \alpha)x^2\Psi(x)\Psi' dx  = 0\,   \quad \Longrightarrow
\end{equation*}
\begin{equation}
 \int\limits_0^\infty x  \Psi'(x)^2 \, dx   + \int\limits_0^\infty  ( \beta^3 x  + \beta^2 \alpha)x^2\Psi(x)\Psi' dx  = 0\, .
\end{equation}
\end{enumerate}
In the above relations, we assume that the parameter $\beta$ is
such that $|\arg(\beta)| < \pi/3$ so that $\Psi(+\infty) =
\Psi'(+\infty) = 0$ (otherwise the integrals could diverge).
So we get a linear system of five equations satisfied by $\int \Psi^2
dx$, $\int x \Psi^2 dx$, $\int x^2 \Psi^2 dx$, $\int \Psi'(x)^2 dx$
and $\int x \Psi'(x)^2 dx$.  Solving this system, we obtain
\begin{eqnarray}
\label{eq:Psi_int}
\int\limits_0^\infty \Psi^2(x) dx &=& \beta^{-3} [\Psi'(0)]^2 - \alpha \beta^{-1} [\Psi(0)]^2  \nonumber\\
&=& \frac{[\Ai'(\alpha)]^2 - \alpha [\Ai(\alpha)]^2}{\beta} , \\
\int\limits_0^\infty x \Psi^2(x) dx &=& \frac{1}{3\beta^3} \biggl(- \Psi'(0) \Psi(0) - 2\alpha \beta^2  \int\limits_0^\infty \Psi^2(x) dx \biggr) \nonumber \\
& =& - \frac{\Ai(\alpha) \Ai'(\alpha) + 2\alpha [\Ai'(\alpha)]^2 - 2\alpha^2 [\Ai(\alpha)]^2}{3\beta^2}\, , \\
\int\limits_0^\infty x^2 \Psi^2(x) dx &=& \frac{1}{5\beta^3} \biggl([\Psi(0)]^2 - 4\alpha \beta^2  \int\limits_0^\infty x~\Psi^2(x) dx \biggr) \nonumber \\
& =& \frac{[\Ai(\alpha)]^2 + \frac43 \alpha \bigl(\Ai(\alpha) \Ai'(\alpha) + 2\alpha [\Ai'(\alpha)]^2 - 2\alpha^2 [\Ai(\alpha)]^2\bigr)}{5\beta^3}\, .
\end{eqnarray}
These relations allow one to compute the normalization constant
$c_n^\#$ of quasimodes and the contribution $\lambda_4^\#$ to the
eigenvalue.  We consider successively Dirichlet, Neumann, Robin, and
Transmission cases.

\subsection*{Dirichlet case} 

The function $\psi_0^D(\tau)$ from (\ref{eq:psi0D}) corresponds to
$\alpha = a_n$ and $\beta = |\, v_{01}|^{\frac13}
\exp\left(\frac{i\pi}{6}\sign \, v_{01}\right)$ so that $\Ai(\alpha) =
0\,$.  The normalization constant $c^D_n$ in (\ref{eq:psi0N}) is then
\begin{equation}
\label{eq:CD}
(c^D_n)^{-2} = \frac{[\Ai'(a_n)]^2}{\beta}\,.
\end{equation}
Using (\ref{eq:Psi_int}), one gets
\begin{eqnarray}
I_1^D & =& \int\limits_0^\infty \tau [\psi_0^D(\tau)]^2 d\tau = - \frac{2a_n}{3\beta} \, ,\\
I_2^D & =& \int\limits_0^\infty \tau^2 [\psi_0^D(\tau)]^2 d\tau = \frac{8 a_n^2}{15 \beta^2} \, . 
\end{eqnarray}
Using (\ref{eq:lambda2}) and (\ref{eq:eta}), we obtain 
\begin{eqnarray}
\label{eq:lambda4D}
\lambda_4^{D,(n)} & =& i \frac{v_{11}^2 a_n^2}{9 v_{20} \beta^2} - \frac{\curv(0)}{2} [\psi_0^D(0)]^2
+ i v_{02} \frac{8 a_n^2}{15 \beta^2}  \nonumber \\
& =& \frac{i a_n^2}{|\, v_{01}|^{\frac23} \exp\left(\frac{i\pi}{3}\sign \, v_{01}\right)} 
\biggl(\frac{1}{9} ~ \frac{v_{11}^2}{v_{20}} +  \frac{8}{15} v_{02} \biggr),  
\end{eqnarray}
where we used $\psi_0^D(0) = 0\,$.

\subsection*{Neumann case}   

The function $\psi_0^N(\tau)$ from (\ref{eq:psi0N}) corresponds to
$\alpha = a'_n$ and $\beta = |\, v_{01}|^{\frac13}
\exp\left(\frac{i\pi}{6}\sign \, v_{01}\right)$ so that $\Ai'(\alpha)
= 0$.  The normalization constant $c^N_n$ in (\ref{eq:psi0N}) is then
\begin{equation}
\label{eq:CN}
(c^N_n)^{-2} = \frac{[\Ai'(\alpha)]^2 - \alpha [\Ai(\alpha)]^2}{\beta} = - \frac{a'_n [\Ai(a'_n)]^2}{\beta}.
\end{equation}
Using (\ref{eq:Psi_int}), one gets
\begin{eqnarray}
I_1^N &=& \int\limits_0^\infty \tau [\psi_0^N(\tau)]^2 d\tau = - \frac{2a'_n}{3\beta} \, ,  \\   
I_2^N &=& \int\limits_0^\infty \tau^2 [\psi_0^N(\tau)]^2 d\tau  = \frac{8(a'_n)^3 - 3}{15 a'_n \beta^2} \, , 
\end{eqnarray}
from which
\begin{equation}
\label{eq:lambda4N}
\lambda_4^{N,(n)} = \frac{i}{|\, v_{01}|^{\frac23} \exp\left(\frac{i\pi}{3}\sign \, v_{01}\right)} 
\biggl( - \frac{(a'_n)^2}{ 18} ~ \frac{v_{11}^2}{v_{20}} + \frac{1}{2a'_n} \curv(0) \, v_{01} 
 + \frac{8(a_n')^3 - 3}{15 a'_n} v_{02} \biggr). 
\end{equation}

\subsection*{Robin case}

The function $\psi_0^R(\tau)$ from (\ref{eq:psi0R}) corresponds to
$\beta = |\, v_{01}|^{\frac13} \delta$ and $\alpha = a_n^R(\kappa)$ so
that $\Ai'(\alpha) = \hat\kappa \, \Ai(\alpha)\,$, with $\hat\kappa =
\kappa/(\delta \, |v_{01}|^{\frac13})$ and $\delta =
\exp\left(\frac{i\pi}{6}\, \sign \, v_{01}\right)\,$.  The normalization
constant $c^R_n$ in (\ref{eq:psi0R}) is then
\begin{equation}
\label{eq:CR}
(c^R_n)^{-2} = \frac{[\Ai(a_n^R(\kappa))]^2}{\beta} \bigl[{\hat\kappa}^2 - a_n^R(\kappa)\bigr] 
= [\Ai(a_n^R(\kappa))]^2 ~ \frac{\kappa^2 + \lambda_0^R}{iv_{01}}  \, ,
\end{equation}
where we used (\ref{eq:lambda0R}) for $\lambda_0^R$.

Using (\ref{eq:Psi_int}), one gets
\begin{eqnarray}
I_1^R = \int\limits_0^\infty \tau [\psi_0^R(\tau)]^2 d\tau 
&=& - \frac{\hat\kappa + 2\hat\kappa^2 a_n^R(\kappa) - 2[a_n^R(\kappa)]^2}{3\beta [\hat\kappa^2 - a_n^R(\kappa)]}  \nonumber \\
&=& \frac{2\lambda_0^R}{3iv_{01}} - \frac{\kappa}{3(\kappa^2 + \lambda_0^R)} \, , \\
I_2^R = \int\limits_0^\infty \tau^2 [\psi_0^R(\tau)]^2 d\tau 
&=& \frac{1 + \frac43 a_n^R(\kappa) \bigl[\hat\kappa + 2\hat\kappa^2 a_n^R(\kappa) - 2[a_n^R(\kappa)]^2\bigr]}{5 \beta^2 [\hat\kappa^2 - a_n^R(\kappa)]} \nonumber \\
&=& \frac{1}{5(\kappa^2 + \lambda_0^R)} - \frac{8 [\lambda_0^R]^2}{15\, v_{01}^2} - \frac{4 \kappa \lambda_0^R}{15\, iv_{01} (\kappa^2 + \lambda_0^R)} \, . 
\end{eqnarray}
Using (\ref{eq:lambda2}) and  (\ref{eq:eta}), we obtain
\begin{eqnarray}
\label{eq:lambda4R}
\lambda_4^{R,(n)} &=& - i\frac{v_{11}^2 [I_1^R]^2}{4v_{20}} - \frac{\curv(0)}{2} [\psi_0^R(0)]^2 + i v_{02} I_2^R  \nonumber \\
&=& -i \frac{v_{11}^2 [I_1^R]^2}{4v_{20}} - \frac{\curv(0)}{2} \frac{iv_{01}}{\kappa^2 + \lambda_0^R} + i v_{02} I_2^R \,.
\end{eqnarray}

\begin{remark}\label{RemRkappa=0}
It is clear from the computation that $\lambda_4^{R,(n)}$ belongs to
$C^\infty$ in a neighborhood of $0\,$.  In particular, we recover
\begin{equation}
\lambda_4^{R,(n)}(0) =\lambda_4^{N,(n)}\,.
\end{equation}
\end{remark}

\subsection*{Transmission case}

In order to compute the above integrals for the transmission case, we
note that (\ref{eq:an_equation}) can be written as
\begin{equation}
\Ai'(a_n^+) \, \Ai'(a_n^-) = - \frac{\kappa}{2\pi |\, v_{01}|^{\frac13}} \,,
\end{equation}
while the Wronskian for Airy functions yields another relation:
\begin{equation}
\bar \delta \Ai'(a_n^-) \Ai(a_n^+) + \delta \Ai'(a_n^+) \Ai(a_n^-) = - \frac{1}{2\pi}\, ,
\end{equation}
where $\delta = \exp\bigl(\frac{\pi i}{6}\, \sign \, v_{01}\bigr)$,
and $a_n^\pm = a_n^\pm(\kappa)$ are given by (\ref{eq:ankappa}).\\  
From (\ref{eq:psi0T}), we then obtain
\begin{equation}
(c^T_n)^{-2} = \frac{a_n^+ \bar\delta}{2\pi |\, v_{01}|^{\frac13}} \biggl(\bar \delta \Ai'(a_n^-) \Ai(a_n^+) - \delta \Ai'(a_n^+) \Ai(a_n^-)\biggr)\,.
\end{equation}

Using (\ref{eq:Psi_int}), we get
\begin{eqnarray}
I_1^T &=& \int\limits_{-\infty}^\infty \tau \psi_0^T(\tau)^2 d\tau = \frac{(c_n^T)^2}{3|\, v_{01}|^{\frac23}} 
\biggl(\frac{\kappa \delta^3}{4\pi^2 |\, v_{01}|^{\frac13}}  \nonumber \\
&+& \frac{(a_n^+)^2 \delta^4}{\pi} \bigl(\bar \delta \Ai'(a_n^-) \Ai(a_n^+) - \delta \Ai'(a_n^+) \Ai(a_n^-)\bigr) \biggr)  \nonumber \\
&=& (c_n^T)^2 \frac{\kappa i}{12\pi^2 \, v_{01}} - \frac{2a_n^+}{3\delta |\, v_{01}|^{\frac13}} \,,   \\
I_2^T &=& \int\limits_{-\infty}^\infty \tau^2 \psi_0^T(\tau)^2 d\tau  = 
\frac{(c_n^T)^2}{5|\, v_{01}|} \biggl(\frac{\kappa a_n^+ \bar\delta^4}{3\pi^2 |\, v_{01}|^{\frac13}}  \nonumber \\
&-& \frac{8[a_n^+]^3-3}{6\pi} \delta^3 \bigl(\bar \delta \Ai'(a_n^-) \Ai(a_n^+) - \delta \Ai'(a_n^+) \Ai(a_n^-)\bigr) \biggr) \nonumber \\
&=& (c_n^T)^2\, \frac{\kappa a_n^+ \bar \delta^4 }{15\pi^2 |\, v_{01}|^{\frac43}}
+ \frac{8[a_n^+]^3-3}{15 a_n^+ \delta^2 |\, v_{01}|^{\frac23}} \,.
\end{eqnarray}
Finally, we compute the coefficient in front of $\frac12\curv(0)$ in
(\ref{eq:lambda4_aux1}):
\begin{equation}
I_0^T := \int \partial_\tau [\psi_0^T(\tau)]^2 = [\psi_0^-(0)]^2 - [\psi_0^+(0)]^2 = \frac{|\, v_{01}|^{\frac13} 
\exp\left(\frac{i\pi}{6}\, \sign \, v_{01}\right)}{a_n^+} \,.
\end{equation}
We conclude that
\begin{equation}
\label{eq:lambda4T}
\lambda_4^{T,(n)} = -i \frac{v_{11}^2 [I_1^T]^2}{4 v_{20}} + \curv(0) \frac{|\, v_{01}|^{\frac13} 
\exp\left(\frac{i\pi}{6}\, \sign \, v_{01}\right)}{2 a_n^+} + i \,v_{02} \, I_2^T\, .
\end{equation}

\begin{remark}\label{RemTkappa=0}
It is clear from the computation that $\lambda_4^{T,(n)}(\kappa)$
belongs to $C^\infty$ in a neighborhood of $0$.  In particular, we
recover
\begin{equation}
\lambda_4^{T,(n)}(0) =\lambda_4^{N,(n)}\,.
\end{equation}
\end{remark}

\subsection{Evaluation of the derivative $(\mu^T_n)'(0)$}   \label{sec:muT}

The asymptotic relation (\ref{newscaling}) involves the derivative of
$\mu^\#_n(\kappa)$ with respect to $\kappa$ at $\kappa = 0$.  In this
subsection, we provide its explicit computation for the transmission
case.  According to (\ref{eq:ankappa}), we have
\begin{equation}
\mu^T_n(\kappa) = - a^+_n(\kappa) = - \lambda^T_n(\kappa/|v_{01}|^{\frac13})\, \exp\left(\frac{2\pi i}{3}\, \sign v_{01}\right)\,,
\end{equation}
where $\lambda^T_n$ satisfies (\ref{eq:an_equation}).\\
  The
derivative with respect to $\kappa$ at $\kappa = 0$ reads
\begin{equation}
(\mu^T_n)'(0) = \left(\frac{\partial}{\partial\kappa} \mu^T_n(\kappa)\right)_{\kappa = 0} 
= - (\lambda^T_n)'(0)\, \frac{1}{|v_{01}|^{\frac13}} \, \exp\left(\frac{2\pi i}{3}\, \sign v_{01}\right) \,.
\end{equation}
In turn, $(\lambda^T_n)'(0)$ can be obtained by differentiating
 (\ref{eq:an_equation}) with respect to $\kappa$
\begin{equation}
\label{eq:Aauxil12}
\begin{split}
2\pi \frac{(\lambda^T_n)'(0)}{|v_{01}|^{\frac13}} & \biggl[ e^{-i\alpha} \, \lambda^T_n(0) \, \Ai'(e^{-i\alpha} \lambda^T_n(0))\, \Ai(e^{i\alpha}\lambda^T_n(0)) \\
& + e^{i\alpha} \, \lambda^T_n(0) \, \Ai'(e^{i\alpha} \lambda^T_n(0)) \, \Ai(e^{-i\alpha}\lambda^T_n(0)) \biggr] = - \frac{1}{|v_{01}|^{\frac13}}\, , \\
\end{split}
\end{equation}
where we used the Airy equation: $\Ai''(z) = z\Ai(z)$, and a shortcut
notation $\alpha = 2\pi/3\,$.\\
At $\kappa = 0\,$, (\ref{eq:an_equation}) admits two solutions,
$\lambda^T_n(0) = e^{i \alpha} \, a'_n$ and $\lambda^T_n(0) = e^{-i
\alpha} \, a'_n\,$, that correspond to $v_{01} < 0$ and $v_{01} >
0\,$, respectively.

When $v_{01} < 0$, the first term in (\ref{eq:Aauxil12}) vanishes (as
$\Ai'(e^{-i\alpha} \lambda^T_n(0)) = 0\, $), while the second term can
be expressed by using the Wronskian,
\begin{equation}
e^{-i\alpha} \Ai'(e^{-i\alpha} z) \Ai(e^{i\alpha} z) - e^{i\alpha} \Ai'(e^{i\alpha} z) \Ai(e^{-i\alpha} z) = \frac{i}{2\pi}  \qquad \forall~ z\in\C\,. 
\end{equation}
We get then 
\begin{equation*}
(\lambda^T_n)'(0) = \frac{i}{\lambda^T_n(0)} = \frac{i}{a'_n e^{i\alpha}} \,.
\end{equation*}
In turn, when $v_{01} > 0$, the second term in (\ref{eq:Aauxil12})
vanishes, while the first term yields 
\begin{equation*}
(\lambda^T_n)'(0) = \frac{-i}{\lambda^T_n(0)} = \frac{-i}{a'_n e^{-i\alpha}} \,.
\end{equation*}  
Combining these relations, we obtain
\begin{equation}
(\mu^T_n)'(0) = - \frac{1}{a'_n \, |v_{01}|^{\frac13}} \, \exp\left(-\frac{\pi i}{6}\, \sign v_{01}\right).
\end{equation}

\end{document}